\documentclass[12pt,reqno]{amsart}
\usepackage{amssymb}
\usepackage{mathtools}
\usepackage{txfonts}
\usepackage{eucal}
\usepackage{extarrows}

\usepackage[all]{xy}

\usepackage{tikz}

\usepackage{graphicx}
\usepackage{float}

\usepackage{xspace}

\usepackage{amsmath}
\usepackage{xcolor}
\usepackage{amsmath,amssymb}
\usepackage{tikz}

\usepackage{amssymb}
\usepackage{multirow} 
\usepackage{array}    

\usepackage[a4paper,body={16.3cm,22.8cm},centering]{geometry}
\usepackage[colorlinks,final,backref=page,hyperindex,pdftex]{hyperref}

\usepackage{tikz-cd}

\newcommand{\nc}{\newcommand}
\newcommand{\delete}[1]{}

	\nc{\mlabel}[1]{\label{#1}}  
	\nc{\mcite}[1]{\cite{#1}}  
	\nc{\mref}[1]{\ref{#1}}  
	\nc{\mbibitem}[1]{\bibitem{#1}} 

\delete{
	\nc{\mlabel}[1]{\label{#1}  
		{\hfill \hspace{1cm}{\small\tt{{\ }\hfill(#1)}}}}
	\nc{\mcite}[1]{\cite{#1}{\small{\tt{{\ }(#1)}}}}  
	\nc{\mref}[1]{\ref{#1}{{\tt{{\ }(#1)}}}}  
	\nc{\mbibitem}[1]{\bibitem[\bf #1]{#1}} 
}
\newtheorem{theorem}{Theorem}[section]
\newtheorem{prop}[theorem]{Proposition}
\newtheorem{lemma}[theorem]{Lemma}
\newtheorem{coro}[theorem]{Corollary}

\theoremstyle{definition}
\newtheorem{defn}[theorem]{Definition}
\newtheorem{remark}[theorem]{Remark}
\newtheorem{exam}[theorem]{Example}
\newtheorem{prop-def}{Proposition-Definition}[section]

\newcommand\cal[1]{\mathcal{#1}}

\newcommand\alphlist{a,b,c,d,e,f,g,h,i,j,k,l,m,n,o,p,q,r,s,t,u,v,w,x,y,z}
\newcommand\Alphlist{A,B,C,D,E,F,G,H,I,J,K,L,M,N,O,P,Q,R,S,T,U,V,W,X,Y,Z}
\newcommand\getcmds[3]{\expandafter\newcommand\csname #2#1\endcsname{#3{#1}}}
\makeatletter
\@for\x:=\alphlist\do{\expandafter\getcmds\expandafter{\x}{frak}{\mathfrak}}
\@for\x:=\Alphlist\do{\expandafter\getcmds\expandafter{\x}{frak}{\mathfrak}}
\makeatother

\nc{\bfk}{{\bf k}}
\font\cyr=wncyr10

\newfont{\scyr}{wncyr10 scaled 550}
\nc{\sha}{\mbox{\cyr X}}
\nc{\ssha}{\mbox{\bf \scyr X}}
\nc{\Id}{\mathrm{Id}}
\nc{\lbar}[1]{\overline{#1}}
\nc{\ot}{\otimes}
\nc{\dep}{\mathrm{dep}}

\nc{\tred}[1]{\textcolor{red}{#1}} \nc{\tgreen}[1]{\textcolor{green}{#1}}
\nc{\tblue}[1]{\textcolor{blue}{#1}} \nc{\tpurple}[1]{\textcolor{purple}{#1}}

\nc{\li}[1]{\tpurple{\underline{Li:}#1 }}
\nc{\liadd}[1]{\tpurple{#1}}
\nc{\xing}[1]{\tblue{\underline{Xing:}#1 }}
\nc{\dominique}[1]{\tblue{\underline{Dominique: }#1 }}
\nc{\yuan}[1]{\tred{\underline{Yuan:}#1 }}
\nc{\markus}[1]{\tred{\underline{Markus:} #1}}
\nc\hu[1]{\tgreen{\underline{Huhu:}#1}}

\usetikzlibrary{calc}
\newlength\xch
\newsavebox\dbox
\sbox\dbox{\tikz{\fill (0,0) circle (0.05cm);}}
\newif\ifqdd
\newif\ifzdd
\nc{\dnx}{\Delta_n A} \nc{\dx}{\Delta A} \nc{\dgp}{{\rm deg_{P}}}
\nc{\dgt}{{\rm deg_{T}}} \nc{\dg}{{\rm deg}} \nc{\ida}{ID($A$)} \nc{\tu}{\tilde{u}} \nc{\tv}{\tilde{v}}
\nc{\nr}{\calr_n} \nc{\nz}{\calz_n} \nc{\fun}{\cala_{n,d}}
\nc{\fbase}{\calb} \nc{\LF}{\mathrm{RF}} \nc{\FFA}{\mathrm{LF}} \nc{\irr}{\mathrm{Irr}}
\nc{\result}{\bfk\mathrm{Irr}(S_n)}  \nc{\I}{I_{\mathrm{ID},n}^0}
\nc{\nrs}{\calr_n^\star} \nc{\ii}{\mathrm{I}} \nc{\iii}{\mathrm{II}}
\nc{\intl}{{\rm int}}\nc{\ws}[1]{{#1}}\nc{\deleted}[1]{\delete{#1}}\nc{\plas}{placements\xspace}

\nc{\bim}[1]{#1}  \nc{\shaop}{\sha_{\Omega}^{+}}  \nc{\shao}{\sha_{\Omega}}
\nc{\bbim}[2]{#1 #2} \nc{\bbbim}[2]{#1,\, #2} \nc{\RBF}{{\rm RBF}}
\nc{\frb}{F_{\RB}} \nc{\shaf}{\ssha_{\tiny{\Omega}}} \nc{\sham}{\diamond_{\tiny{\Omega}}}
\nc{\lf}{\lfloor} \nc{\rf}{\rfloor} \nc{\shan}{\ssha_{\lambda}}
\nc{\rlex}{{\rm {lex}}} \nc{\bb}{\Box} \nc{\ra}{\rightarrow}
\nc{\e}{{\rm {e}}}
\nc{\DDF}{\mathrm{DD}(X,\,\Omega)}\nc{\DTF}{\mathrm{DT}(X,\,\Omega)} \nc{\DT}{\mathrm{DT}'(\Omega,\,V)}
\nc{\bra}{\mathrm{bra}} \nc{\bre}{\mathrm{bre}}
\nc{\dec}{\mathrm{dec}} \nc{\diamondw}{\diamond_{w}}
\nc{\type}{\mathrm{type}}

\nc\caF[1]{\cal{F}_{#1}(X,\,\Omega)}
\nc\calt{\cal{T}(X,\,\Omega)} \nc\caltn{\cal{T}_n(X,\,\Omega)}
\nc\caltbin{\cal{T}_b(X,\,\Omega)}
\nc\calta{\cal{T}_0(X,\,\Omega)}
\nc\caltb{\cal{T}_1(X,\,\Omega)}
\nc\caltc{\cal{T}_2(X,\,\Omega)}
\nc\caltd{\cal{T}_3(X,\,\Omega)}
\nc\caltm{\cal{T}_m(X,\,\Omega)}
\nc\calf{\cal{F}(X,\,\Omega)}
\nc\fram{\frak{M}(\Omega,\, X)}
\nc\shaw{\sha^{NC}_w(\Omega,\, X)}
\nc\dw{\diamond_w} \nc\dl{\diamond_\ell}
\nc\shal{\sha^{NC}_\ell(X,\, \Omega)} \nc\shav{\sha^{NC}_w(\Omega,\, V)} \nc\shat{\sha^{NC,1}_w(\Omega,\, T^{+}(V))}
\nc{\cfo}{\cal{F}(X,\,\Omega)}

\nc{\lar}{\varinjlim}
\nc\XO{(X,\,\Omega)}
\def\cxo#1#2;{\cal{#1}#2\XO}
\def\cxob#1#2;{\cal{#1}#2_b\XO}
\nc\lrf[2]{B_{#2}^+(#1)}
\nc{\fd}{\mathrm{\text{typed angularly decorated planar rooted trees}}}
\nc{\rb}{\mathrm{RBFWs}} \nc{\dfw}{\mathrm{DFW{(X)}}} \nc{\tfw}{\mathrm{TFW{(X)}}}
\nc{\tfv}{\mathrm{TFW{(V)}}} \nc{\rbf}{\mathrm{RBF}}
\nc{\db}{\mathrm{db}}
\nc{\st}{\mathrm{st}}

\def\Ve#1,#2,#3;{\vee_{#1,\,(#2,\,#3)}}
\def\bigv#1;#2;#3;{\bigvee\nolimits_{#1}^{#2;\,#3}}

\nc{\Irr}{\mathrm{Irr}} \nc{\lc}{\lfloor} \nc{\rc}{\rfloor}
\nc{\rswx}{\frak{M}( \Omega_R\sqcup \Omega_S, X)}
\nc{\rswxs}{\frak{M}^\star( \Omega_R\sqcup \Omega_S, X)}
\nc{\Dl}{\leq_{_{{\rm Dl}}}} \nc{\Dll}{<_{_{{\rm Dl}}}} \nc{\bbs}{\mathbb{S}}
\nc{\orbsa}{$\Omega$-Rota-Baxter system\xspace}
\nc{\orbsas}{$\Omega$-Rota-Baxter systems\xspace}
\nc{\mrbs}{matching Rota-Baxter system\xspace}
\nc{\mrbss}{matching Rota-Baxter systems\xspace}

\nc\prbsla[4]{{R}_{#1}\left(#3\right)R_{#2}\left(#4\right)}
\nc\prbsra[4]{{R}_{#1\rightarrow#2}\left(R_{#1\rhd#2}\left(#3\right)#4\right)+{R}_{#1\leftarrow#2}\left(#3S_{#1\lhd#2}\left(#4\right)\right)}
\nc\prbslb[4]{{S}_{#1}\left(#3\right)S_{#2}\left(#4\right)}
\nc\prbsrb[4]{{S}_{#1\rightarrow#2}\left(R_{#1\rhd#2}\left(#3\right)#4\right)+{S}_{#1\leftarrow#2}\left(#3S_{#1\lhd#2}\left(#4\right)\right)}

\nc\rbsla[4]{\lc #3 \rc ^{R}_{#1} \lc #4 \rc ^R_{#2}}
\nc\rbslb[4]{\lc #3\rc ^{S}_{#1}  \lc #4 \rc ^S_{#2}}
\nc\rbsray[4]{\lc \lc #3 \rc^R_{#1\rhd#2} #4 \rc ^{R}_{#1\rightarrow#2}}
\nc\rbsraz[4]{\lc #3 \lc #4\rc^S_{#1\lhd#2}\rc ^{R}_{#1\leftarrow#2}}
\nc\rbsrby[4]{\lc \lc #3\rc ^R_{#1\rhd#2}#4\rc ^{S}_{#1\rightarrow#2}}
\nc\rbsrbz[4]{\lc #3 \lc #4 \rc ^S_{#1\lhd#2}\rc ^{S}_{#1\leftarrow#2}}
\nc\rbsrac[4]{\lc #3 \lc #4\rc^R_{#1\lhd#2}\rc ^{R}_{#1\leftarrow#2}}

\nc\rbslq[4]{\lc #3 \rc ^{Q}_{#1} \lc #4 \rc ^Q_{#2}}
\nc\rbslt[4]{\lc #3\rc ^{T}_{#1}  \lc #4 \rc ^T_{#2}}
\nc\rbsrqy[4]{\lc \lc #3 \rc^R_{#1\rhd#2} #4 \rc ^{Q}_{#1\rightarrow#2}}
\nc\rbsrqz[4]{\lc #3 \lc #4\rc^S_{#1\lhd#2}\rc ^{Q}_{#1\leftarrow#2}}
\nc\rbsrty[4]{\lc \lc #3\rc ^R_{#1\rhd#2}#4\rc ^{T}_{#1\rightarrow#2}}
\nc\rbsrtz[4]{\lc #3 \lc #4 \rc ^S_{#1\lhd#2}\rc ^{T}_{#1\leftarrow#2}}

\nc{\obr}[1]{\lc #1 \rc_\omega^R} \nc{\obs}[1]{\lc #1 \rc_\omega^S} \nc{\obq}[1]{\lc #1 \rc_\omega^*}
\nc{\obqa}[1]{\lc #1 \rc_\alpha^*} \nc{\obqb}[1]{\lc #1 \rc_\beta^*}
\nc{\obra}[1]{\lc #1 \rc_\alpha^R} \nc{\obrb}[1]{\lc #1 \rc_\beta^R}
\nc{\obsa}[1]{\lc #1 \rc_\alpha^S} \nc{\obsb}[1]{\lc #1 \rc_\beta^S}

\nc{\ad}{\mathrm{ad}}
\nc{\id}{\mathrm{id}}

\begin{document}

\title[Nijenhuis BiHom-Lie bialgebras and differential Lie bialgebras]{Nijenhuis BiHom-Lie bialgebras and differential Lie bialgebras}

\author{Jiaqi Liu
}
\address{School of Mathematics and Statistics, Lanzhou University, Lanzhou, 730000, P.\,R. China}
\email{jiaqiliu2025@lzu.edu.cn}

\author{Lin Gao
}
\address{School of Mathematics and Statistics, Henan University, Henan, Kaifeng 475004, P.\,R. China}
\email{gaolin821000@163.com
}

\author{Yuanyuan Zhang$^*$
}
\footnotetext{*Corresponding author.}
\address{School of Mathematics and Statistics, Henan University, Henan, Kaifeng 475004, P.\,R. China}
\email{zhangyy17@henu.edu.cn
}

\date{\today}

\begin{abstract}
In this paper, we first introduce the concept of Nijenhuis BiHom-Lie algebras. We then establish the equivalence relations among Manin triples of Nijenhuis BiHom-Lie algebras, Nijenhuis BiHom-Lie bialgebras, and matched pairs of Nijenhuis BiHom-Lie algebras. Furthermore, we establish the connection between (Nijenhuis) BiHom-Lie bialgebras and the corresponding Yang-Baxter equations. Moreover, we apply the above theoretical framework to differential Lie algebras and obtain a series of related results. Finally, we illustrate the relationships among all structures investigated in this paper via two commutative diagrams.
\end{abstract}

\makeatletter
\@namedef{subjclassname@2020}{\textup{2020} Mathematics Subject Classification}
\makeatother
\subjclass[2020]{
	12H05, 
	16T25, 
    17A30, 
    17B62,  
    17D99. 
}

\keywords{Nijenhuis BiHom-Lie bialgebras, differential Lie bialgebras, Manin triple, Matched pair, classical Yang-Baxter equations}

\maketitle

\tableofcontents

\setcounter{section}{0}

\allowdisplaybreaks

\section{Introduction}
\subsection{BiHom-Lie bialgebras}
Lie algebras, as a class of non-associative algebraic structures, have been widely applied in diverse fields since the late 19th century~\cite{Liealg}. With the growing depth of research on Lie algebras, scholars have not only explored the various fundamental properties of these algebraic objects, but also extended the original framework to develop generalized and operator-versions on this foundation, such as Nijenhuis Lie algebras~\cite[Definition~3.6]{Coho-Nij-Lie}, differential Lie algebras, BiHom-Lie algebras~\cite{gmmp} (see Definition~\ref{defn-BHL}), and Lie bialgebras.

BiHom-algebras were first introduced in \cite{gmmp} and constitute a natural extension of Hom-algebras. Up to now, numerous BiHom-type algebraic structures have been proposed, including BiHom-associative algebras \cite{gmmp}, BiHom-Lie algebras \cite{gmmp}, BiHom-(tri)dendriform algebras \cite{bhtriden}, BiHom-pre-Lie algebras \cite{bhpL}, and BiHom-PostLie algebras \cite{bhPL}.
 BiHom Lie algebras weaken the antisymmetry and Jacobi identity conditions of Lie algebras by introducing two commuting linear maps $\alpha$ and $\beta$. They not only naturally incorporate classical structures such as Hom Lie algebras~\cite{homlie} (when $\alpha$=$\beta$) and Lie algebras (when $\alpha$=$\beta$=id) as special cases, but also demonstrated unique value in the research of the deformation theory of quantum groups, discrete integrable systems and non-commutative geometry.

In recent years, the representation theory of BiHom-Lie algebras has also attracted much attention. Unlike the representation of Lie algebras, the representation of BiHom-Lie algebras must be compatible with structure maps $\alpha~\text{and}~\beta$. In~\cite{BHL}, Cheng and Qi first introduced the fundamental definition of representations of BiHom-Lie algebras.

To date, many bialgebraic structures have been proposed, including Lie bialgebras~\cite{liebialg}, pre-Lie bialgebras~\cite{preLie}, Nijenhuis pre-Lie bialgebra~\cite{Nijprbialg}, Rota-Baxter bialgebras~\cite{RBbialg}, Hom-Lie bialgebras~\cite{Hom-Liebialg} and so on. In this paper, we introduce the concept of BiHom-Lie bialgebras, which simultaneously satisfy the compatibility conditions between the algebraic structure (BiHom-Lie algebras) and the coalgebra structure (BiHom-Lie coalgebras).

\subsection{Operators on (generalized) Lie algebras}
The notion of Nijenhuis operators originates from complex geometry in the study of complex structures~\cite{Nij}, where Nijenhuis~\cite{Nijoperator} introduced the notion of a Nijenhuis tensor. Within Lie algebras, a triple $(L,[-,-], N)$ is a Nijenhuis Lie algebra~\cite{NijLie-1,NijLie-2}, if $(L,[-,-])$ is a Lie algebra, and $N$ is a linear operator on $L$ satisfying
\begin{align}\label{Nij-identity}
	[N(x),N(y)]=N\big([N(x),y]+[x,N(y)]-N[x,y]\big),\qquad \qquad \text{ for any }x,y\in L.
\end{align}
In particular, $(L,[-,-],N)$ reduces to a Lie algebra when taking $N=\id_{L}$. In this framework, Nijenhuis operators serve as powerful tools in describing integrable deformations of algebraic structures and in constructing trivial cohomology classes~\cite{NIj}. Later, there are more and more researchers study the properties of Nijenhuis Lie algebras, for instance, Das~\cite{Coho-Nij-Lie} introduced the cohomologies of Nijenhuis Lie algebras and Nijenhuis Lie bialgebras and recently Song~et al.~\cite{Nij-song} investigated homotopy Nijenhuis structures on $ L_\infty $-algebras from the perspective of operad theory.

When considering the Nijenhuis operators on BiHom-Lie algebras, one must take into account the twisting effect of the two commuting homomorphisms on the Nijenhuis identity. A preliminary definition of BiHom-Nijenhuis operators has appeared in~\cite{BH-Nij-operator}, but the connections with BiHom-Lie algebra cohomology
remain unexplored. This forms the first motivation of the present work: we define the Nijenhuis BiHom-Lie bialgebra by the coboundary operator.

As an algebraic structure about differential operators, differential algebra~\cite{diffalg-1,diffalg-2} is the study of differentiation and nonlinear differential equations by purely algebraic means, without using an underlying topology. So far, the algebraic structures about differential operators that have been proposed include differential Lie algebras~\cite{diffLie}, differential Rota-Baxter algebras~\cite{diff-RBalg}, differential ($q$-tri)dendriform algebras~\cite{diff-dendri} etc. However, there are some fundamental problems remain open, including the interplay between differential Lie bialgebras, matched pairs of differential Lie algebras and manin triples of differential Lie algebras. These gaps provide the second major motivation of this work.

\subsection{Manin triples, matched pair and bialgebras structures}
The concepts of matched pairs and Manin triples form central geometric tools in the theory of Lie bialgebras~\cite{D1983}.
 A matched pair of Lie algebras encodes compatible mutual actions for the construction of a semi-direct product Lie algebra;  Manin triple, namely a Lie algebra with a nondegenerate invariant bilinear form and a decomposition into two complementary isotropic subalgebras, give an equivalent description of Lie bialgebras. Specifically, let's recall the following conclusions:

\begin{theorem}\cite[Theorem~5.2.6]{double-B}\label{th-eq-liebialg}
Let $(L,[-,-])$ and $(L^*,[-,-]_{*})$ be two Lie algebras. Then the following conditions are equivalent.
\begin{enumerate}
\item $(L\oplus L^*, L,L^*)$ is a standard Manin triple of Lie algebras with the bilinear form (\ref{eq-B-d}).
\item $(L,L^*,ad^*,\mathfrak{ad}^{*})$ is a matched pair of Lie algebras.
\item $(L,L^*)$ is a Lie bialgebra.
\end{enumerate}
\end{theorem}

\noindent Extending these concepts to the BiHom-Lie setting poses two challenges:
\begin{enumerate}
 \item the mutual actions in a matched pair must be adjusted to become BiHom-actions, i.e., compatible with the structure maps;
 \item in the case of a BiHom-Manin triple, the bilinear form must satisfy twisted invariance conditions.
 \end{enumerate}
Resolving these two challenges constitutes the third motivation of this paper. Furthermore, we generalize Theorem~\ref{th-eq-liebialg} to the Nijenhuis BiHom framework.

\subsection{Bialgebras structures and classical Yang-Baxter equations}
The Yang-Baxter equations (abbr.~YBEs) first introduced by Baxter et al.~\cite{YBE1,YBE2,YBE3} in statistical mechanics. The classical Yang-Baxter equations (abbr.~CYBEs), as a special case of the YBEs, originates from the quantum inverse scattering method~\cite{CYBE3}. Subsequently, the work in~\cite{CYBE} established the relationship between CYBEs and Lie bialgebras as follows.

\begin{defn}\cite{CYBE}\label{CYBE-bi}
	Let $(L,[-,-])$ be a Lie algebra. 
	For a given element $r= r_{1}\otimes r_{2}\in L\otimes L,$  we define
	\begin{align}\label{eq-delta-r}
		\delta_{r}(x):=(\ad_{x}\otimes \id+\id\otimes \ad_{x})(r)= [x,r_{1}]\otimes r_{2}+r_{1}\otimes [x,r_{2}], \qquad \text{ for any }x\in L.
	\end{align}
	\begin{enumerate}
		\item We call
		\begin{align*}
			\operatorname{CYB}(r):=[r_{12},r_{13}]+[r_{13},r_{23}]+[r_{12},r_{23}]=0
		\end{align*}
		a {\bf classical Yang-Baxter equation~(abbr.~CYBE)} in $(L,[-,-])$, where for $  r_{1}\otimes r_{2}=r=\bar{r}= \bar{r}_{1}\otimes \bar{r}_{2}\in L\otimes L$,
		\begin{align*}
			[r_{12},r_{13}]=&   [r_{1},\bar{r}_{1}]\otimes r_{2}\otimes \bar{r}_{2},\\
			[r_{13},r_{23}]=&  r_{1}\otimes \bar{r}_{1}\otimes [r_{2},\bar{r}_{2}],\\
			[r_{23},r_{12}]=&  r_{1}\otimes [r_{2},\bar{r}_{1}]\otimes \bar{r}_{2}.
		\end{align*}
		\item If $r$ is an antisymmetric solution of the CYBE, then $(L,[-,-],\delta_{r})$ is a Lie bialgebra, where $\delta_{r}$ is defined by Eq.~(\ref{eq-delta-r}). 
	\end{enumerate}
\end{defn}

Since then, numerous researchers have studied the connections between different types of bialgebras and the corresponding Yang-Baxter equations, such as Hom-Lie bialgebras and classical Hom Yang-Baxter equations~\cite{Hom-Liebialg}, Nijenhuis-Lie bialgebras and classical $P$-Nijenhuis Yang-Baxter equations~\cite{NijYBE}
, BiHom-antisymmetric infinitesimal bialgebras and associative BiHom-classical Yang-Baxter equation~\cite{assoBHCYBE}, etc. In this paper, we give the similar results about BiHom-Lie version~(see Theorem~\ref{thm-BHYBE-BHLbi}) and Nijenhuis BiHom-Lie-version~(see Theorem~\ref{thm-BHLbi-NBHLbi}).

$\textbf{The paper is organized as follows.}$ Section~\ref{sec-2} is devoted to the study of Nijenhuis BiHom-Lie algebras and their related structures. We begin by considering BiHom-Lie algebras and BiHom-Lie coalgebras, and based on this, we introduce the notion of BiHom-Lie bialgebras (Definition~\ref{def-BHLbialg}). Subsequently, we define Nijenhuis BiHom-Lie algebras together with their representations, and investigate their fundamental properties. Furthermore, we introduce the concepts of Manin triples of Nijenhuis BiHom-Lie algebras (Definition~\ref{def-Manintriple}), Nijenhuis BiHom-Lie bialgebras (Definition~\ref{def-NBHL bialg}), and matched pairs of Nijenhuis BiHom-Lie algebras (Definition~\ref{MPOBHL}), and prove that these three concepts are equivalent (see Theorems~\ref{th}, \ref{th-23}, \ref{th-123} for details). Finally, we give the concept of $S$-Nijenhuis classical BiHom Yang-Baxter equations~(abbr.~SN-CBHYBEs)~(Definition~\ref{def-Nij-BiHom-YBE}) and prove that this equation induce a Nijenhuis BiHom-Lie bialgebra~(see Theorem~\ref{thm-NBHYBE-NBHLbi}). Section~\ref{sec-3} focuses on the differential setting. We first study differential Lie algebras and their representations. We then introduce Manin triples of differential Lie algebras, differential Lie bialgebras, and matched pairs of differential Lie algebras, and establish the equivalence among these three structures (see Theorem~\ref{th-equiva-diff-Lie}). Lastly, we  introduce the concept of $D$-differential classical Yang-Baxter equations~(abbr.~DD-CYBEs)~(Definition~\ref{def-diff-YBE}) and show that this equation gives rise to a differential-Lie bialgebra (see Theorem~\ref{thm-DYBE-diffLiebi}). Section~\ref{sec-4} mainly summarizes the relationships among various structures in this paper via two commutative diagrams

$\mathbf{Notation.}$ Throughout this paper, we will always work over a base field \bfk. All algebras, linear spaces etc. will be over \bfk. Denote by $\mathbb{N^{+}}$ the set of positive integers, unless otherwise specified. For the composition of two maps $ p $ and $ q $, we will write either $ p\circ q $ or simply $ pq$. In this paper, all considered vector space are finite dimensional, and we denote by $\mathfrak{gl}(V)$ the general linear Lie algebra consisting of all linear maps from vector space $V$ to $V$. Denote by $\langle \cdot, \cdot \rangle$ the natural pair between the vector space $L$ and its dual space $L^*$. For convenience, we write $r=\sum r_{1}\otimes r_{2}$ as $r=r_{1}\otimes r_{2}$.

\section{Nijenhuis BiHom-Lie algebras and Nijenhuis BiHom-Lie bialgebras}\label{sec-2}

\subsection{BiHom-Lie bialgebras}
In this section, we first recall the definitions of BiHom-Lie algebras and BiHom-Lie coalgebras, then we give some results of them. Finally, we propose the concept of BiHom-Lie bialgebras.

\begin{defn}\cite{gmmp}\label{defn-BHL}
A {\bf BiHom-Lie algebra} is a $4$-tuple $(L,[-,-],\alpha,\beta)$, where $L$ is a vector space, $[-,-]: L\otimes L\rightarrow L$ is a bilinear map and $\alpha,\,\beta$ are two commuting linear maps on $L$, satisfying
\begin{align}
\alpha([x,y])&=[\alpha(x),\alpha(y)],\quad \beta([x,y])=[\beta(x),\beta(y)],\qquad \text{(multiplicativity)}\label{eq-BH-multi}\\
[\beta(x),\alpha(y)]&=-[\beta(y),\alpha(x)],\qquad \text{(BiHom-antisymmetry)}\label{eq-BH-anti}\\
[\beta^2(x),[\beta(y),\alpha(z)]]&+[\beta^2(y),[\beta(z),\alpha(x)]]+[\beta^2(z),[\beta(x),\alpha(y)]]=0,\qquad \text{(BiHom-Jacobi condition)}\label{eq-BH-Jacobi}
\end{align}
for any $x,y,z\in L$. Further, a BiHom-Lie algebra $(L, [-,-],\alpha,\beta)$ is called an {\bf involutive BiHom-Lie algebra} if $\alpha^2=\beta^2=\mathrm{id}_{L}$.
\end{defn}

\begin{remark}
	Let $(L,[-,-],\alpha,\beta)$ be a BiHom-Lie algebra, and the linear map $\alpha$ be invertible. For any $x,y_{1},y_{2},\cdots,y_{n}\in L,\, n\in \mathbb{N^{+}}$, we define the {\bf adjoint diagonal action} $\ad^{(n)}_{x}: L^{\otimes n}\rightarrow L^{\otimes n}$ as follows.

\[
\mathrm{ad}_x^{(n)}\left(y_1 \otimes y_2 \otimes \cdots \otimes y_n\right) :=
\begin{cases}
\left[ x, y_1 \right], \;\; n=1,\\
\sum\limits_{i=1}^{n} \beta(y_1) \otimes \cdots \otimes \beta(y_{i-1}) \otimes [\alpha^{-i}{\beta^i(x)},y_i] \otimes \beta(y_{i+1}) \otimes \cdots \otimes \beta(y_n),\;\; n\geq 2.
\end{cases}
\]

	For example, when $n=2$, we have
	\[\ad^{(2)}_{x}(y_{1}\otimes y_{2})=[\alpha^{-1}\beta(x),y_{1}]\otimes\beta(y_{2})+ \beta(y_{1})\otimes [\alpha^{-2}\beta^{2}(x),y_{2}],\quad \text{for any }y_{1},\,y_{2}\in L.\]
\end{remark}

We now extend the notion of BiHom-Lie algebras to BiHom-Lie bialgebras. First we recall the concept of BiHom-Lie coalgebras.

\begin{defn}\cite{BHL-coalg}\label{defn-BHLco}
	A {\bf BiHom-Lie coalgebra} is a $4$-tuple $(L,\Delta,\alpha,\beta)$, where $L$ is a vector space, $\Delta: L\rightarrow L\otimes L$ is a linear map and $\alpha,\,\beta$ are two commuting linear map on $L$, satisfying
	\begin{align}
		&\hspace{2cm}\Delta\circ \alpha=(\alpha\otimes\alpha)\circ \Delta,\quad \Delta\circ \beta=(\beta\otimes \beta)\circ\Delta \qquad \text{(comultiplicativity)}\label{eq-comulti}\\
		&\hspace{2cm}(\beta\otimes \alpha)\circ\Delta+\tau(\beta\otimes \alpha)\circ \Delta=0\qquad \text{(BiHom-co-antisymmetry)}\label{eq-anti}\\
	&(\id\otimes\beta\otimes \alpha)(\beta^{2}\otimes \Delta)\circ \Delta+(\tau\otimes \id)(\id\otimes \tau)(\id\otimes\beta\otimes \alpha)(\beta^{2}\otimes \Delta)\circ \Delta\nonumber\\
	&+(\id \otimes \tau)(\tau\otimes \id)(\id\otimes\beta\otimes \alpha)(\beta^{2}\otimes \Delta)\circ \Delta=0,\qquad\text{(BiHom-coJacobi condition)}\label{eq-cojacobi}
	\end{align}
where $\tau(x\otimes y):=y\otimes x$, for any $x,y\in L$.
\end{defn}

As is well known, there exists the duality between Lie algebras and Lie coalgebras. For BiHom-Lie algebras and BiHom-Lie coalgebras, we have the following results.

\begin{prop}\label{prop-dual-BHL}
Let $L$ be a finite-dimensional vector space and $L^*$ denote its dual space. Then the $4$-tuple $(L^*,[-,-]_*,\alpha^*,\beta^*)$ is a BiHom-Lie algebra if and only if the triple $(L,\Delta,\alpha,\beta)$ is a BiHom-Lie coalgebra. Here, the linear map $\Delta: L \to L \otimes L$ is the dual of the bracket $[-,-]_*: L^* \otimes L^* \to L^*$, while $\alpha$ and $\beta$ are the linear duals of the linear maps $\alpha^*$ and $\beta^*$ respectively, i.e., $[-,-]_{*}=\Delta^*$.
\end{prop}

\begin{proof}
	First of all, for any $a^*,\,b^*,\,c^*\in L^*,\;x\in L$, it is obvious that the linear maps $\alpha,\,\beta$ are commutative if and only if $\alpha^*,\,\beta^*$ are commutative.
	Next, we consider Eq.~\eqref{eq-comulti} as follows.
	\begin{align*}
		&\langle  (\alpha^*\Delta^*-\Delta^*\circ(\alpha^*\otimes\alpha^*))(a^*\otimes b^*),x\rangle  \\
=&\ \langle  \Delta^*(a^*\otimes b^*),\alpha(x)\rangle  -\langle  (\alpha^*\otimes \alpha^*)(a^*\otimes b^*),\Delta(x)\rangle  \\
=&\ \langle  a^*\otimes b^*,\Delta\alpha(x)-(\alpha\otimes \alpha)\Delta(x)\rangle.
	\end{align*}
Thus, \[\alpha^*\circ \Delta^*=\Delta^*\circ(\alpha^*\otimes \alpha^*)\Leftrightarrow\Delta\circ \alpha=(\alpha\otimes \alpha)\circ \Delta.\]

\noindent Similarly, \[\beta^*\circ \Delta^*=\Delta^*\circ(\beta^*\otimes \beta^*)\Leftrightarrow \Delta\circ \beta=(\beta\otimes \beta)\circ \Delta.\]

 For Eq.~\eqref{eq-anti}, we have
	\begin{align*}
        &\langle  [\beta^*(a^*),\alpha^*(b^*)]_{*}+[\beta^*(b^*),\alpha^*(a^*)]_{*},x\rangle  \\
		=&\ \langle  (\Delta^*\circ(\beta^*\otimes \alpha^*)+\Delta^*\circ(\beta^*\otimes\alpha^*)\circ\tau)(a^*\otimes b^*),x\rangle  \\
		=&\ \langle  (\beta^*\otimes \alpha^*)(a^*\otimes b^*),\Delta(x)\rangle  +\langle  (\beta^*\otimes\alpha^*)\circ\tau(a^*\otimes b^*),\Delta(x)\rangle  \\
		=&\ \langle  a^*\otimes b^*,(\beta\otimes \alpha)\Delta(x)+\tau\circ(\beta\otimes\alpha)\Delta(x)\rangle.
	\end{align*}
Thus, $[\beta^*(a^*),\alpha^*(b^*)]_{*}=-[\beta^*(b^*),\alpha^*(a^*)]_{*}$ if and only if $(\beta\otimes \alpha)\circ\Delta+\tau\circ(\beta\otimes\alpha)\circ\Delta=0.$

For Eq.~\eqref{eq-cojacobi}, we have
\begin{align*}
    &\langle  [(\beta^*)^2(a^*),[\beta^*(b^*),\alpha^*(c^*)]_{*}]_{*}+[(\beta^*)^2(b^*),[\beta^*(c^*),\alpha^*(a^*)]_{*}]_{*}+[(\beta^*)^2(c^*),[\beta^*(a^*),\alpha^*(b^*)]_{*}]_{*},x\rangle  \\
	=&\ \langle  (\Delta^*((\beta^*)^2\otimes \Delta^*)(\id\otimes\beta^*\otimes\alpha^*)+\Delta^*((\beta^*)^2\otimes \Delta^*)(\id\otimes\beta^*\otimes\alpha^*)(\id\otimes\tau)(\tau\otimes \id)\\
	&+\Delta^*((\beta^*)^2\otimes \Delta^*)(\id\otimes\beta^*\otimes\alpha^*)(\tau\otimes \id)(\id\otimes \tau))(a^*\otimes b^*\otimes c^*),x\rangle  \\
	=&\ \langle  (((\beta^*)^2\otimes \Delta^*)(\id\otimes\beta^*\otimes\alpha^*)+((\beta^*)^2\otimes \Delta^*)(\id\otimes\beta^*\otimes\alpha^*)(\id\otimes\tau)(\tau\otimes \id)\\
	&+((\beta^*)^2\otimes \Delta^*)(\id\otimes\beta^*\otimes\alpha^*)(\tau\otimes \id)(\id\otimes \tau))(a^*\otimes b^*\otimes c^*),\Delta(x)\rangle  \\
	=&\ \langle  a^*\otimes b^*\otimes c^*,(\id\otimes\beta\otimes \alpha)(\beta^{2}\otimes \Delta)\circ \Delta(x)+(\tau\otimes \id)(\id\otimes \tau)(\id\otimes\beta\otimes \alpha)(\beta^{2}\otimes \Delta)\circ \Delta(x)\\
	&+(\id \otimes \tau)(\tau\otimes \id)(\id\otimes\beta\otimes \alpha)(\beta^{2}\otimes \Delta)\circ \Delta(x)\rangle.
\end{align*}
Thus, Eqs.~(\ref{eq-BH-multi})-(\ref{eq-BH-Jacobi}) hold for $(L^*, [-,-]_{*},\alpha^*,\beta^*)$ if and only if Eqs.~(\ref{eq-comulti})-(\ref{eq-cojacobi}) hold for $(L,\Delta,\alpha,\beta)$. By Definition~\ref{defn-BHL} and Definition~\ref{defn-BHLco}, we complete this proof.
\end{proof}

\begin{coro}\label{dual-Liealg-Liecoalg}
Under the same conditions as in Proposition~\ref{prop-dual-BHL}, the pair $(L^{*},[-,-]_{*})$ is a Lie algebra if and only if $(L,\Delta)$ is a Lie coalgebra, where $\Delta^{*}=[-,-]_{*}$.
\end{coro}
\begin{proof}
We complete this proof by taking $\alpha=\beta=\id_{L}$ in Proposition~\ref{prop-dual-BHL}.
\end{proof}

In order to introduce the Yau-twist property of Lie coalgebras, let us first recall the Yau-twist property of Lie algebras.

\begin{prop}\cite[Proposition~3.16]{gmmp}\label{Yautwist-Liealg}
    Let $(L, [-,-])$ be a Lie algebra, and let $\alpha,\beta: L\to L$ be two commuting linear maps satisfying the compatibility conditions
    \[
    \alpha([a,b]) = [\alpha(a),\alpha(b)],\qquad
    \beta([a,b]) = [\beta(a),\beta(b)]
    \]
    for all $a,b\in L$. Define a bilinear operation $\{-,-\}: L\otimes L\to L$ by
    \[
    \{a,b\} := [\alpha(a),\beta(b)]
    \]
    for all $a,b\in L$. Then $(L,\{-,-\},\alpha,\beta)$ is a BiHom-Lie algebra.
\end{prop}

\begin{prop}\label{prop-Yautwist-Liecoalg}
    Let $(L, \Delta)$ be a Lie coalgebra, and let $\alpha,\beta: L\to L$ be two commuting Lie coalgebra homomorphisms, i.e.,
    \[
    \Delta\circ\alpha = (\alpha\otimes\alpha)\circ\Delta,\qquad
    \Delta\circ\beta = (\beta\otimes\beta)\circ\Delta.
    \]
    Define a new coproduct $\blacktriangle: L\to L\otimes L$ on $L$ by
    \[
    \blacktriangle := (\alpha\otimes\beta)\circ\Delta.
    \]
    Then $(L, \blacktriangle, \alpha, \beta)$ is a BiHom-Lie coalgebra.
\end{prop}

\begin{proof}
Based on the known results, we first present the following commutative diagram which illustrates the core relationship of this proof, where our main objective is to verify the validity of the dashed morphism.

\begin{center}
\begin{tikzpicture}[node distance=7cm]
    \node (C) {$\text{Lie algebra}(L^*, [-,-]_*)$};
    \node (D) [right of=C] {$\text{Lie coalgebra}(L, \triangle)$};

    \node (A) [below of=C, node distance=2.5cm] {$\begin{aligned}
        &\text{BiHom-Lie algebra} \\
        &\;(L^*, \{-,-\}_*, \alpha^*, \beta^*)
    \end{aligned}$};
    \node (B) [below of=D, node distance=2.5cm] {$\begin{aligned}
        &\text{BiHom-Lie coalgebra} \\
        &\;\quad (L, \blacktriangle, \alpha, \beta)
    \end{aligned}$};

    \draw[<->] (C) -- node[above] {$\text{Corollary~\ref{dual-Liealg-Liecoalg}}$} (D);
    \draw[<->] (A) -- node[above] {$\text{Proposition~\ref{prop-dual-BHL}}$} (B);
    \draw[->] (C) -- node[left] {$\text{Proposition~\ref{Yautwist-Liealg}}$} (A);
    \draw[dashed, ->] (D) -- (B);
    \draw[<->] (A) -- 
                  node[below, font=\small] {$\{-,-\}_{*}=\blacktriangle^{*}$} (B);
\draw[->] (C) -- 
                  node[right, font=\small] {$\{a^*,b^*\}_{*}=[\alpha^*(a^*),\beta^*(b^*)]_{*}$} (A);
\draw[<->] (C) -- node[below] {$[-,-]_{*}=\Delta^{*}$} (D);

\end{tikzpicture}

\end{center}

For all $x\in L$ and $a^{*}, b^{*}\in L^{*}$, by the duality pairing and the above identities, we derive:
\begin{align*}
&\ \langle \blacktriangle(x),a^*\otimes b^*\rangle  \\
=&\ \langle x,\{a^*,b^*\}_{*}\rangle  \hspace{1cm}(\text{by }\{-,-\}_{*}=\blacktriangle^{*})\\
=&\ \langle x,[\alpha^*(a^*),\beta^*(b^*)]_{*}\rangle  \hspace{1cm}(\text{by }\{a^*,b^*\}_{*}=[\alpha^*(a^*),\beta^*(b^*)]_{*})\\
=&\ \langle \Delta(x),(\alpha^*\otimes\beta^*)(a^*\otimes b^*)\rangle  \hspace{1cm}(\text{by }[-,-]_{*}=\Delta^{*})\\
=&\ \langle (\alpha\otimes \beta)\circ \Delta(x),a^*\otimes b^*\rangle  .
\end{align*}

By applying Corollary~\ref{dual-Liealg-Liecoalg}, Proposition~\ref{Yautwist-Liealg} and Proposition~\ref{prop-dual-BHL} in turn,  we complete this proof.

\end{proof}

We now recall the definition of Lie bialgebras~\cite{liebialg}.

\begin{defn}
\begin{enumerate}
	\item A {\bf Lie bialgebra} is a triple $(L,[-,-],\Delta)$, where $(L,[-,-])$ is a Lie algebra, $(L,\Delta)$ is a Lie coalgebra and satisfy the following compatibility condition:
	\begin{align}\label{eq-Liebialg}
		\Delta([x,y])=(\ad_{x}\otimes \id+\id\otimes \ad_{x})\Delta(y)-(\ad_{y}\otimes \id+\id\otimes \ad_{y})\Delta(x),
	\end{align}
	where $\ad_{x}(y)=[x,y]$, for any $x,\,y\in L$.

\item Let $(L,[-,-],\Delta)$ and $(L',[-,-]',\Delta')$ be two Lie bialgebras, then linear map $f:L\rightarrow L'$ is a {\bf homomorphism of Lie bialgebras} if for any $x,y\in L$,
\begin{align*}
	f([x,y])=[f(x),f(y)]',\quad \Delta'\circ f(x)=(f\otimes f)\circ \Delta(x).
\end{align*}
\end{enumerate}
\end{defn}

Next, by the definition of coboundary operator about BiHom-Lie algebras in~\cite{BHL}, we give the concept of BiHom-Lie bialgebras.

\begin{defn}\label{def-BHLbialg}
	A {\bf BiHom-Lie bialgebra} is a 5-tuple $(L,[-,-],\Delta,\alpha,\beta)$, where
 \begin{enumerate}
 \item $(L,[-,-],\alpha,\beta)$ is a BiHom-Lie algebra;
 \item $(L,\Delta,\alpha,\beta)$ is a BiHom-Lie coalgebra;
 \item linear map $\alpha$ is invertible;
  \item $\Delta$ is a 1-cocycle of BiHom-Lie algebra $(L,[-,-],\alpha,\beta)$ with values in $L\otimes L$, i.e.,
	\begin{align}\label{eq-BHLiebialg}
		\Delta([\alpha^{-1}\beta(x),y])=&\ (\mathrm{ad}_{\beta(x)}\otimes \beta+\beta\otimes \mathrm{ad}_{\alpha^{-1}\beta^2(x)})\Delta(y)-(\mathrm{ad}_{\beta(y)}\otimes\beta+\beta\otimes \mathrm{ad}_{\alpha^{-1}\beta^2(y)})\Delta(x),
\end{align}
	 for any $x,\,y\in L$.
\end{enumerate}
\end{defn}

Now we show that BiHom-Lie bialgebras can be obtained from Lie bialgebras as follows.

\begin{prop}\label{Yau-Lie bialg}
    Let $(L, [-,-], \Delta)$ be a Lie bialgebra, and let $\alpha,\beta: L\to L$ be two commuting Lie bialgebra homomorphisms with $\alpha$ invertible.
    Define new operations on $L$ as follows:
    \[
    \blacktriangle := (\alpha\otimes\beta)\circ\Delta, \qquad
    \{-,-\} := [-,-]\circ(\alpha\otimes\beta),
    \]
    and denote the adjoint map associated to $\{-,-\}$ by
    \[
    \mathrm{Ad}_x(y) := \{x,y\},\quad \text{for all}\,x,y\in L.
    \]
    Then $(L, \{-,-\}, \blacktriangle, \alpha, \beta)$ is a BiHom-Lie bialgebra, which is called the Yau twist of the Lie bialgebra $(L, [-,-], \Delta)$.
\end{prop}

\begin{proof}
	By~Proposition~\ref{Yautwist-Liealg} and Proposition~\ref{prop-Yautwist-Liecoalg}, we only need to prove that Eq.~(\ref{eq-BHLiebialg}) holds for $(L, \{-,-\}, \blacktriangle, \alpha, \beta)$. For any $x,y\in L$, we have
	\begin{align}
		\mathrm{Ad}_{x}(y)=&\ \{x,y\}=[\alpha(x),\beta(y)]=\mathrm{ad}_{\alpha(x)}\beta(y).\label{eq-Ad}\\
        \mathrm{Ad}_{x}^{(2)}=&\ \mathrm{Ad}_{\alpha^{-1}\beta(x)}\otimes \beta +\beta \otimes \mathrm{Ad}_{\alpha^{-2}\beta^{2}(x)}.\label{eq-Ad2}
	\end{align}
Thus
	\begin{align*}
		&\ \blacktriangle(\{\alpha^{-1}\beta(x),y\})\\
=&\ (\alpha\otimes \beta)\circ\Delta([\beta(x),\beta(y)])\\
		=&\ (\alpha\otimes\beta)\big((\mathrm{ad}_{\beta(x)}\otimes \mathrm{id}+\mathrm{id}\otimes \mathrm{ad}_{\beta(x)})\Delta\beta(y)-(\mathrm{ad}_{\beta(y)}\otimes \mathrm{id}+\mathrm{id}\otimes \mathrm{ad}_{\beta(y)})\Delta\beta(x)\big)\quad(\text{by Eq.~(\ref{eq-Liebialg})})\\
        =&\ (\alpha\otimes\beta)\big((\mathrm{ad}_{\beta(x)}\otimes \mathrm{id}+\mathrm{id}\otimes \mathrm{ad}_{\beta(x)})(\beta\otimes\beta)\Delta(y)-(\mathrm{ad}_{\beta(y)}\otimes \mathrm{id}+\mathrm{id}\otimes \mathrm{ad}_{\beta(y)})(\beta\otimes\beta)\Delta(x)\big)\\
        &\hspace{1cm}(\text{by Eq.~(\ref{eq-comulti})})\\
		=&\ (\alpha\otimes \beta)\big((\mathrm{ad}_{\beta(x)}\beta\otimes \beta+\beta\otimes \mathrm{ad}_{\beta(x)}\beta)\Delta(y)-(\mathrm{ad}_{\beta(y)}\beta\otimes\beta+\beta\otimes \mathrm{ad}_{\beta(y)}\beta)\Delta(x)\big)\\
		=&\ (\mathrm{ad}_{\alpha\beta(x)}\alpha\beta\otimes \beta^2+\alpha\beta\otimes \mathrm{ad}_{\beta^{2}(x)}\beta^2)\Delta(y)-(\mathrm{ad}_{\alpha\beta(y)}\alpha\beta\otimes\beta^2+\alpha\beta\otimes \mathrm{ad}_{\beta^2(y)}\beta^2)\Delta(x)\\
		=&\ (\mathrm{Ad}_{\beta(x)}\alpha\otimes\beta^2+\alpha\beta\otimes \mathrm{Ad}_{\alpha^{-1}\beta^2(x)}\beta)\Delta(y)-(\mathrm{Ad}_{\beta(y)}\alpha\otimes \beta^2+\alpha\beta\otimes \mathrm{Ad}_{\alpha^{-1}\beta^2(y)}\beta)\Delta(x)\quad(\text{by Eq.~(\ref{eq-Ad})})\\
		=&\ (\mathrm{Ad}_{\beta(x)}\otimes \beta+\beta\otimes \mathrm{Ad}_{\alpha^{-1}\beta^2(x)})(\alpha\otimes \beta)\Delta(y)-(\mathrm{Ad}_{\beta(y)}\otimes \beta+\beta\otimes \mathrm{Ad}_{\alpha^{-1}\beta^2(y)})(\alpha\otimes \beta)\Delta(x)\\
        =&\ (\mathrm{Ad}_{\beta(x)}\otimes \beta+\beta\otimes \mathrm{Ad}_{\alpha^{-1}\beta^2(x)})\blacktriangle(y)-(\mathrm{Ad}_{\beta(y)}\otimes \beta+\beta\otimes \mathrm{Ad}_{\alpha^{-1}\beta^2(y)})\blacktriangle(x)\\
		=&\ \mathrm{Ad}_{\alpha(x)}^{(2)}\blacktriangle(y)-\mathrm{Ad}_{\alpha(y)}^{(2)}\blacktriangle(x).\hspace{1cm}(\text{by Eq.~(\ref{eq-Ad2})})
	\end{align*}
	This completes the proof.
\end{proof}

Conversely, we give the following result.

\begin{prop}\label{Yau-Lie bialg-2}
Let $(L,\{-,-\},\blacktriangle,\alpha,\beta)$ be a BiHom-Lie bialgebra with invertible linear maps $\alpha$ and $\beta$. Define the linear maps $\Delta:L\rightarrow L\otimes L$ and $[-,-]:L\otimes L\rightarrow L$ by
	\[\Delta:=(\alpha^{-1}\otimes \beta^{-1})\circ \blacktriangle,\quad [-,-]:=\{-,-\}\circ (\alpha^{-1}\otimes \beta^{-1}).\]
	Then $(L,[-,-],\Delta)$ is a Lie bialgebra.
\end{prop}
\begin{proof}
	This is straightforward.
\end{proof}

\begin{coro}
	Let $(L,[-,-],\Delta)$ be a Lie bialgebra and $\alpha$ be a Lie bialgebra homomorphism. Then $L_{\alpha}=(L,[-,-]_{\alpha},\Delta_{\alpha})$ is a Hom-Lie bialgebra, where $[-,-]_{\alpha}:=\alpha\circ [-,-]$ and $\Delta_{\alpha}:=\Delta\circ\alpha.$
\end{coro}

\begin{proof}
	For any $x,y\in L,$ we have
	\begin{align*}
		\mathrm{Ad}_{x}(y)=[x,y]_{\alpha}=\alpha([x,y])=[\alpha(x),\alpha(y)]=\ad_{\alpha(x)}\alpha(y).
	\end{align*}
Hence, this proof follows from Proposition~\ref{Yau-Lie bialg} and the definition of adjoint diagonal action on Hom-Lie algebras~\cite{double-C}.
\end{proof}

\subsection{Nijenhuis BiHom-Lie algebras and their representations}
In this subsection, we mainly study the representation of Nijenhuis BiHom-Lie algebras.

First of all, let's recall the definition of representation of BiHom-Lie algebras.

\begin{defn}~\cite{BHL}\label{repre-BHL}
	A {\bf representation of BiHom-Lie algebra} $(L,[-,-],\alpha,\beta)$ is a $4$-tuple $(V,\rho,p,q)$, where $V$ is a vector space, $p,q: V\rightarrow V$ are two commuting linear maps and $\rho: L\rightarrow \mathfrak{gl}(V)$ is a linear map satisfying for any $x,y\in L,\, v\in V$,
	\begin{align}
		&p\big(\rho(x)v\big)=\rho\big(\alpha(x)\big)p(v),\quad q\big(\rho(x)v\big)=\rho\big(\beta(x)\big)q(v),\nonumber \\
		&\rho\big([\beta(x),y]\big)q(v)=\rho\big(\alpha\beta(x)\big)\rho(y)v-\rho\big(\beta(y)\big)
\rho\big(\alpha(x))v.\label{eq-BHLrepre-2}
	\end{align}
\end{defn}

Next, we give a specific representation of BiHom-Lie algebras.

\begin{prop}\label{BHL_repre-ad2}
	Let $(L, [-,-],\alpha,\beta)$ be an involutive BiHom-Lie algebra. Define the operation $\ad^{(2)}: L\rightarrow \mathfrak{gl}(L\otimes L)$ by
	\[\ad^{(2)}(x):=\ad^{(2)}_{x}:=\ad_{\alpha^{-1}\beta(x)}\otimes \beta+\beta\otimes \ad_{x},\qquad\text{ for any } x\in L.\]
	Then $(L\otimes L, \ad^{(2)},\alpha\otimes\alpha,\beta\otimes\beta)$ is a representation of BiHom-Lie algebra $L$.
\end{prop}
\begin{proof}
	For any $a\otimes b\in L\otimes L$ and $x\in L$, we have
	\begin{align*}
		(\alpha\otimes \alpha)\ad^{(2)}_{x}(a\otimes b)=&\ (\alpha\otimes \alpha)(\ad_{\alpha^{-1}\beta(x)}\otimes \beta+\beta\otimes \ad_{x})(a\otimes b)\\
		=&\ (\alpha \ad_{\alpha^{-1}\beta(x)}\otimes \alpha\beta+\alpha\beta\otimes \alpha \ad_{x})(a\otimes b)\\
		=&\ (\ad_{\beta(x)}\alpha\otimes \beta\alpha+\alpha\beta\otimes \ad_{\alpha^{-1}(x)}\alpha)(a\otimes b)\\
		=&\ (\ad_{\beta(x)}\otimes \beta+\beta\otimes \ad_{\alpha^{-1}(x)})(\alpha\otimes \alpha)(a\otimes b)\\
		=&\ \ad^{(2)}_{\alpha(x)}(\alpha\otimes \alpha)(a\otimes b).\hspace{1cm}(\text{by }\alpha^{2}=\id_{L})
	\end{align*}
	
	Similarly, we get $(\beta\otimes \beta)\ad^{(2)}_{x}(a\otimes b)=\ad^{(2)}_{\beta(x)}(\beta\otimes \beta)(a\otimes b).$
	
	So, by Definition~\ref{repre-BHL}, we only need to prove that Eq.~(\ref{eq-BHLrepre-2}) holds for $(L\otimes L, \ad^{(2)},\alpha\otimes\alpha,\beta\otimes\beta)$. For any $x,y\in L$, we have
	\begin{align}\label{LLad2-repre}
		&\ \ad^{(2)}_{\alpha\beta(x)}\ad^{(2)}_{y}(a\otimes b)-\ad^{(2)}_{\beta(y)}\ad^{(2)}_{\alpha(x)}(a\otimes b)\nonumber\\
		=&\ \ad^{(2)}_{\alpha\beta(x)}\big([\alpha^{-1}\beta(y),a]\otimes \beta(b)+\beta(a)\otimes [x,b]\big)-\ad^{(2)}_{\beta(y)}\big([\beta(x),a]\otimes \beta(b)+\beta(a)\otimes [\alpha(x),b]\big)\nonumber\\
		=&\ \big([x,[\alpha^{-1}\beta(y),a]]-[\alpha^{-1}(y),[\beta(x),a]]\big)\otimes b+a\otimes \big([\alpha^{-1}\beta(x),[y,b]]\nonumber\\
		& -[\beta(y),[\alpha^{-1}\beta^{2}(x),b]]\big)+[\alpha^{-1}(y),\beta(a)]\otimes [\alpha^{-1}\beta(x),\beta(b)]+[x,\beta(a)]\otimes [\beta(y),\beta(b)]\nonumber\\
		&-[x,\beta(a)]\otimes [\beta(y),\beta(b)]-[\alpha^{-1}(y),\beta(a)]\otimes [\alpha^{-1}\beta(x),\beta(b)]\hspace{1cm}(\text{by }\alpha^{2}=\beta^{2}=\id_{L})\nonumber\\
		=&\ \Big([x,[\alpha^{-1}\beta(y),a]]+[\alpha^{-1}(y),[\alpha^{-1}\beta(a),\alpha(x)]]\Big)\otimes b
		+a \otimes \Big([\alpha^{-1}\beta(x),[y,b]]\nonumber\\
		&+[\beta(y),[\alpha^{-1}\beta(b),\beta(x)]]\Big)\hspace{1cm}(\text{by Eq.~(\ref{eq-BH-anti})})\nonumber\\
		=&\ -[\beta^{2}\alpha^{-1}(a),[\beta^{-1}(x),y]]\otimes \beta^{2}(b)-\beta^{2}(a)\otimes [\alpha^{-1}\beta^{2}(b),[\alpha^{-1}\beta^{2}(x),\alpha^{-1}\beta(y)]]\nonumber\\
		& \hspace{1cm}(\text{by Eq.~(\ref{eq-BH-Jacobi})~and }\beta^{2}=\id_{L})\nonumber\\
		=&\ [[\alpha^{-1}(x),\alpha^{-1}\beta(y)],\beta(a)]\otimes \beta^{2}(b)+\beta^{2}(a)\otimes[[\alpha^{-2}\beta^{3}(x),\alpha^{-2}\beta^{2}(y)],\beta(b)]\hspace{1cm}(\text{by Eq.~(\ref{eq-BH-anti})})\nonumber\\
		=&\ [[\alpha^{-1}\beta^{2}(x),\alpha^{-1}\beta(y)],\beta(a)]\otimes \beta^{2}(b)+\beta^{2}(a)\otimes [[\alpha^{-2}\beta^{3}(x),\alpha^{-2}\beta^{2}(y)],\beta(b)]\hspace{1cm}(\text{by }\beta^{2}=\id_{L}) \nonumber\\
		=&\ \Big(\ad_{\alpha^{-1}\beta([\beta(x),y])}\otimes \beta+\beta\otimes \ad_{\alpha^{-2}\beta^{2}([\beta(x),y])}\Big)(\beta\otimes \beta)(a\otimes b)\nonumber\\
		=&\ \ad^{(2)}_{[\beta(x),y]}(\beta\otimes \beta)(a\otimes b).\hspace{1cm}(\text{by }\alpha^{2}=\beta^{2}=\id_{L})
	\end{align}
	Thus, Eq.~(\ref{eq-BHLrepre-2}) holds. We complete this proof.
\end{proof}

Next, we introduce the concept of Nijenhuis BiHom-Lie algebras.

\begin{defn}\label{nijenhuis-bihom-lie}
A {\bf Nijenhuis BiHom-Lie algebra} is a 5-tuple $(L,[-,-],N,\alpha,\beta)$, where $(L,[-,-],\\
\alpha,\beta)$ is a BiHom-Lie algebra and linear map $N:L\rightarrow L$ satisfies for any $x,y\in L$,
\[\alpha\circ N=N\circ \alpha,\quad \beta\circ N=N\circ\beta,\]
\begin{align}
	[N(x),N(y)]=N([N(x),y]+[x,N(y)])-N^2([x,y]).\label{eq-Nij}
\end{align}
\end{defn}

\begin{exam}
Let $A$ be a vector space over a field $\bfk$, with a basis $\{e_1,e_2\}$.
Let $\alpha,\beta: A\to A$ be linear maps defined on the basis by
\begin{equation*}
\alpha(e_1)=e_1, \qquad
\alpha(e_2)=\frac{1}{m}e_1+\frac{n-1}{n}e_2,\qquad 
\beta(e_1)=e_1, \qquad
\beta(e_2)=e_2, 
\end{equation*}
and let the bilinear bracket $[-,-]: A\otimes A\to A$ be given by
\begin{equation*}
[e_1,e_1]=0, \quad [e_1,e_2]=m e_2 - n e_1, \quad
[e_2,e_1]=(n-1)e_1 - \frac{m(n-1)}{n}e_2,
\quad
[e_2,e_2]=-\frac{n}{m}e_1 + e_2, \mlabel{eq:lieb}
\end{equation*}
where $m,n\in\bfk$ are parameters satisfying $m\neq0$, $n\neq0$ and $n\neq1$.
Then $(A,[-,-],\alpha,\beta)$ is a BiHom-Lie algebra.
Furthermore, suppose that $N: A\to A$ is a linear map defined on the basis as
\begin{equation*}
N(e_1)=e_1, \qquad
N(e_2)=\frac{n(1-m)}{m}e_1+me_2, \mlabel{eq:nijen}
\end{equation*}
then $(A,[-,-],N,\alpha,\beta)$ is a Nijenhuis BiHom-Lie algebra.
\end{exam}

\begin{defn}\label{repre-NBL}
	A {\bf representation of a Nijenhuis BiHom-Lie algebra} $(L,[-,-],N,\alpha,\beta)$ is a 5-tuple $(V,\rho,\eta,p,q)$, where $(V,\rho,p,q)$ is a representation of the BiHom-Lie algebra $(L,[-,-],\alpha,\beta)$ and $\eta: V\rightarrow V$ is a linear map satisfying for any $ x\in L,\; v\in V,$
	\[\eta\circ p=p\circ \eta,\quad\eta\circ q=q\circ\eta,\]
	\begin{align}
		\rho\big(N(x)\big)\eta(v)=\eta\big(\rho(N(x))v\big)+\eta\big(\rho(x)\eta(v)\big)-\eta^2\big(\rho(x)v\big).
\label{eq-NijBHL-repre}
	\end{align}
\end{defn}

\begin{remark}
In Definition~\ref{repre-NBL}, if the linear maps $\alpha,\;\beta:L\rightarrow L$ satisfy $\alpha=\beta=\id_{L}$, then $L$ reduces to a Nijenhuis Lie algebra. Naturally, the triple $(V,\rho,\eta)$ is a representation of Nijenhuis Lie algebra $(L,[-,-],N).$
\end{remark}

Now we characterize the representations of Nijenhuis BiHom-Lie algebras via semi-direct products.
\begin{prop}\label{prop-semi-direct}
Let $(L,[-,-],N,\alpha,\beta)$ be a Nijenhuis BiHom-Lie algebra, and let $(V,p,q)$ be a vector space equipped with two commuting linear maps $p,q:V\to V$. Let $\rho: L\rightarrow \mathfrak{gl}(V)$ and $\eta:V\rightarrow V$ be linear maps. Suppose that $\alpha$ and $q$ are bijective. Then $(V,\rho,\eta,p,q)$ is a representation of $(L,[-,-],N,\alpha,\beta)$ if and only if the semi-direct product $L\ltimes V:=(L\oplus V,[-,-]_{L\oplus V},N+\eta,\alpha+p,\beta+q)$ is a Nijenhuis BiHom-Lie algebra, where the operations on $L\oplus V$ are defined for all $x,y\in L$ and $a,b\in V$ by:
\begin{align*}
	[x+a,y+b]_{L\oplus V} :=& [x,y]+\rho(x)b-\rho\big(\alpha^{-1}\beta(y)\big)\big(pq^{-1}(a)\big),\\
(N+\eta)(x+a) := N(x)+\eta(a),\;&
(\alpha+p)(x+a) := \alpha(x)+p(a),\;
(\beta+q)(x+a) := \beta(x)+q(a).
\end{align*}
\end{prop}

\begin{proof}
By \cite[Proposition4.4]{BHL}, $(L\oplus V,[-,-]_{L\oplus V},\alpha+p,\beta+q)$ is a BiHom-Lie algebra if and only if $(V,\rho,p,q)$ is a representation of the BiHom-Lie algebra $(L,[-,-],\alpha,\beta)$. We only need to verify the Nijenhuis condition (\ref{eq-Nij}) for $L\ltimes V$.

Assume $(V,\rho,\eta,p,q)$ is a representation of $(L,[-,-],N,\alpha,\beta)$. For all $x,y\in L$ and $a,b\in V$, compute the left-hand side (LHS) and right-hand side (RHS) of Eq.~(\ref{eq-Nij}) for $L\ltimes V$:
\begin{align*}
\text{LHS} =&\ (N+\eta)\big([(N+\eta)(x+a),y+b]_{L\oplus V}+[x+a,(N+\eta)(y+b)]_{L\oplus V}\big)-(N+\eta)^2[x+a,y+b]_{L\oplus V}\\
		=&\ (N+\eta)\big([N(x)+\eta(a),y+b]_{L\oplus V}+[x+a,N(y)+\eta(b)]_{L\oplus V}\big)-(N+\eta)^2[x+a,y+b]_{L\oplus V}\\
		=&\ (N+\eta)\big([N(x),y]+\rho(N(x))b-\rho\big(\alpha^{-1}\beta(y)\big)\big(pq^{-1}\eta(a)\big)+[x,N(y)]+\rho(x)\eta(b)\\
		&-\rho\big(\alpha^{-1}\beta N(y)\big)\big(pq^{-1}(a)\big)\big)-N^2([x,y])-\eta^2\big(\rho(x)b-\rho\big(\alpha^{-1}\beta(y)\big)\big(pq^{-1}(a)\big)\big)\\
		=&\ N\big([N(x),y]+[x,N(y)]\big)-N^2([x,y])+\eta\big(\rho(N(x))b\big)+\eta\big(\rho(x)\eta(b)\big)-\eta^2\big(\rho(x)b\big)\\
		& -\eta\bigg(\rho\big(\alpha^{-1}\beta(y)\big)\eta\big(pq^{-1}(a)\big)\bigg)-\eta\bigg(N\big(\alpha^{-1}\beta(y)\big)\big(pq^{-1}(a)\big)\bigg)+\eta^2\bigg(\rho\big(\alpha^{-1}\beta(y)\big)\big(pq^{-1}(a)\big)\bigg)\\
		=&\ [N(x),N(y)]+\rho\big(N(x)\big)\eta(b)-\rho\big(N\big(\alpha^{-1}\beta(y)\big)\big)\eta\big(pq^{-1}(a)\big)\hspace{1cm}(\text{by Eqs.~(\ref{eq-Nij})-(\ref{eq-NijBHL-repre})})\\
		=&\ [N(x),N(y)]+\rho\big(N(x)\big)\eta(b)-\rho\big(\alpha^{-1}\beta N(y)\big)\big(pq^{-1}\eta(a)\big)\\
		=&\ [N(x)+\eta(a),N(y)+\eta(b)]_{L\oplus V}\\
		=&\ [(N+\eta)(x+a),(N+\eta)(y+b)]_{L\oplus V}= \text{RHS}.
	\end{align*}
Thus, $L\ltimes V$ satisfies the Nijenhuis condition (\ref{eq-Nij}).

Conversely, assume $L\ltimes V$ is a Nijenhuis BiHom-Lie algebra. Then it satisfies the BiHom-Lie algebra axioms and the Nijenhuis condition:
\begin{align}
(\alpha+p)(\beta+q)(x+a) &= (\beta+q)(\alpha+p)(x+a), \label{eq-1} \\
[(\beta+q)(x+a),(\alpha+p)(y+b)]_{L\oplus V} &=-[(\beta+q)(y+b),(\alpha+p)(x+a)]_{L\oplus V}, \label{eq-2} \\
[(\beta+q)^2(x+a),[(\beta+q)(y+b),(\alpha+p)(z+c)]_{L\oplus V}]_{L\oplus V} \nonumber \\
+[(\beta+q)^2(y+b),[(\beta+q)(z+c),(\alpha+p)(x+a)]_{L\oplus V}]_{L\oplus V} \nonumber \\
+[(\beta+q)^2(z+c),[(\beta+q)(x+a),(\alpha+p)(y+b)]_{L\oplus V}]_{L\oplus V} &= 0, \label{eq-3} \\
(\alpha+p)\big([x+a,y+b]_{L\oplus V}\big) &= [\alpha(x)+p(a),\alpha(y)+p(b)]_{L\oplus V}, \label{eq-4} \\
(\beta+q)\big([x+a,y+b]_{L\oplus V}\big) &= [\beta(x)+q(a),\beta(y)+q(b)]_{L\oplus V}, \label{eq-5} \\
(\alpha+p)(N+\eta)(x+a) &= (N+\eta)(\alpha+p)(x+a), \label{eq-6} \\
(\beta+q)(N+\eta)(x+a) &= (N+\eta)(\beta+q)(x+a), \label{eq-7} \\
[(N+\eta)(x+a),(N+\eta)(y+b)]_{L\oplus V} &= (N+\eta)\big([(N+\eta)(x+a),y+b]_{L\oplus V} \nonumber \\
+ [x+a,(N+\eta)(y+b)]_{L\oplus V}\big) - (N+\eta)^2[x+a,y+b]_{L\oplus V}. \label{eq-8}
\end{align}
Substituting specific values of $x,y,z\in L$ and $a,b,c\in V$  into the above equations yields the desired conclusions, summarized in the following table: 
\begin{table}[h]
\centering
\footnotesize 
\setlength{\tabcolsep}{5pt}  
\renewcommand{\arraystretch}{1.0}
\begin{tabular}{|p{7.5cm}|p{8.0cm}|}
\hline
\textbf{Substitution Conditions} & \textbf{Derived Results} \\
\hline
$a=0$ in Eq.~(\ref{eq-1}), Eqs.~(\ref{eq-6})-(\ref{eq-7}); & \multirow{3}{*}{$(L, [-,-], N, \alpha, \beta)$ is a Nijenhuis BiHom-Lie algebra.} \\
$a=b=c=0$ in Eq.~(\ref{eq-3}); & \\
$a=b=0$ in Eq.~(\ref{eq-2}), Eqs.~(\ref{eq-4})-(\ref{eq-5}) and Eq.~(\ref{eq-8}). & \\
\hline
$x=0$ in Eq.~(\ref{eq-1}), Eqs.~(\ref{eq-6})-(\ref{eq-7}); &\\
$a=b=z=0$ in Eq.~(\ref{eq-3}); & \multirow{1}{*}{$(V, \rho, \eta, p, q)$ is a representation of $(L, [-,-], N, \alpha, \beta)$}. \\
$a=y=0$ in Eqs.~(\ref{eq-4})-(\ref{eq-5}) and Eq.~(\ref{eq-8}). & \\
\hline
\end{tabular}
\end{table}

\noindent This completes the proof.
\end{proof}

Inspired by Definition~\ref{repre-NBL}, we introduce the dual representation as follows.
Let $L$ and $V$ be vector spaces, and let
\[
\rho: L \to \mathfrak{gl}(V), \qquad \eta, p, q: V \to V
\]
be linear maps. Their duals on the dual space $V^* = \operatorname{Hom}(V,\mathbf{k})$ are defined by
\begin{align}
    \langle \rho^*(x)(v^*), v \rangle &= - \langle v^*, \rho(x)(v) \rangle, \label{star-rho}\\
    \langle \eta^*(v^*), v \rangle &= \langle v^*, \eta(v) \rangle, \nonumber\\
    \langle p^*(v^*), v \rangle &= \langle v^*, p(v) \rangle, \nonumber\\
    \langle q^*(v^*), v \rangle &= \langle v^*, q(v) \rangle,\nonumber
\end{align}
for all $x \in L$, $v \in V$, and $v^* \in V^*$. Here $\rho^*: L \to \mathfrak{gl}(V^*)$ is the dual representation induced by $\rho$.

 \begin{prop}\label{prop-Nijrepre-iff-repre+eq}
 	Let $(L,[-,-],N,\alpha,\beta)$ be an involutive Nijenhuis BiHom-Lie algebra. Then $(V^*,\rho^*,\eta^*,p^*,q^*)$ is a representation of Nijenhuis BiHom-Lie algebra $(L,[-,-],N,\alpha,\beta)$ if and only if $(V,\rho,p,q)$ is a representation of BiHom-Lie algebra $(L,[-,-],\alpha,\beta)$ and
 	\begin{align}\label{eq-admissible}
 		\eta\big(\rho(N(x))u\big)+\rho(x)\big(\eta^2(u)\big)=\rho\big(N(x)\big)\eta(u)+\eta\big(\rho(x)\eta(u)\big),
 	\end{align}
 for any $x\in L,\,u\in V.$
 \end{prop}

 \begin{proof}
For any $x,y\in L,\;u\in V,\; v^*\in V^*$, we first consider the duality of representations of BiHom-Lie algebras. If $(V,\rho,p,q)$ is a representation of $(L,[-,-],\alpha,\beta)$, then we have
		\begin{align*}
		\langle \rho^*\big(\alpha(x)\big)p^*(v^*)-p^*\big(\rho^*(x)v^*\big),u\rangle
		=&\ -\langle p^*(v^*),\rho\big(\alpha(x)\big)u\rangle  -\langle \rho^*(x)v^*,p(u)\rangle  \\
		=&\ \langle v^*,\rho(x)p(u)-p\big(\rho(\alpha(x))u\big)\rangle  \\
		=&\ \langle v^*,\rho(x)p(u)-\rho\big(\alpha^2(x)\big)p(u)\rangle  \\
		=&\ \langle v^*,\rho(x)p(u)-\rho(x)p(u)\rangle  =0.
	\end{align*}
Thus, $\rho^*\big(\alpha(x)\big)p^*(v^*)=p^*\big(\rho^*(x)v^*\big)$. Similarly, we get $\rho^*\big(\beta(x)\big)q^*(v^*)=q^*\big(\rho^*(x)v^*\big)$.

\noindent Also, we have
\[
\begin{split}
\langle \rho^*([\beta(x), y]) q^*(v^*), u \rangle
=&\ -\langle q^*(v^*), \rho([\beta(x), y]) u \rangle \hspace{1cm}(\text{by Eq.~(\ref{star-rho})})\\
=&\ -\langle v^*, q(\rho([\beta(x), y]) u) \rangle \\
=&\ -\langle v^*,\rho([\beta^{2}(x),\beta(y)])q(u)\rangle\\
=&\ -\langle v^*, \rho(\alpha\beta^{2}(x)) \rho(\beta(y))u - \rho(\beta^{2}(y)) \rho(\alpha\beta(x)) u \rangle \hspace{1cm} (\text{by Eq.~\eqref{eq-BHLrepre-2} for } V) \\
=&\ \langle \rho^*(\alpha\beta(x)) \rho^*(y) v^* - \rho^*(\beta(y)) \rho^*(\alpha(x)) v^*, u \rangle.\hspace{1cm}(\text{by }\beta^{2}=\mathrm{id}_{L})
\end{split}
\]
So, \(\rho^*([\beta(x), y]) q^*(v^*) = \rho^*(\alpha\beta(x)) \rho^*(y)v^* - \rho^*(\beta(y)) \rho^*(\alpha(x))v^*\). Therefore, $(V^*,\rho^*,p^*,q^*)$ is a representation of $(L,[-,-],\alpha,\beta)$.
Conversely, if \((V^*, \rho^*, p^*, q^*)\) is a representation of $L$, reverse the chains (using \(V = (V^*)^*\)) to recover the axioms of \(V\).

Next, we consider the Nijenhuis condition, for Eq.~(\ref{eq-NijBHL-repre}) and Eq.~(\ref{eq-admissible}), we have
\begin{align*}
	&\ \langle\rho^*\big(N(x)\big)\eta^*(v^*)+(\eta^*)^2\big(\rho^*(x)v^*\big)-\eta^*\big(\rho^*\big(N(x)\big)v^*\big)-\eta^*\big(\rho^*(x)\eta^*(v^*)\big),u\rangle\\
	=&\ -\langle\eta^*(v^*),\rho\big(N(x)\big)u\rangle+\langle\rho^*(x)v^*,\eta^2(u)\rangle-\langle\rho^*\big(N(x)\big)v^*,\eta(u)\rangle-\langle\rho^*(x)\eta^*(v^*),\eta(u)\rangle\\
	=&\ -\langle v^*,\eta\Big(\rho\big(N(x)\big)u\Big)\rangle-\langle v^*,\rho(x)\big(\eta^2(u)\big)\rangle+\langle v^*,\rho\big(N(x)\big)\eta(u)\rangle+\langle\eta^*(v^*),\rho(x)\eta(u)\rangle\\
	=&\ \langle v^*,-\eta\Big(\rho\big(N(x)\big)u\Big)-\rho(x)\big(\eta^2(u)\big)+\rho\big(N(x)\big)\eta(u)
+\eta\big(\rho(x)\eta(u)\big)\rangle.
\end{align*}
That is to say, Eq.~(\ref{eq-NijBHL-repre}) holds for $(V^*,\rho^*,\eta^*,p^*,q^*)$ if and only if Eq.~(\ref{eq-admissible}) holds.
This completes the proof.
\end{proof}
%
%
%


\begin{coro}\label{coro-dual-repre-Lie}
	Let $(L,[-,-])$ be a Lie algebra. Then $(V^*,\rho^*)$ is a representation of Lie algebra $(L,[-,-])$ if and only if $(V,\rho)$ is a representation of Lie algebra $L$.
\end{coro}

\begin{proof}
	We complete this proof by taking $\alpha=\beta=N=\id_{L},\, \eta^{*}=p^*=q^*=\id_{V^*},\,\eta=p=q=\id_{V}$ in Proposition~\ref{prop-Nijrepre-iff-repre+eq}.
\end{proof}

In particular, when the representation $(V,\rho,p,q)$ in Proposition~\ref{prop-Nijrepre-iff-repre+eq} is taken as the adjoint representation $(L, \mathrm{ad}, \alpha,\beta)$ with $\mathrm{ad}(x)(y):=\mathrm{ad}_{x}(y)=[x,y]$ for any $x,y\in L$, we can draw the following conclusion:

\begin{coro}\label{coro-adjoint}
Let $(L,[-,-],N,\alpha,\beta)$ be an involutive Nijenhuis BiHom-Lie algebra, and let $S:L\rightarrow L$ be a linear map.
Then  $(L^*,\mathrm{ad}^*,S^*,\alpha^*,\beta^*)$ is a representation of $(L,[-,-],N,\alpha,\beta)$ if and only if $S$ satisfies
\begin{equation}\label{eq-adjoint-admissible}
S\big([N(x),y]\big) + [x,S^2(y)] = [N(x),S(y)] + S\big([x,S(y)]\big), \qquad \text{for all}\,x,y\in L.
\end{equation}
\end{coro}

\begin{proof}
Since $(L,[-,-],\alpha,\beta)$ is a BiHom-Lie algebra, the adjoint map $\mathrm{ad}:L\rightarrow\mathfrak{gl}(L)$ gives a representation $(L,\mathrm{ad},\alpha,\beta)$.
By Proposition~\ref{prop-Nijrepre-iff-repre+eq}, its dual $(L^*,\mathrm{ad}^*,\alpha^*,\beta^*)$ is also a representation of $(L,[-,-],\alpha,\beta)$.

It remains to verify the compatibility condition between $\mathrm{ad}^*$ and $S^*$ induced by the Nijenhuis operator $N$.
Applying Proposition~\ref{prop-Nijrepre-iff-repre+eq} with $\rho=\mathrm{ad}$ and $\eta=S$, we say that $(L^*,\mathrm{ad}^*,S^*,\alpha^*,\beta^*)$ satisfies the Nijenhuis compatibility condition if and only if
\begin{align}\label{26-ad}
S\big(\mathrm{ad}(N(x))y\big) + \mathrm{ad}(x)S^2(y) = \mathrm{ad}(N(x))S(y) + S\big(\mathrm{ad}(x)S(y)\big)
\end{align}
holds for all $x,y\in L$.
Rewriting $\mathrm{ad}(x)y = [x,y]$, Eq.~(\ref{26-ad}) is precisely Eq.~\eqref{eq-adjoint-admissible}. This completes the proof.
\end{proof}

\begin{defn}\label{def-admissible}
	Let $(L,[-,-],N,\alpha,\beta)$ be an involutive Nijenhuis BiHom-Lie algebra, $(V,\rho,p,q)$ be a representation of BiHom-Lie algebra $(L,[-,-],\alpha,\beta)$ and $\eta: V\rightarrow V$ be a linear map.
	\begin{enumerate}
		\item If $(V^*,\rho^*,\eta^*,p^*,q^*)$ is a representation of $(L,[-,-],N,\alpha,\beta)$, i.e., Eq.~(\ref{eq-admissible}) holds, then we call that {\bf $\eta$ is admissible to the Nijenhuis BiHom-Lie algebra $(L,[-,-],N,\alpha,\beta)$ on $(V,\rho,p,q)$}, or {\bf $(L,[-,-],N,\alpha,\beta)$ is $\eta$-admissible on $(V,\rho,p,q)$}.
		\item If there is a linear map $S: L\rightarrow L$ satisfying Eq.~(\ref{eq-adjoint-admissible}) such that $(L^*,\ad^*,S^*,\alpha^*,\beta^*)$ is a representation of $(L,[-,-],N,\alpha,\beta)$, then we say that {\bf $S$ is adjoint-admissible to $(L,[-,-],N,\alpha,\beta)$}, or {\bf $(L,[-,-],N,\alpha,\beta)$ is $S$-adjoint-admissible}.
	\end{enumerate}
\end{defn}

In order to better understand the concept of Nijenhuis BiHom-Lie algebras, we give the definition of Nijenhuis BiHom-Lie coalgebras as follows.

\begin{defn}\label{def-NijBHLcoalg}
    A \textbf{Nijenhuis BiHom-Lie coalgebra} is a $5$-tuple $(L, \Delta, S, \alpha, \beta)$, where $(L, \Delta, \alpha, \beta)$ is a BiHom-Lie coalgebra, and the linear map $S: L\to L$ satisfies
    \begin{align}
        (S\otimes S)\circ\Delta + \Delta\circ S^{2} = (S\otimes \mathrm{id})\circ\Delta\circ S + (\mathrm{id}\otimes S)\circ\Delta\circ S,
        \label{NBL-coalg}
    \end{align}
    for all $x\in L$.
\end{defn}

Next, we find the following duality between Nijenhuis BiHom-Lie algebras and Nijenhuis BiHom-Lie coalgebras.

\begin{prop}\label{prop-dual-Nij}
	Let $\Delta: L\rightarrow L\otimes L$ be the linear dual of $[-,-]_{*}:L^{*}\otimes L^{*}\rightarrow L^{*}$. Then $(L^{*},[-,-]_{*},N^{*},\alpha^{*},\beta^{*})$ is a Nijenhuis BiHom-Lie algebra if and only if $(L,\Delta,N,\alpha,\beta)$ is a Nijenhuis BiHom-Lie coalgebra.
\end{prop}

\begin{proof}
	For any $x\in L,\; a^*,b^ *\in L^{*}$, we have
	\begin{align*}
		\langle N^*\alpha^*(a^*)-\alpha^{*}N^{*}(a^{*}),x\rangle =&\ \langle N^*\alpha^*(a^*),x\rangle-\langle\alpha^*N^*(a^*),x\rangle\\
		=&\ \langle\alpha^*(a^*),N(x)\rangle-\langle N^*(a^*),\alpha(x)\rangle\\
		=&\ \langle a^*,\alpha N(x)-N\alpha(x)\rangle.
	\end{align*}
Similarly, $\langle N^*\beta^*(a^*)-\beta^*N^*(a^*),x\rangle=\langle a^*,\beta N(x)-N\beta(x)\rangle$. Thus, \[N^*\circ \alpha^*=\alpha^*\circ N^*\;\text{ and }\;N^*\circ \beta^*=\beta^*\circ N^*\] equivalent to \[\alpha\circ N=N\circ \alpha\;\text{ and }\; \beta\circ N=N\circ \beta.\]

Then by Proposition~\ref{prop-dual-BHL}, we only need to prove that $N^*$ satisfies Eq.~(\ref{eq-Nij}) if and only if $N$ satisfies Eq.~(\ref{NBL-coalg}). We have
\begin{align*}
	&\ \langle[N^*(a^*),N^*(b^*)]_{*}+(N^*)^2([a^*,b^*]_{*})-N^*\big([N^*(a^*),b^*]_{*}+[a^*,N^*(b^*)]_{*}\big),x\rangle\\
	=&\ \langle[N^*(a^*),N^*(b^*)]_{*},x\rangle+\langle[a^*,b^*]_{*},N^2(x)\rangle-\langle[N^*(a^*),b^*]_{*}-[a^*,N^*(b^*)]_{*},N(x)\rangle\\
	=&\ \langle N^*(a^*)\otimes N^*(b^*),\Delta(x)\rangle+\langle a^*\otimes b^*,\Delta N^2(x)\rangle-\langle N^*(a^*)\otimes b^*-a^*\otimes N^*(b^*),\Delta N(x)\rangle\\
	=&\ \langle a^*\otimes b^*,(N\otimes N)\Delta(x)+\Delta(N^2(x))-(N\otimes \mathrm{id})\Delta N(x)-(\mathrm{id}\otimes N)\Delta N(x)\rangle.
\end{align*}
That is, \[[N^*(a^*),N^*(b^*)]_{*}+(N^*)^2([a^*,b^*]_{*})=N^*\big([N^*(a^*),b^*]_{*}+[a^*,N^*(b^*)]_{*}\big)\] if and only if \[(N\otimes N)\Delta(x)+\Delta(N^2(x))=(N\otimes \mathrm{id})\Delta N(x)+(\mathrm{id}\otimes N)\Delta N(x).\]
This completes the proof.
\end{proof}

\begin{defn}
A bilinear form $\mathcal{B}$ on BiHom-Lie algebra $(L,[-,-],\alpha,\beta)$ is called {\bf invariant} if
\begin{align}
	\mathcal{B}\big(\alpha(x),\alpha(y)\big)=&\ \mathcal{B}\big(x,y\big),\nonumber\\ \mathcal{B}\big(\beta(x),\beta(y)\big)=&\ \mathcal{B}\big(x,y\big),\nonumber\\
	\mathcal{B}\big([x,\alpha(y)],\beta(z)\big)=&\ \mathcal{B}\big(\alpha(x),[\beta(y),z]\big),\label{BH-B}
\end{align}
for any $x,y,z\in L$.
\end{defn}

\begin{defn}
	Let $(L,[-,-],N,\alpha,\beta)$ be a Nijenhuis BiHom-Lie algebra, and let $\mathcal{B}$ be a nondegenerate invariant bilinear form on $L$. A linear map $\widetilde{N}: L\to L$ is called an \textbf{adjoint linear map} of $N$ with respect to $\mathcal{B}$ if
	\begin{align}\label{eq-adjoint-map}
		\mathcal{B}\bigl(N(x),y\bigr)=\mathcal{B}\bigl(x,\widetilde{N}(y)\bigr)
	\end{align}
	holds for all $x,y\in L$.
\end{defn}

\begin{prop}\label{prop-star}
    Let $(L,[-,-],N,\alpha,\beta)$ be an involutive Nijenhuis BiHom-Lie algebra, and let $\mathcal{B}$ be a nondegenerate invariant bilinear form on $L$. Suppose that $\widetilde{N}: L\to L$ is the adjoint linear map of $N$ with respect to $\mathcal{B}$. Then the linear map $\widetilde{N}$ is adjoint-admissible to $(L,[-,-],N,\alpha,\beta)$.
\end{prop}

\begin{proof}
For any $x,y,z\in L$, we have
	\begin{align*}
		&\ \mathcal{B}\big(x,\widetilde{N}\big([N(y),z]\big)+[y,(\widetilde{N})^2(z)]-\widetilde{N}\big([y,\widetilde{N}(z)]\big)-[N(y),\widetilde{N}(z)]\big)\\
		=&\ \mathcal{B}\big(N(x),[N(y),z]\big)+\mathcal{B}\big(x,[y,(\widetilde{N})^2(z)]\big)-\mathcal{B}\big(N(x),[y,\widetilde{N}(z)]\big)-\mathcal{B}\big(x,[N(y),\widetilde{N}(z)]\big)\\
		=&\ \mathcal{B}\big([\alpha N(x),\alpha\beta N(y)],\beta(z)\big)+\mathcal{B}\big([\alpha(x),\alpha\beta(y)],\beta(\widetilde{N})^2(z)\big)-\mathcal{B}\big([\alpha N(x),\alpha\beta (y)],\beta\widetilde{N}(z)\big)\\
		&\ -\mathcal{B}\big([\alpha(x),\alpha\beta N(y)],\beta\widetilde{N}(z)\big)\hspace{1cm}(\text{by Eq.~(\ref{BH-B}) and }\alpha^{2}=\beta^{2}=\id_{L})\\
		=&\ \mathcal{B}\big([N\alpha(x),N\alpha\beta(y)]+N^2([\alpha(x),\alpha\beta(y)])-N\big([N\alpha(x),\alpha\beta(y)]+[\alpha(x),N\alpha\beta(y)]\big),\beta(z)\big)\\
		=&\ 0.\hspace{1cm}(\text{by Eq.~(\ref{eq-Nij})})
	\end{align*}
Thus, Eq.~(\ref{eq-adjoint-admissible}) holds, i.e., $\widetilde{N}$ is adjoint-admissible to $(L,[-,-],N,\alpha,\beta)$. This completes the proof.
\end{proof}

\subsection{Manin triples of Nijenhuis BiHom-Lie algebras}
In this subsection, we mainly study the relation between Nijenhuis BiHom-Lie bialgebras and the Manin triples of Nijenhuis BiHom-Lie algebras.

To start with, we recall the concept of standard Manin triples of Lie algebras~\cite{liebialg}.

\begin{defn}
	A (standard) {\bf Manin triple of Lie algebras} is a triple $\big((L\oplus L^*,[-,-]_{\oplus}),L,L^*\big)$, where $(L,[-,-])$ and $(L^*,[-,-]_{*})$ are Lie subalgebras of the Lie algebra $(L\oplus L^*,[-,-]_{\oplus})$ and the natural nondegenerate symmetric bilinear form $\mathcal{B}_{d}$ on $(L\oplus L^*,[-,-]_{\oplus})$ is defined by
	\begin{align}\label{eq-B-d}
		\mathcal{B}_{d}(x+a^*,y+b^*):=\langle x,b^*\rangle+\langle a^*,y\rangle,
	\end{align}
for any $x,y\in L,\;a^*,b^*\in L^*$.
\end{defn}

Now we extend this notation to Nijenhuis BiHom-Lie algebras.

\begin{defn}\label{def-Manintriple}
	Let $(L,[-,-],N,\alpha,\beta)$ and $(L^*,[-,-]_{*},S^*,\alpha^*,\beta^*)$ be Nijenhuis BiHom-Lie subalgebras of $(L\oplus L^*,[-,-]_{\oplus}^{\alpha\beta},N+S^*,\alpha_{\oplus},\beta_{\oplus})$, where $\alpha_{\oplus}=\alpha+\alpha^*,\;\beta_{\oplus}=\beta+\beta^*$ and the natural nondegenerate symmetric bilinear form $\mathcal{B}_{d}$ on $L\oplus L^*$ is defined by Eq.~(\ref{eq-B-d}). Then the triple $\big((L\oplus L^*,[-,-]_{\oplus}^{\alpha\beta},N+S^*,\alpha_{\oplus},\beta_{\oplus}),(L,[-,-],N,\alpha,\beta),(L^*,[-,-]_{*},S^*,\alpha^*,\beta^*)\big)$ is called a {\bf Manin triple of Nijenhuis BiHom-Lie algebras} with respect to the symmetric invariant bilinear form. Simply denote by $\big((L\oplus L^*,[-,-]_{\oplus}^{\alpha\beta},N+S^*,\alpha_{\oplus},\beta_{\oplus}),L,L^*\big)$.
\end{defn}

\begin{remark}
	The triple $\big((L\oplus L^*,[-,-]_{\oplus}^{\alpha\beta},N+S^*,\alpha_{\oplus},\beta_{\oplus}),L,L^*\big)$ is also called a double construction of Nijenhuis BiHom-Lie algebra with a nondegenerate symmetric invariant form associated to $(L,[-,-],N,\alpha,\beta)$ and $(L^*,[-,-]_{*},S^*,\alpha^*,\beta^*)$.
\end{remark}

Now, for the linear map $N: L\rightarrow L$, we consider the sufficient and necessary condition that $N^*$ is adjoint-admissible to involutive Nijenhuis BiHom-Lie algebra $(L^*,[-,-]_{*},S^*,\alpha^*,\beta^*)$.

\begin{lemma}
Let $(L^*,[-,-]_{*},S^*,\alpha^*,\beta^*)$ be an involutive Nijenhuis BiHom-Lie algebra. If for any $x\in L$,
\begin{align}
	(S\otimes \id)\Delta N(x)+(\id \otimes N^2)\Delta(x)=(S\otimes N)\Delta(x)+(\id\otimes N)\Delta N(x),\label{eq-N*-ad}
\end{align}
then the Nijenhuis BiHom-Lie algebra $(L^*,[-,-]_{*},S^*,\alpha^*,\beta^*)$ is $N^*$-adjoint-admissible.
\end{lemma}

\begin{proof}
By Definition~\ref{def-admissible}, for any $x\in L,\;a^*,b^*\in L^*,$ we need to prove
\[N^*\big([S^*(a^*),b^*]_{*}\big)+[a^*,(N^*)^2(b^*)]_{*}=[S^*(a^*),N^*(b^*)]_{*}+N^*\big([a^*,N^*(b^*)]_{*}\big).\]
Then, through the arbitrariness of $x$, and
\begin{align*}
	&\ \langle  N^*\big([S^*(a^*),b^*]_{*}\big)+[a^*,(N^*)^2(b^*)]_{*}-[S^*(a^*),N^*(b^*)]_{*}-N^*\big([a^*,N^*(b^*)]_{*}\big),x\rangle  \\
	=&\ \langle  [S^*(a^*),b^*]_{*},N(x)\rangle  +\langle  [a^*,(N^*)^2(b^*)]_{*}-[S^*(a^*),N^*(b^*)]_{*},x\rangle  -\langle  [a^*,N^*(b^*)]_{*},N(x)\rangle  \\
	=&\ \langle  S^*(a^*)\otimes b^*,\Delta N(x)\rangle  +\langle  a^*\otimes (N^*)^2(b^*)-S^*(a^*)\otimes N^*(b^*),\Delta(x)\rangle  -\langle  a^*\otimes N^*(b^*),\Delta N(x)\rangle  \\
	=&\ \langle  a^*\otimes b^*,(S\otimes \id)\Delta N(x)+(\id \otimes N^2)\Delta(x)-(S\otimes N)\Delta(x)-(\id\otimes N)\Delta N(x)\rangle  \\
    =&\ \langle  a^*\otimes b^*,0\rangle  \hspace{1cm}(\text{by Eq.~(\ref{eq-N*-ad})})\\
    =&\ 0.
\end{align*}
We get
\[N^*\big([S^*(a^*),b^*]_{*}\big)+[a^*,(N^*)^2(b^*)]_{*}-[S^*(a^*),N^*(b^*)]_{*}-N^*\big([a^*,N^*(b^*)]_{*}\big)=0.\]
This completes the proof.
\end{proof}

\begin{lemma}\label{lemm-**}
	Let $\big((L\oplus L^*,[-,-]_{\oplus}^{\alpha\beta},N+S^*,\alpha_{\oplus},\beta_{\oplus}),L,L^*\big)$ be a Manin triple of involutive Nijenhuis BiHom-Lie algebras. Then we have the following statements.
	\begin{enumerate}
		\item \label{lemma-1} The adjoint $\widetilde{N+S^*}$ of linear map $N+S^*$ with respect to $\mathcal{B}_{d}$ is $S+N^*$. In addition, $S+N^*$ is adjoint-admissible to Nijenhuis BiHom-Lie algebra $(L\oplus L^*,[-,-]_{\oplus}^{\alpha\beta},N+S^*,\alpha_{\oplus},\beta_{\oplus})$.
		\item \label{lemma-2} The linear map $S$ is adjoint-admissible to $(L,[-,-],N,\alpha,\beta)$.
		\item \label{lemma-3} The linear map $N^*$ is adjoint-admissible to $(L^*,[-,-]_{*},S^*,\alpha^*,\beta^*)$.
	\end{enumerate}
\end{lemma}

\begin{proof}
	\ref{lemma-1}. For any $x,y\in L$ and $a^*,b^*\in L^*$, we have
	\begin{align*}
		\mathcal{B}_{d}\big(x+a^*,\widetilde{N+S^*}(y+b^*)\big)=&\ \mathcal{B}_{d}\big((N+S^*)(x+a^*),y+b^*\big)\\
		=&\ \mathcal{B}_{d}\big(N(x)+S^*(a^*),y+b^*\big)\\
		=&\ \langle  N(x),b^*\rangle  +\langle  S^*(a^*),y\rangle  \\
		=&\ \langle  x,N^*(b^*)\rangle  +\langle  a^*,S(y)\rangle  \\
		=&\ \mathcal{B}_{d}\big(x+a^*,S(y)+N^*(b^*)\big)\\
		=&\ \mathcal{B}_{d}\big(x+a^*,(S+N^*)(y+b^*)\big).
	\end{align*}
So we get $\widetilde{N+S^*}=S+N^*$, namely, the adjoint $\widetilde{N+S^*}$ of linear map $N+S^*$ with respect to $\mathcal{B}_{d}$ is $S+N^*$. By Proposition~\ref{prop-star}, we obtain that $\widetilde{N+S^*}$ is adjoint-admissible to $(L\oplus L^*,[-,-]_{\oplus}^{\alpha\beta},N+S^*,\alpha_{\oplus},\beta_{\oplus})$. Thus, $S+N^*$ is adjoint-admissible to $L\oplus L^*$.

\ref{lemma-2}. By~\ref{lemma-1}, we know that $S+N^*$ satisfies Eq.~(\ref{eq-adjoint-admissible}), i.e.,
\begin{align}
	&\ (S+N^*)\big([(N+S^*)(x+a^*),y+b^*]_{\oplus}^{\alpha\beta}\big)+[x+a^*,(S+N^*)^2(y+b^*)]_{\oplus}^{\alpha\beta}\nonumber\\
	=&\ [(N+S^*)(x+a^*),(S+a^*)(y+b^*)]_{\oplus}^{\alpha\beta}+(S+N^*)\big([x+a,(S+N^*)(y+b^*)]_{\oplus}^{\alpha\beta}\big).\label{eq-s-N*}
\end{align}
Let $a^*=b^*=0$, we have
\[S\big([N(x),y]\big)+[x,S^2(y)]=[N(x),S(y)]+S\big([x,S(y)]\big).\]
By Definition~\ref{def-admissible}, linear map $S$ is adjoint-admissible to $(L,[-,-],N,\alpha,\beta)$.

\ref{lemma-3}. By $(L^*,[-,-]_{*},S^*,\alpha^*,\beta^*)$ is a sub-Nijenhuis BiHom-Lie algebra, we have
\[[a^*,b^*]_{\oplus}^{\alpha\beta}=[a^*,b^*]_{*}, \quad \text{for any}\, a^*,b^*\in L^*.\]
Next, taking $x=y=0$ in Eq.~(\ref{eq-s-N*}), we have
\[N^*\big([S^*(a^*),b^*]_{*}\big)+[a^*,(N^*)^2(b^*)]_{*}=[S^*(a^*),N^*(b^*)]_{*}+N^*\big([a^*,N^*(b^*)]_{*}\big).\]
By Definition~\ref{def-admissible}, linear map $N^*$ is adjoint-admissible to $(L^*,[-,-]_{*},S^*,\alpha^*,\beta^*)$. This completes the proof.
\end{proof}

Now we give the definition of Nijenhuis BiHom-Lie bialgebras as follows.

\begin{defn}\label{def-NBHL bialg}
A \textbf{Nijenhuis BiHom-Lie bialgebra} is a $7$-tuple $(L,[-,-],N,\Delta,S,\alpha,\beta)$, where $L$ is a vector space, the bracket $[-,-]: L\otimes L\rightarrow L$ is a bilinear map, and $N, S, \alpha, \beta: L\rightarrow L$ are linear maps, such that the following conditions hold:
\begin{enumerate}
    \item $(L,[-,-],\Delta,\alpha,\beta)$ is a BiHom-Lie bialgebra;
    \item $(L,[-,-],N,\alpha,\beta)$ is an involutive Nijenhuis BiHom-Lie algebra;
    \item $(L,\Delta,S,\alpha,\beta)$ is a Nijenhuis BiHom-Lie coalgebra;
    \item the linear map $S$ is adjoint-admissible to $(L,[-,-],N,\alpha,\beta)$;
    \item the linear map $N^*$ is adjoint-admissible to $(L^*,\Delta^*,S^*,\alpha^*,\beta^*)$.
\end{enumerate}
\end{defn}

By~\cite{extending}, we give the specific definition of the operation $[-,-]_{\oplus}^{\alpha\beta}: L\oplus L^*\rightarrow L\oplus L^*$ as follows.
\begin{align}
\label{eq-bigspace}
[x+a^*,\,y+b^*]_{\oplus}^{\alpha\beta}
:=&\ [x,y]
    +\mathfrak{ad}^*(a^*)\,y
    -\mathfrak{ad}^*\big(\alpha^*\beta^*(b^*)\big)\alpha^{-1}\beta(x) \nonumber \\
 &+\,[a^*,b^*]_{*}
    +\mathrm{ad}^*(x)\,b^*
    -\mathrm{ad}^*\big(\alpha\beta(y)\big)(\alpha^*)^{-1}\beta^*(a^*).
\end{align}
for any $x,y\in L,\; a^*,b^*\in L^*$, where linear map $\mathfrak{ad}^*: L^*\rightarrow \mathfrak{gl}(L)$ is a representation of Lie algebra $(L^*,[-,-]_{*})$ on representation space $L$ and
\begin{equation}\label{eq-dual}
\bigl\langle \mathfrak{ad}^*(a^*)\,x ,\, b^* \bigr\rangle
=
\bigl\langle x ,\, [a^*,b^*]_{*} \bigr\rangle, \quad \text{for any }x\in L,\; a^*,b^*\in L^*.
\end{equation}

\vspace{0.3cm}

Let $\alpha, \beta: L\rightarrow L$ be invertible linear maps. We establish the following lemma.

\begin{lemma}\label{lemma1}
Let $(L,[-,-],\alpha,\beta)$ be a BiHom-Lie algebra, and let $(L^*,[-,-]_{*},\alpha^*,\beta^*)$ be a BiHom-Lie algebra structure on the dual space $L^*$. Suppose that $\Delta: L\rightarrow L\otimes L$ is the linear dual map of $[-,-]_{*}: L^*\otimes L^*\rightarrow L^*$. Then the triple $\bigl((L\oplus L^*,[-,-]_{\oplus}^{\alpha\beta},\alpha_{\oplus},\beta_{\oplus}),\, L,\, L^*\bigr)$ is a Manin triple of BiHom-Lie algebras if and only if $(L,[-,-],\Delta,\alpha,\beta)$ is a BiHom-Lie bialgebra.
\end{lemma}

\begin{proof}
This proof relies on the transformation between BiHom-Lie and Lie structures, as well as the well-established equivalence of Lie bialgebras and Manin triples. The argument proceeds in three steps:

\noindent\textbf{Step 1:}
First, we link BiHom-Lie bialgebras on $L$ to  Lie bialgebras via the invertible maps $\alpha, \beta$.

Suppose that $(L,[-,-],\Delta,\alpha,\beta)$ is a BiHom-Lie bialgebra. Define twisted operations on $L$:
\[
\begin{aligned}
\nabla :\; & L \longrightarrow L\otimes L,\\
           & x \longmapsto \nabla(x):=(\alpha^{-1}\otimes\beta^{-1})\Delta(x),
\end{aligned}
\]
\[
\begin{aligned}
[[ -,- ]] :\; & L\otimes L \longrightarrow L,\\
              & (x,y)\longmapsto [[x,y]] := [\alpha^{-1}(x),\beta^{-1}(y)].
\end{aligned}
\]
By Proposition~\ref{Yau-Lie bialg-2}, $(L,[[ -,- ]],\nabla)$ is a Lie bialgebra.

Conversely, if $(L,[[ -,- ]],\nabla)$ is a Lie bialgebra, define BiHom-twisted operations (reversing the above twist):
\[
\begin{aligned}
\Delta :\; & L \longrightarrow L\otimes L,\\
           & x \longmapsto \Delta(x):=(\alpha\otimes\beta)\nabla(x),
\end{aligned}
\]
\[
\begin{aligned}
[-,-] :\; & L\otimes L \longrightarrow L,\\
          & (x,y)\longmapsto [x,y]:=[[\alpha(x),\beta(y)]].
\end{aligned}
\]
By Proposition~\ref{Yau-Lie bialg}, $(L,[-,-],\Delta,\alpha,\beta)$ is a BiHom-Lie bialgebra.

We conclude:
\[
(L,[-,-],\Delta,\alpha,\beta) \text{ is a BiHom-Lie bialgebra} \iff (L,[[ -,- ]],\nabla) \text{ is a Lie bialgebra}.
\]

\noindent\textbf{Step 2:}
We extend the twist to the direct sum $L \oplus L^*$. Define the twisted bracket on $L \oplus L^*$:
\[
[[u, v]]_{\oplus} := [\alpha_{\oplus}^{-1}(u), \beta_{\oplus}^{-1}(v)]_{\oplus}^{\alpha\beta}, \quad \text{for any } u, v \in L \oplus L^*,
\]
where $\alpha_{\oplus} = \alpha \oplus \alpha^*$, $\beta_{\oplus} = \beta \oplus \beta^*$.

By the same twisting logic as Step 1:
$
\big((L\oplus L^*,[-,-]_{\oplus}^{\alpha\beta},\alpha_{\oplus},\beta_{\oplus}),L,L^*\big)$  is a Manin triple of BiHom-Lie algebras if and only if $(L\oplus L^*,[[-,-]]_{\oplus},L,L^*)$ is a Manin triple of Lie algebras.

\noindent\textbf{Step 3:}
By~\cite{liebialg}, for any finite-dimensional Lie algebra $L$, the triple $\big((L\oplus L^*,[[-,-]]_{\oplus}),L,L^*\big)$ is a Manin triple of Lie algebras if and only if $(L,[[-,-]],\nabla)$ is a Lie bialgebra.

Combining Steps 1-3, we get the equivalence:
$
\big((L\oplus L^*,[-,-]_{\oplus}^{\alpha\beta},\alpha_{\oplus},\beta_{\oplus}),L,L^*\big)$ is a Manin triple of BiHom-Lie algebras if and only if $(L,[-,-],\Delta,\alpha,\beta)$ is a BiHom-Lie bialgebra.
\end{proof}

\begin{theorem}\label{th}
Let $(L,[-,-],N,\alpha,\beta)$ be an involutive Nijenhuis BiHom-Lie algebra, and let $(L^*,[-,-]_{*},S^*,\alpha^*,\beta^*)$ be a Nijenhuis BiHom-Lie algebra structure on the dual space $L^*$. Suppose that $\Delta: L\rightarrow L\otimes L$ is the linear dual map of $[-,-]_{*}: L^*\otimes L^*\rightarrow L^*$. Then the triple $\bigl((L\oplus L^*,[-,-]_{\oplus}^{\alpha\beta},N+S^*,\alpha_{\oplus},\beta_{\oplus}),\, L,\, L^*\bigr)$ is a Manin triple of Nijenhuis BiHom-Lie algebras if and only if $(L,[-,-],N,\Delta,S,\alpha,\beta)$ is a Nijenhuis BiHom-Lie bialgebra.
\end{theorem}

\begin{proof}
If $\big((L\oplus L^*,[-,-]_{\oplus}^{\alpha\beta},N+S^*,\alpha_{\oplus},\beta_{\oplus}),L,L^*\big)$ is a Manin triple of Nijenhuis BiHom-Lie algebras, then by Proposition~\ref{prop-dual-Nij},  Lemma~\ref{lemm-**}, Definition~\ref{def-NBHL bialg} and Lemma~\ref{lemma1}, we know that $(L,[-,-],N,\Delta,S,\alpha,\beta)$ is a Nijenhuis BiHom-Lie bialgebra.

Conversely, by Definition~\ref{def-Manintriple}, we only need to prove that $(L\oplus L^*,[-,-]_{\oplus}^{\alpha\beta},N+S^*,\alpha_{\oplus},\beta_{\oplus})$ is a Nijenhuis BiHom-Lie algebra.

By Lemma~\ref{lemma1}, $(L\oplus L^*,[-,-]_{\oplus}^{\alpha\beta},\alpha_{\oplus},\beta_{\oplus})$ is a BiHom-Lie algebra, so it suffices to prove
\begin{align}\label{eq-Th}
&\ [(N+S^*)(x+a^*),(N+S^*)(y+b^*)]_{\oplus}^{\alpha\beta}+(N+S^*)^2([x+a^*,y+b^*]_{\oplus}^{\alpha\beta})\nonumber\\
=&\ (N+S^*)\big([(N+S^*)(x+a^*),y+b^*]_{\oplus}^{\alpha\beta}+[x+a^*,(N+S^*)(y+b^*)]_{\oplus}^{\alpha\beta}\big),
\end{align}
for any $x,y\in L,\; a^*,b^*\in L^*$.

Indeed, we calculate
\begin{align*}
&\ \big[(N+S^*)(x+a^*),(N+S^*)(y+b^*)\big]_{\oplus}^{\alpha\beta}
      +(N+S^*)^2\big([x+a^*,y+b^*]_{\oplus}^{\alpha\beta}\big) \\
=&\ [N(x),N(y)]
   +\mathfrak{ad}^*(S^*(a^*))\,N(y)
   -\mathfrak{ad}^*\!\big(\alpha^*\beta^*S^*(b^*)\big)\alpha^{-1}\beta (N(x))\\
&+[S^*(a^*),S^*(b^*)]_{*}
   +\mathrm{ad}^*(N(x))S^*(b^*)
   -\mathrm{ad}^*(\alpha\beta N(y))(\alpha^*)^{-1}\beta^*(S^*(a^*)) \\
&+N^2\Big(
      [x,y]
      +\mathfrak{ad}^*(a^*)y
      -\mathfrak{ad}^*\!\big(\alpha^*\beta^*b^*\big)\alpha^{-1}\beta(x)
      \Big) \\
&+(S^*)^2\Big(
      [a^*,b^*]_{*}
      +\mathrm{ad}^*(x)b^*
      -\mathrm{ad}^*(\alpha\beta(y))(\alpha^*)^{-1}\beta^*(a^*)
      \Big), \hspace{1cm}(\text{by Eq.~(\ref{eq-bigspace})})
\end{align*}
and
\begin{align*}
&\ (N+S^*)\Big(
      [(N+S^*)(x+a^*),y+b^*]_{\oplus}^{\alpha\beta}
      +[x+a^*,(N+S^*)(y+b^*)]_{\oplus}^{\alpha\beta}
      \Big) \\
&\ = N\Big(
      [N(x),y]
      +\mathfrak{ad}^*(S^*(a^*))(y)
      -\mathfrak{ad}^*\!\big(\alpha^*\beta^*(b^*)\big)\alpha^{-1}\beta(N(x))
      \Big) \\
&\quad +S^*\Big(
      [S^*(a^*),b^*]_{*}
      +\mathrm{ad}^*(N(x))b^*
      -\mathrm{ad}^*(\alpha\beta(y))(\alpha^*)^{-1}\beta^*S^*(a^*)
      \Big) \\
&\quad +N\Big(
      [x,N(y)]
      +\mathfrak{ad}^*(a^*)N(y)
      -\mathfrak{ad}^*\!\big(\alpha^*\beta^*S^*(b^*)\big)\alpha^{-1}\beta(x)
      \Big) \\
&\quad +S^*\Big(
      [a^*,S^*(b^*)]_{*}
      +\mathrm{ad}^*(x)S^*(b^*)
      -\mathrm{ad}^*(\alpha\beta N(y))(\alpha^*)^{-1}\beta^*(a^*)
      \Big).\hspace{1cm}(\text{by Eq.~(\ref{eq-bigspace})})
\end{align*}

We need to prove that Eq.~(\ref{eq-Th}) holds on the basis of $L\oplus L^*$, so we divide it into the following four cases.

\textbf{Case 1.}Taking $a^*=b^*=0$, Eq.~(\ref{eq-Th}) just is Eq.~(\ref{eq-Nij}).

\textbf{Case 2.}Taking $x=y=0$, Eq.~(\ref{eq-Th}) holds if $S^*$ satisfies Eq.~(\ref{eq-Nij}) on Lie algebra $(L^*,[-,-]_{*})$.

\textbf{Case 3.}Taking $a^*=y=0$, we need some additional minor calculations as follows.

For any $x\in L,\;b^*,c^*\in L^*$, we have
\begin{align*}
&\ \Big\langle
   \mathfrak{ad}^*\!\big(S^*\alpha^*\beta^*(b^*)\big)\,N\alpha^{-1}\beta(x)
   +N^2\Big(\mathfrak{ad}^*\!\big(\alpha^*\beta^*(b^*)\big)\alpha^{-1}\beta(x)\Big) \\
&-N\Big(\mathfrak{ad}^*\!\big(\alpha^*\beta^*(b^*)\big)N\alpha^{-1}\beta(x)\Big)
   -N\Big(\mathfrak{ad}^*\!\big(S^*\alpha^*\beta^*(b^*)\big)\alpha^{-1}\beta(x)\Big),
   c^*
 \Big\rangle \\
=&\ -\big\langle
      N(\alpha^{-1}\beta(x)),
      [S^*\alpha^*\beta^*(b^*),c^*]_{*}
    \big\rangle
   -\big\langle
      \alpha^{-1}\beta(x),
      [\alpha^*\beta^*(b^*),(N^*)^2(c^*)]_{*}
    \big\rangle\\
&\ +\big\langle
      N\alpha^{-1}\beta(x),
      [\alpha^*\beta^*(b^*),N^*(c^*)]_{*}
    \big\rangle
   +\big\langle
      \alpha^{-1}\beta(x),
      [S^*\alpha^*\beta^*(b^*),N^*(c^*)]_{*}
    \big\rangle \hspace{1cm}(\text{by Eq.~(\ref{eq-dual})}) \\
=&\ -\Big\langle
      \alpha^{-1}\beta(x),
      N^*\big[ S^*\alpha^*\beta^*(b^*),c^* \big]_{*}
      +[\alpha^*\beta^*(b^*),(N^*)^2(c^*)]_{*} \\
&\ -N^*\big([\alpha^*\beta^*(b^*),N^*(c^*)]_{*}\big)
      -[S^*\alpha^*\beta^*(b^*),N^*(c^*)]_{*}
   \Big\rangle .
\end{align*}
and
\begin{align*}
&\ \bigl\langle
    \mathrm{ad}^*(N(x))S^*(b^*)
    +(S^*)^2 \mathrm{ad}^*(x)b^*
    -S^*\mathrm{ad}^*(N(x))b^*
    -S^*\mathrm{ad}^*(x)S^*(b^*),
    \, z
\bigr\rangle \\[0.5em]
=&\
\bigl\langle \mathrm{ad}^*(x)b^*, S^2(z) \bigr\rangle
-\bigl\langle S^*(b^*), [N(x),z] \bigr\rangle
-\bigl\langle \mathrm{ad}^*(N(x))b^*, S(z) \bigr\rangle
-\bigl\langle \mathrm{ad}^*(x)S^*(b^*), S(z) \bigr\rangle \quad\text{(by Eq.~(\ref{eq-dual}))}\\[0.5em]
=&\
\bigl\langle
    b^*,
    \, -S([N(x),z])-[x,S^2(z)]+[N(x),S(z)]
    +S([x,S(z)])
\bigr\rangle .
\end{align*}

Hence, Eq.~(\ref{eq-Th}) holds if Eq.~(\ref{eq-adjoint-admissible}) and Eq.~(\ref{eq-N*-ad}) hold.

\textbf{Case 4.}Taking $x=b^*=0$, the argument is exactly the same as in Case 3.

Therefore, we conclude that Eq.~(\ref{eq-Th}) holds if
\begin{enumerate}
	\item  $(L,[-,-], N, \alpha,\beta)$ is a Nijenhuis BiHom-Lie algebra;
	\item $(L,\Delta, S, \alpha,\beta)$ is a Nijenhuis BiHom-Lie coalgebra;
	\item linear map $S$ is adjoint-admissible to $(L,[-,-],N,\alpha,\beta)$;
	\item linear map $N^*$ is adjoint-admissible to $(L^*,[-,-]_{*},S^*,\alpha^*,\beta^*)$.
\end{enumerate}
By Definition~\ref{def-NBHL bialg}, we complete this proof.
\end{proof}

\begin{coro}\label{coro-nijenhuis-lie}
Let $(L,[-,-],N)$ be a Nijenhuis Lie algebra, and let $(L^*,[-,-]_{*},S^*)$ be a Nijenhuis Lie algebra structure on the dual space $L^*$. Suppose that $\Delta: L\rightarrow L\otimes L$ is the linear dual map of $[-,-]_{*}: L^*\otimes L^*\rightarrow L^*$. Then the triple $\bigl((L\oplus L^*,[-,-]_{\oplus},N+S^*),\, L,\, L^*\bigr)$ is a Manin triple of Nijenhuis Lie algebras if and only if $(L,[-,-],N,\Delta,S)$ is a Nijenhuis Lie bialgebra.
\end{coro}

\begin{proof}
	This proof follows from Theorem~\ref{th} by taking $\alpha=\beta=\mathrm{id}_{L},\;\alpha^*=\beta^*=\mathrm{id}_{L^*}$.
\end{proof}

\subsection{Matched pairs of Nijenhuis BiHom-Lie algebras}
In this subsection, we mainly study the relation between Nijenhuis BiHom-Lie bialgebras and the matched pairs of Nijenhuis BiHom-Lie algebras.

In order to define the matched pairs of Nijenhuis BiHom-Lie algebras, we first introduce the following results.
\begin{theorem}\label{th-mp}
	Let $(L,[-,-],N,\alpha,\beta)$ and $(V,[-,-]_{V},S,p,q)$ be two involution Nijenhuis BiHom-Lie algebras. Suppose that $\rho: L\rightarrow \mathfrak{gl}(V)$ and $h: V\rightarrow \mathfrak{gl}(L)$ are two linear maps such that $(V,\rho,S,p,q)$ is a representation of $L$, and $(L,h,N,\alpha,\beta)$ is a representation of $V$. For all $x,y\in L$ and $a,b,c,z\in V$, the following identities hold:
{\small
\begin{equation}
[y,\,h(q(c))\alpha(x)] - [x,\,h(q(c))\alpha(y)] - h\big(\rho(\alpha(y))q(c)\big)\alpha\beta(x) + h\big(\rho(\alpha(x))q(c)\big)\alpha\beta(y)
 +h(c)\big([\beta(x),\alpha(y)]\big)
=0,\label{M-3}
\end{equation}
\begin{equation}
[b,\,\rho(\beta(z))p(a)]_{V}-[a,\,\rho(\beta(z))p(b)]_{V}-\rho \big(h(p(b))\beta(z)\big)pq(a) +\rho\big(h(p(a))\beta(z)\big)pq(b)
 +\rho(z)\big([q(a),p(b)]_{V}\big)
=0,\label{M-4}
\end{equation}
}
Then there is a Nijenhuis BiHom-Lie algebra structure $(L\oplus V,[-,-]_{\oplus}^{\alpha\beta},N+S,\alpha+p,\beta+q)$ on the direct sum $L\oplus V$ of the underlying vector space of $L$ and $V$ given by
\begin{align}
	&\hspace{2cm}(N+S)(x+a):=N(x)+S(a),\nonumber\\
	&(\alpha+p)(x+a):=\alpha(x)+p(a),\quad (\beta+q)(x+a):=\beta(x)+q(a),\nonumber\\
	&[x+a,y+b]_{\oplus}^{\alpha\beta}:=[x,y]+h(a)y-h(pq(b))\alpha\beta(x)+[a,b]_{V}+\rho(x)b
-\rho(\alpha\beta(y))pq(a),\label{eq-oplus}
\end{align}
for any $x,y\in L,\, a,b\in V$. We denote this Nijenhuis BiHom-Lie algebra by $L\Join_{h,\alpha,\beta}^{\rho,p,q}V$ or simply $L\Join V$.

Conversely, if $(L\oplus V,[-,-]_{\oplus}^{\alpha\beta},N+S,\alpha+p,\beta+q)$ is a Nijenhuis BiHom-Lie algebra, then $(V,\rho,S,p,q)$ is a representation of $L$, $(L,h,N,\alpha,\beta)$ is a representation of $V$ and Eqs.~(\ref{M-3})-(\ref{M-4}) hold.
\end{theorem}

\begin{proof}
By Definition~\ref{nijenhuis-bihom-lie}, we know that the $5$-tuple $(L\oplus V,[-,-]_{\oplus}^{\alpha\beta},N+S,\alpha+p,\beta+q)$ is a Nijenhuis BiHom-Lie algebra if and only if for any $x,y,z\in L,\,a,b,c\in V$,
\begin{align}
\label{BHL-LV-1}
[(\beta+q)(x+a),(\alpha+p)(y+b)]_{\oplus}^{\alpha\beta}
+[(\beta+q)(y+b),(\alpha+p)(x+a)]_{\oplus}^{\alpha\beta}
=0,
\end{align}
\begin{align}
\label{BHL-LV-2}
&[(\beta+q)^2(x+a),
  [(\beta+q)(y+b),(\alpha+p)(z+c)]_{\oplus}^{\alpha\beta}]_{\oplus}^{\alpha\beta} \nonumber\\
&\quad
+[(\beta+q)^2(y+b),
  [(\beta+q)(z+c),(\alpha+p)(x+a)]_{\oplus}^{\alpha\beta}]_{\oplus}^{\alpha\beta} \nonumber\\
&\quad
+[(\beta+q)^2(z+c),
  [(\beta+q)(x+a),(\alpha+p)(y+b)]_{\oplus}^{\alpha\beta}]_{\oplus}^{\alpha\beta}
=0,
\end{align}
\begin{align}
\label{BHL-LV-3}
&\quad [(N+S)(x+a),(N+S)(y+b)]_{\oplus}^{\alpha\beta}
+(N+S)^2\big([x+a,y+b]_{\oplus}^{\alpha\beta}\big) \nonumber\\
&
=(N+S)\Big(
   [(N+S)(x+a),y+b]_{\oplus}^{\alpha\beta}
  +[x+a,(N+S)(y+b)]_{\oplus}^{\alpha\beta}
 \Big) .
\end{align}
Since $L\oplus V$ is the direct sum of $L$ and $V$ as vector spaces, Eqs.~(\ref{BHL-LV-1})-(\ref{BHL-LV-3}) hold if and only if they are satisfied for all generators of $L\oplus V$, i.e., for elements of the form $x+0$ (with $x\in L$) and $0+a$ (with $a\in V$).
By the definition of the bracket operation $[-,-]_{\oplus}^{\alpha\beta}$ (cf. Eq.~(\ref{eq-oplus})), if $(L,[-,-],N,\alpha,\beta)$ and $(V,[-,-]_{V},S,p,q)$ are involutive Nijenhuis BiHom-Lie algebras, then Eq.~(\ref{BHL-LV-1}) holds trivially. For Eqs.~(\ref{BHL-LV-2})-(\ref{BHL-LV-3}), the validity relies on specific generator combinations, which are summarized in the following table:

\begin{table}[htbp]
    \centering
    \begin{tabular}{|p{0.60\textwidth}|p{0.35\textwidth}|}
        \hline
        \textbf{Conditions} & \textbf{Derived Results} \\
        \hline
        Eqs.~(\ref{M-3})-(\ref{M-4}) hold; &\\
        $(V,\rho,p,q)$ is a representation of $(L,[-,-],\alpha,\beta)$; & Eq.~(\ref{BHL-LV-2}) holds. \\
        $(L,h,\alpha,\beta)$ is a representation of $(V,[-,-]_{V},p,q).$ & \\
        \hline
        $(L,h,N,\alpha,\beta)$ is a representation of $(V,[-,-]_{V},S,p,q)$; & Eq.~(\ref{BHL-LV-3}) holds.\\
        $(V,\rho,S,p,q)$ is a representation of $(L,[-,-],N,\alpha,\beta).$ & \\
        \hline
    \end{tabular}
\end{table}

\noindent Conversely, through directly computing, we obtain the following table:
\begin{table}[ht]
\centering
\renewcommand{\arraystretch}{1.25}
\begin{tabular}{|p{0.55\textwidth}|p{0.40\textwidth}|}
\hline
\textbf{Conditions (Special values substitution)} & \textbf{Derived Results} \\
\hline

\begin{minipage}[t]{\linewidth}
Eq.~(\ref{BHL-LV-2}) holds for \(x+0,\ y+0,\ 0+a\) on \(L\).
\end{minipage}
&
\begin{minipage}[t]{\linewidth}
Eq.~(\ref{M-3}) holds.
\end{minipage}
\\
\hline

\begin{minipage}[t]{\linewidth}
Eq.~(\ref{BHL-LV-2}) holds for \(x+0,\ 0+a,\ 0+b\) on \(V\).
\end{minipage}
&
\begin{minipage}[t]{\linewidth}
Eq.~(\ref{M-4}) holds.
\end{minipage}
\\
\hline

\begin{minipage}[t]{\linewidth}
   Eq.~(\ref{BHL-LV-2}) holds for \(x+0,\ 0+a,\ 0+b\) on \(V\);\\
   Eq.~(\ref{BHL-LV-3}) holds for \(x+0,\ 0+b\) on \(L\).
\end{minipage}
&
\begin{minipage}[t]{\linewidth}
\((L,h,N,\alpha,\beta)\) is a representation of \(V\).
\end{minipage}
\\
\hline

\begin{minipage}[t]{\linewidth}
 Eq.~(\ref{BHL-LV-2}) holds for \(x+0,\ y+0,\ 0+a\) on \(V\);\\
 Eq.~(\ref{BHL-LV-3}) holds for \(x+0,\ 0+a\) on \(L\).
\end{minipage}
&
\begin{minipage}[t]{\linewidth}
\((V,\rho,S,p,q)\) is a representation of \(L\).
\end{minipage}
\\
\hline
\end{tabular}
\end{table}

\noindent
This completes the proof.
\end{proof}

\begin{coro}
	Let $(L,[-,-],\alpha,\beta)$ and $(V,[-,-]_{V},p,q)$ be two BiHom-Lie algebras. If linear maps $\rho: L\rightarrow \mathfrak{gl}(V)$ and $h: V\rightarrow \mathfrak{gl}(L)$ satisfy Eqs.~(\ref{M-3})-(\ref{M-4}) and
	\begin{enumerate}
		\item $(V,\rho,p,q)$ is a representation of BiHom-Lie algebra $L$ on $V$,
		\item $(L,h,\alpha,\beta)$ is a representation of BiHom-Lie algebra $V$ on $L$,
	\end{enumerate}
then $(L\oplus V,[-,-]_{\oplus}^{\alpha\beta},\alpha+p,\beta+q)$ is a BiHom-Lie algebra.
\end{coro}

In this case, we say that $(L,V,h,\rho)$ is a matched pair of BiHom-Lie algebras  
if $(V,\rho,p,q)$, $(L,h,\alpha,\beta)$ are representations of BiHom-Lie algebras $(L,[-,-],\alpha,\beta)$ and $(V,[-,-]_{V},p,q)$ respectively and Eqs.~(\ref{M-3})-(\ref{M-4}) hold.

\begin{defn}\label{MPOBHL}
A {\bf matched pair of Nijenhuis BiHom-Lie algebras $(L,[-,-],N,\alpha,\beta)$ and $(V,[-,-]_{V},S,p,q)$} is a $6$-tuple $(L,V,N,S,h,\rho)$, where
    \begin{enumerate}
    \item $(L,V,h,\rho)$ is a matched pair of BiHom-Lie algebras $(L,[-,-],\alpha,\beta)$ and $(V,[-,-]_{V},p,q)$;
        \item $(V,\rho,S,p,q)$ is a representation of the Nijenhuis BiHom-Lie algebra $(L,[-,-],N,\alpha,\beta)$;
        \item $(L,h,N,\alpha,\beta)$ is a representation of the Nijenhuis BiHom-Lie algebra $(V,[-,-]_{V},S,p,q)$.
    \end{enumerate}
   \end{defn}

\begin{lemma}\label{lemma-mp}
    Let $(L, [-,-], \alpha, \beta)$ be a BiHom-Lie algebra, and suppose that its dual space $L^*$ admits a BiHom-Lie algebra structure $(L^*, [-,-]_*, \alpha^*, \beta^*)$.
    Then the $5$-tuple $(L, [-,-], \Delta, \alpha, \beta)$ is a BiHom-Lie bialgebra if and only if $(L,L^*,\mathfrak{ad}^*,\mathrm{ad}^*)$ is a matched pair of BiHom-Lie algebras.
\end{lemma}

\begin{proof}
	This is straightforward.
\end{proof}

\begin{theorem}\label{th-23}
    Let $(L, [-,-], N, \alpha, \beta)$ be an involutive Nijenhuis BiHom-Lie algebra, and suppose that its dual space $L^*$ carries a Nijenhuis BiHom-Lie algebra structure $(L^*, [-,-]_*, S^*, \alpha^*, \beta^*)$.
    Denote by $\Delta^* = [-,-]_*$, and let $\mathfrak{ad}^*: L^* \to \mathfrak{gl}(L)$ and $\mathrm{ad}^*: L \to \mathfrak{gl}(L^*)$ be the corresponding coadjoint representations.
    Then $(L, [-,-], \Delta, N, S, \alpha, \beta)$ is a Nijenhuis BiHom-Lie bialgebra if and only if $(L, L^*, N, S^*, \mathfrak{ad}^*, \mathrm{ad}^*)$ is a matched pair of Nijenhuis BiHom-Lie algebras.
\end{theorem}

\begin{proof}
	Suppose that $(L,[-,-],\Delta,N,S,\alpha,\beta)$ is an involutive Nijenhuis BiHom-Lie bialgebra. By Definition~\ref{def-NBHL bialg}, it implies that $(L^*,\mathrm{ad}^*,S^*,\alpha^*,\beta^*)$ is a representation of $(L,[-,-],N,\alpha,\beta)$, and that $(L,\mathfrak{ad}^*,N,\alpha,\beta)$ is a representation of $(L^*,[-,-]_{*},S^*,\alpha^*,\beta^*)$. Moreover, by Lemma~\ref{lemma-mp}, the tuple $(L,L^*,\mathfrak{ad}^*,\mathrm{ad}^*)$ is a matched pair of BiHom-Lie algebras if and only if $(L,[-,-],\Delta,\alpha,\beta)$ is a BiHom-Lie bialgebra.
    Combining the above results with Definition~\ref{MPOBHL}, we conclude that $(L,L^*,N,S^*,\mathfrak{ad}^*,\mathrm{ad}^*)$ is a matched pair of Nijenhuis BiHom-Lie algebras.

	Conversely, if $(L,L^*,N,S^*,\mathfrak{ad}^*,\mathrm{ad}^*)$ is a matched pair of Nijenhuis BiHom-Lie algebras, then $(L,L^*,\mathfrak{ad},\mathrm{ad}^*)$ is a matched pair of BiHom-Lie algebras and $S,\, N^*$ are adjoint-admissible to $(L,[-,-],N,\alpha,\beta)$ and $(L^*,\Delta^*,S^*,\alpha^*,\beta^*)$ respectively. Thus, $(L,[-,-],\Delta,N,S,\alpha,\beta)$ is a Nijenhuis BiHom-Lie bialgebra by  Proposition~\ref{prop-dual-Nij},  Definition~\ref{def-NBHL bialg} and Lemma~\ref{lemma-mp}.
\end{proof}

\begin{theorem}\label{th-123}
	Let $(L,[-,-],N,\alpha,\beta)$ and $(L^{*},[-,-]_{*},S^{*},\alpha^{*},\beta^{*})$ be two Nijenhuis BiHom-Lie algebras. Then the following  conditions are equivalent.
	\begin{enumerate}
		\item $((L\oplus L^{*},[-,-]_{\oplus}^{\alpha\beta},N+S^{*},\alpha_{\oplus},\beta_{\oplus}),L,L^{*})$ is a Manin triple of Nijenhuis BiHom-Lie algebras.
		\item $(L,[-,-],\Delta,N,S,\alpha,\beta)$ is a Nijenhuis BiHom-Lie bialgebra, where $\Delta: L\rightarrow L\otimes L$ is the linear dual of $[-,-]_{*}: L^{*}\otimes L^{*}\rightarrow L^{*}$.
		\item $(L,L^{*},N,S^{*},\mathfrak{ad}^{*},\mathrm{ad}^{*})$ is a matched pair of Nijenhuis BiHom-Lie algebras.
	\end{enumerate}
\end{theorem}
\begin{proof}
	It follows from Theorem~\ref{th} and Theorem~\ref{th-23}.
\end{proof}

\subsection{Nijenhuis BiHom-Lie bialgebras and SN-CBHYBEs}
In this subsection, we mainly study the corresponding Yang-Baxter equations for  (Nijenhuis) BiHom-Lie bialgebras. We almost suppose that the BiHom-Lie algebra $(L,[-,-],\alpha,\beta)$ is an involutive BiHom-Lie algebra with invertible linear map $\alpha$ until the end of this section.

To find the BiHom-version of classical Yang-Baxter equations, we first consider the equivalent conditions for a BiHom-Lie algebra to be a BiHom-Lie bialgebra.

\begin{lemma}\label{lemm-equiv}
	Let $(L,[-,-],\alpha,\beta)$ be a BiHom-Lie algebra and $r\in L\otimes L$ with
		\[r=-\tau(r),\quad (\alpha\otimes\alpha)r=r, \quad (\beta\otimes\beta)r=r.\]  Let $\Delta_{r}: L\rightarrow L\otimes L$ be an operation defined by  \begin{align}\label{eq-coboundary}
		\Delta_{r}(x):=\mathrm{ad}_{\alpha^{-1}\beta(x)}^{(2)}(r)=(\mathrm{ad}_{x}\otimes\beta+\beta\otimes \mathrm{ad}_{\alpha^{-1}\beta(x)})r,\qquad \text{ for any } x\in L.
	\end{align} Then
	\begin{enumerate}
		\item\label{lemm426-1} $\Delta_{r}$ satisfies Eq.~(\ref{eq-BHLiebialg}), i.e., $\Delta_{r}$ is a 1-cocycle of BiHom-Lie algebras;
		\item\label{lemm426-2} $(L,\Delta_{r},\alpha,\beta)$ is a BiHom-Lie coalgebra in Definition~\ref{defn-BHLco} if and only if
		\begin{equation}\label{eq-anti-r}
			\big(\mathrm{ad}_{\beta(x)}\alpha\otimes \mathrm{id} +\mathrm{id} \otimes \mathrm{ad}_{\beta(x)}\alpha\big)(r+\tau(r))=0,
			\qquad \text{for any }x\in L,
		\end{equation}
		and
		\begin{equation}\label{eq-coJacobi-YBE}
			(\mathrm{ad}_{x}\otimes \beta\otimes \beta+\beta\otimes \mathrm{ad}_{x}\otimes \beta+\beta\otimes\beta\otimes \mathrm{ad}_{x})\big([r_{12},r_{13}]_{\alpha\beta}+[r_{12},r_{23}]_{\alpha\beta}+[r_{13},r_{23}]_{\alpha\beta}\big)=0,
		\end{equation}
		where for $r_{1}\otimes r_{2}=r=\bar{r}=\bar{r}_{1}\otimes \bar{r}_{2}\in L\otimes L$,
		\[[r_{12},r_{13}]_{\alpha\beta}:=[r_{1},\alpha^{-1}\beta(\bar{r}_{1})]\otimes \beta(r_{2})\otimes\beta(\bar{r}_{2}),\]
		\[[r_{12},r_{23}]_{\alpha\beta}:=\beta(r_{1})\otimes [r_{2},\alpha^{-1}\beta(\bar{r}_{1})]\otimes \beta(\bar{r}_{2}),\]
		\[[r_{13},r_{23}]_{\alpha\beta}:=\beta(r_{1})\otimes \beta(\bar{r}_{1})\otimes [r_{2},\alpha^{-1}\beta(\bar{r}_{2})].\]
	\end{enumerate}
\end{lemma}
\begin{proof}
\ref{lemm426-1}. By~Proposition~\ref{BHL_repre-ad2}, we know that $(L\otimes L, \mathrm{ad}^{(2)},\alpha\otimes\alpha,\beta\otimes \beta)$ is a representation of BiHom-Lie algebra $(L,[-,-],\alpha,\beta)$, then for any $x,y\in L$, we have
\begin{align*}
	\Delta_{r}([\alpha^{-1}\beta(x),y])=&\ \mathrm{ad}^{(2)}_{\alpha^{-1}\beta([\alpha^{-1}\beta(x),y])}(r)\\
	=&\ \mathrm{ad}^{(2)}_{[x,\alpha^{-1}\beta(y)]}(\beta\otimes\beta)(r)\hspace{1cm}(\text{by }(\beta\otimes\beta)r=r \text{ and }\alpha^{2}=\beta^{2}=\mathrm{id}_{L})\\
	=&\ \mathrm{ad}^{(2)}_{\alpha(x)}\mathrm{ad}^{(2)}_{\alpha^{-1}\beta(y)}(r)-\mathrm{ad}^{(2)}_{\alpha(y)}\mathrm{ad}^{(2)}_{\alpha^{-1}\beta(x)}(r)\hspace{1cm}(\text{by Eq.~(\ref{LLad2-repre})~ and }\alpha^{2}=\beta^{2}=\mathrm{id}_{L})\\
	=&\ \mathrm{ad}^{(2)}_{\alpha(x)}\Delta_{r}(y)-\mathrm{ad}^{(2)}_{\alpha(y)}\Delta_{r}(x)\hspace{1cm}(\text{by Eq.~(\ref{eq-coboundary})})\\
	=&\ (\mathrm{ad}_{\beta(x)}\otimes \beta+\beta\otimes \mathrm{ad}_{\alpha^{-1}\beta^2(x)})\Delta_{r}(y)-(\mathrm{ad}_{\beta(y)}\otimes\beta+\beta\otimes \mathrm{ad}_{\alpha^{-1}\beta^2(y)})\Delta_{r}(x).
\end{align*}
Thus, $\Delta_{r}$ is a 1-cocycle of BiHom-Lie algebras.

\ref{lemm426-2}. By the definition of $\Delta_{r}$, we have
\[\Delta_{r}\circ \alpha=(\alpha\otimes\alpha)\circ\Delta_{r},\quad \Delta_{r}\circ\beta=(\beta\otimes\beta)\circ \Delta_{r}.\]
So we only need to prove that Eqs.~(\ref{eq-anti})-(\ref{eq-cojacobi}) hold for $(L,\Delta_{r},\alpha,\beta)$ if and only if Eqs.~(\ref{eq-anti-r})-(\ref{eq-coJacobi-YBE}) hold.
For any $x\in L$, we first consider the Eq.~(\ref{eq-anti}).
\begin{align*}
	&\ (\beta\otimes\alpha)\circ\Delta_{r}(x)+\tau(\beta\otimes\alpha)\circ \Delta_{r}(x)\\
	=&\ (\beta\otimes\alpha)(\mathrm{ad}_{x}\otimes \beta+\beta\otimes \mathrm{ad}_{\alpha^{-1}\beta(x)})r+\tau(\beta\otimes\alpha)
	(\mathrm{ad}_{x}\otimes \beta+\beta\otimes \mathrm{ad}_{\alpha^{-1}\beta(x)})r\\
	=&\ [\beta(x),\beta(r_{1})]\otimes \alpha\beta(r_{2})+\beta^{2}(r_{1})\otimes [\beta(x),\alpha(r_{2})]+\alpha\beta(r_{2})\otimes [\beta(x),\beta(r_{1})]+[\beta(x),\alpha(r_{2})]\otimes\beta^{2}(r_{1})\\
	=&\ \mathrm{ad}_{\beta(x)}\alpha(r_{1})\otimes r_{2}+r_{1}\otimes \mathrm{ad}_{\beta(x)}\alpha(r_{2})+r_{2}\otimes \mathrm{ad}_{\beta(x)}\alpha(r_{1})+\mathrm{ad}_{\beta(x)}\alpha(r_{2})\otimes r_{1}\\
	&\hspace{1cm}(\text{by }(\alpha\otimes \alpha)r=r,\;(\beta\otimes\beta)r=r, \text{ and } \alpha^{2}=\beta^{2}=\mathrm{id}_{L})\\
	=&\ (\mathrm{ad}_{\beta(x)}\alpha\otimes \mathrm{id})(r_{1}\otimes r_{2}+r_{2}\otimes r_{1})+(\mathrm{id}\otimes \mathrm{ad}_{\beta(x)}\alpha)(r_{1}\otimes r_{2}+r_{2}\otimes r_{1})\\
	=&\ (\mathrm{ad}_{\beta(x)}\alpha\otimes \mathrm{id}+\mathrm{id}\otimes \mathrm{ad}_{\beta(x)}\alpha)(r+\tau(r)).
\end{align*}
Thus, the Eq.~(\ref{eq-anti}) holds for $(L,\Delta_{r},\alpha,\beta)$ if and only if Eq.~(\ref{eq-anti-r}) holds. Next, through direct calculation, we derive
\begin{align*}
	&\ (\mathrm{id}\otimes \beta\otimes \alpha)(\beta^{2}\otimes \Delta_{r})\Delta_{r}(x)\\
	=&\ [x,r_{1}]\otimes \beta \mathrm{ad}_{\beta(r_{2})}(\bar{r}_{1})\otimes \alpha\beta(\bar{r}_{2})+[x,r_{1}]\otimes \beta^{2}(\bar{r}_{1})\otimes \alpha \mathrm{ad}_{\alpha^{-1}\beta^{2}(r_{2})}(\bar{r}_{2})\\
	&+\beta(r_{1})\otimes \beta \mathrm{ad}_{[\alpha^{-1}\beta(x),r_{2}]}(\bar{r}_{1})\otimes \alpha\beta(\bar{r}_{2})+\beta(r_{1})\otimes \beta^{2}(\bar{r}_{1})\otimes \alpha \mathrm{ad}_{[\alpha^{-2}\beta^{2}(x),\alpha^{-1}\beta(r_{2})]}(\bar{r}_{2}),
\end{align*}
so, for Eq.~(\ref{eq-cojacobi}), we have
\begin{align*}
	&\ (\mathrm{id}\otimes \beta\otimes \alpha)(\beta^{2}\otimes \Delta_{r})\Delta_{r}(x)+(\tau\otimes \mathrm{id})(\mathrm{id}\otimes \tau)(\mathrm{id}\otimes \beta\otimes \alpha)(\beta^{2}\otimes \Delta_{r})\Delta_{r}(x)\\
	&+(\mathrm{id}\otimes \tau)(\tau\otimes \mathrm{id})(\mathrm{id}\otimes \beta\otimes \alpha)(\beta^{2}\otimes \Delta_{r})\Delta_{r}(x)\\
	=&\ [x,r_{1}]\otimes [\beta^{2}(r_{2}),\beta(\bar{r}_{1})]\otimes \alpha\beta(\bar{r}_{2})+[x,r_{1}]\otimes \beta^{2}(\bar{r}_{1})\otimes [\beta^{2}(r_{2}),\alpha(\bar{r}_{2})]\\
	&+\beta(r_{1})\otimes [[\alpha^{-1}\beta^{2}(x),\beta(r_{2})],\beta(\bar{r}_{1})]\otimes \alpha\beta(\bar{r}_{2})+\beta(r_{1})\otimes \beta^{2}(\bar{r}_{1})\otimes [[\alpha^{-1}\beta^{2}(x),\beta(r_{2})],\alpha(\bar{r}_{2})]\\
	&+\alpha\beta(\bar{r}_{2})\otimes [x,r_{1}]\otimes [\beta^{2}(r_{2}),\beta(\bar{r}_{1})]+[\beta^{2}(r_{2}),\alpha(\bar{r}_{2})]\otimes [x,r_{1}]\otimes \beta^{2}(\bar{r}_{1})\\
	&+\alpha\beta(\bar{r}_{2})\otimes \beta(r_{1})\otimes [[\alpha^{-1}\beta^{2}(x),\beta(r_{2})],\beta(\bar{r}_{1})]
	+[[\alpha^{-1}\beta^{2}(x),\beta(r_{2})],\alpha(\bar{r}_{2})]\otimes \beta(r_{1})\otimes \beta^{2}(\bar{r}_{1})\\
	&+[\beta^{2}(r_{2}),\beta(\bar{r}_{1})]\otimes \alpha\beta(\bar{r}_{2})\otimes [x,r_{1}]+\beta^{2}(\bar{r}_{1})\otimes [\beta^{2}(r_{2}),\alpha(\bar{r}_{2})]\otimes [x,r_{1}]\\
	&+[[\alpha^{-1}\beta^{2}(x),\beta(r_{2})],\beta(\bar{r}_{1})]\otimes \alpha\beta(\bar{r}_{2})\otimes \beta(r_{1})+\beta^{2}(\bar{r}_{1})\otimes [[\alpha^{-1}\beta^{2}(x),\beta(r_{2})],\alpha(\bar{r}_{2})]\otimes \beta(r_{1})\\
	=&\ -[\beta(\bar{r}_{2}),[\beta(x),\alpha(r_{2})]]\otimes \beta(r_{1})\otimes \beta^{2}(\bar{r}_{1})-[\alpha^{-1}\beta^{2}(\bar{r}_{1}),[\beta(x),\alpha(r_{2})]]\otimes \alpha\beta(\bar{r}_{2})\otimes \beta(r_{1})\\
	&-\beta(r_{1})\otimes [\alpha^{-1}\beta^{2}(\bar{r}_{1}),[\beta(x),\alpha(r_{2})]]\otimes \alpha\beta(\bar{r}_{2})-\beta^{2}(\bar{r}_{1})\otimes [\beta(\bar{r}_{2}),[\beta(x),\alpha(r_{2})]]\otimes\beta(r_{1})\\
	&-\beta(r_{1})\otimes \beta^{2}(\bar{r}_{1})\otimes [\beta(\bar{r}_{2}),[\beta(x),\alpha(r_{2})]]-\alpha\beta(\bar{r}_{2})\otimes \beta(r_{1})\otimes [\alpha^{-1}\beta^{2}(\bar{r}_{1}),[\beta(x),\alpha(r_{2})]]\\
	&+[x,r_{1}]\otimes [r_{2},\beta(\bar{r}_{1})]\otimes \alpha\beta(\bar{r}_{2})+[r_{2},\alpha(\bar{r}_{2})]\otimes [x,r_{1}]\otimes \bar{r}_{1}+[x,r_{1}]\otimes \bar{r}_{1}\otimes [r_{2},\alpha(\bar{r}_{2})]\\
	&+[r_{2},\beta(\bar{r}_{1})]\otimes \alpha\beta(\bar{r}_{2})\otimes [x,r_{1}]+\alpha\beta(\bar{r}_{2})\otimes [x,r_{1}]\otimes [r_{2},\beta(\bar{r}_{1})]+\bar{r}_{1}\otimes [r_{2},\alpha(\bar{r}_{2})]\otimes [x,r_{1}]\\
	&\hspace{1cm}(\text{by Eq.~(\ref{eq-BH-anti}) and  }\alpha^{2}=\beta^{2}=\id_{L})\\
	=&\ \mathrm{ad}_{x}([r_{1},\alpha\beta(\bar{r}_{1})])\otimes r_{2}\otimes \bar{r}_{2}+r_{1}\otimes \mathrm{ad}_{x}([r_{2},\alpha\beta(\bar{r}_{1})])\otimes \bar{r}_{2}+r_{1}\otimes \bar{r}_{1}\otimes \mathrm{ad}_{x}([r_{2},\alpha\beta(\bar{r}_{2})])\\
	&+ \mathrm{ad}_{x}(r_{1})\otimes [r_{2},\alpha^{-1}(\bar{r}_{1})]\otimes \bar{r}_{2}+[r_{1},\alpha(\bar{r}_{1})]\otimes \mathrm{ad}_{x}(r_{2})\otimes \bar{r}_{2}+\mathrm{ad}_{x}(r_{1})\otimes \bar{r}_{1}\otimes [r_{2},\alpha(\bar{r}_{2})]\\
	&+[\beta(r_{1}),\alpha\beta(\bar{r}_{1})]\otimes r_{2}\otimes \mathrm{ad}_{x}(\bar{r}_{2})+r_{1}\otimes \mathrm{ad}_{x}(\bar{r}_{1})\otimes [\beta(r_{2}),\alpha\beta(\bar{r}_{2})]+r_{1}\otimes [\beta(r_{2}),\alpha\beta(\bar{r}_{1})]\otimes \mathrm{ad}_{x}(\bar{r}_{2})\\
	&\hspace{1cm}(\text{by Eqs.~(\ref{eq-BH-anti})-(\ref{eq-BH-Jacobi}), }(\alpha\otimes\alpha)r=r,\;(\beta\otimes\beta)r=r,\; r=-\tau(r) \text{ and }\alpha^{2}=\beta^{2}=\mathrm{id}_{L})\\
	=&\ (\mathrm{ad}_{x}\otimes \beta\otimes\beta)\Big([r_{1},\alpha\beta(\bar{r}_{1})]\otimes \beta(r_{2})\otimes \beta(\bar{r}_{2})+r_{1}\otimes [\beta(r_{2}),\alpha^{-1}\beta(\bar{r}_{1})]\otimes \beta(\bar{r}_{2})\\
	&+r_{1}\otimes \beta(\bar{r}_{1})\otimes [\beta(r_{2}),\alpha\beta(\bar{r}_{2})]\Big)+(\beta\otimes \mathrm{ad}_{x}\otimes \beta)\Big(\beta(r_{1})\otimes [r_{2},\alpha\beta(\bar{r}_{1})]\otimes \beta(\bar{r}_{2})\\
	&+[\beta(r_{1}),\alpha\beta(\bar{r}_{1})]\otimes r_{2}\otimes \beta(\bar{r}_{2})+\beta(r_{1})\otimes \bar{r}_{1}\otimes [r_{2},\alpha(\bar{r}_{2})]\Big)+(\beta\otimes\beta\otimes \mathrm{ad}_{x})\Big(\beta(r_{1})\otimes \beta(\bar{r}_{1})\\
	&\otimes [r_{2},\alpha\beta(\bar{r}_{2})]+[r_{1},\alpha(\bar{r}_{1})]\otimes \beta(r_{2})\otimes \bar{r}_{2}+\beta(r_{1})\otimes [r_{2},\alpha(\bar{r}_{1})]\otimes \bar{r}_{2}\Big)\hspace{1cm}(\text{by }\alpha^{2}=\beta^{2}=\id_{L})\\
	=&\ (\mathrm{ad}_{x}\otimes \beta\otimes \beta+\beta\otimes \mathrm{ad}_{x}\otimes\beta+\beta\otimes\beta\otimes \mathrm{ad}_{x})\Big([r_{12},r_{13}]_{\alpha\beta}+[r_{12},r_{23}]_{\alpha\beta}+[r_{13},r_{23}]_{\alpha\beta}\Big).\\
	&\hspace{1cm}(\text{by $(\beta\otimes\beta)r=r$ and $\alpha^{2}=\beta^{2}=\mathrm{id}_{L}$})
\end{align*}

Hence, the Eq.~(\ref{eq-cojacobi}) holds for $(L,\Delta_{r},\alpha,\beta)$ if and only if Eq.~(\ref{eq-coJacobi-YBE}) holds. This completes the proof.
\end{proof}

\begin{defn}\label{def-BHYBE}
Let $(L,[-,-], \alpha,\beta)$ be a BiHom-Lie algebra. We call $r\in L\otimes L$ a solution of the {\bf classical BiHom Yang-Baxter equation~(abbr.~CBHYBE)} if
\[\operatorname{CBHYB}(r):=[r_{12},r_{13}]_{\alpha\beta}+[r_{12},r_{23}]_{\alpha\beta}+[r_{13},r_{23}]_{\alpha\beta}=0\]
 in $(L,[-,-],\alpha, \beta)$. 
\end{defn}

\begin{remark}\label{prop:relation-current-BiHomYBE-CYBE}
	\begin{enumerate}
		\item When $\alpha=\beta$, the CBHYBE reduces to a classical Hom Yang-Baxter equation~\cite{Hom-Liebialg}.
		\item Furthermore, when $\alpha=\beta=\mathrm{id}_L$, the CBHYBE reduces to the classical Yang-Baxter equation in Definition~\ref{CYBE-bi}.
	\end{enumerate}
\end{remark}

Next, we introduce the relation between CBHYBEs and BiHom-Lie bialgebras.

\begin{theorem}\label{thm-BHYBE-BHLbi}
	Let $(L,[-,-], \alpha,\beta)$ be a BiHom-Lie algebra. If $r$ is a solution of the $\operatorname{CBHYBE}$ and satisfies
		\[r=-\tau(r),\quad (\alpha\otimes\alpha)r=r, \quad (\beta\otimes\beta)r=r,\]
	then $(L,[-,-],\Delta_{r},\alpha,\beta)$ is a BiHom-Lie bialgebra, where $\Delta_{r}$ is defined by Eq.~(\ref{eq-coboundary}).
\end{theorem}
\begin{proof}
	It follows from Definition~\ref{def-BHLbialg} and Lemma~\ref{lemm-equiv}.
\end{proof}

The remainder of this section is devoted to investigating the relationship between the Nijenhuis BiHom-version of the Yang-Baxter equations and Nijenhuis BiHom-Lie bialgebras.

\begin{lemma}\label{lemma_1}
	Let $S$ be an adjoint-admissible linear map to Nijenhuis BiHom-Lie algebra $(L,[-,-],\\N,\alpha,\beta)$. If $(L,[-,-],\Delta,\alpha,\beta)$ is a BiHom-Lie bialgebra, then by the Definition~\ref{def-NBHL bialg}, we have the following results:
	\begin{enumerate}
		\item\label{it-a} the 5-tuple $(L,\Delta,S,\alpha,\beta)$ is a Nijenhuis BiHom-Lie coalgebra if and only if Eq.~(\ref{NBL-coalg}) holds;
		\item\label{it-b} the linear map $N^*$ is adjoint-admissible to $(L^*,\Delta^*,S^*,\alpha^*,\beta^*)$ if and only if
		\begin{align}\label{eq30'}
			(S\otimes \mathrm{id})\Delta\circ N(x)+(\mathrm{id}\otimes N^{2})\Delta(x)=(S\otimes N)\Delta(x)+(\mathrm{id}\otimes N)\Delta\circ N(x),
		\end{align}
		for any $x\in L$.
	\end{enumerate}
\end{lemma}
\begin{proof}
	First, by Definition~\ref{def-BHLbialg} and Definition~\ref{def-NijBHLcoalg},  \ref{it-a} is obviously true. Next, we consider \ref{it-b}. By Definition~\ref{def-admissible}, $N^*$ is adjoint-admissible to $(L^*,\Delta^*,S^*,\alpha^*,\beta^*)$ if and only if
	\begin{small}
	\begin{align}\label{eq-star}
		N^*\circ\Delta^*(S^*\otimes \mathrm{id})(a\otimes b)+\Delta^*(\mathrm{id}\otimes (N^*)^{2})(a\otimes b)=\Delta^*(S^*\otimes N^*)(a\otimes b)+N^*\circ\Delta^*(\mathrm{id}\otimes N^*)(a\otimes b),
	\end{align}
	\end{small}
	for any $a,b\in L^*$.  Owing to for any $x\in L$,
	\begin{align*}
		&\ \langle  N^*\circ \Delta^*(S^*\otimes \mathrm{id})(a\otimes b)+\Delta^*(\mathrm{id}\otimes (N^*)^{2})(a\otimes b)-\Delta^*(S^*\otimes N^*)(a\otimes b)\\
		&\ -N^*\circ \Delta^*(\mathrm{id}\otimes N^*)(a\otimes b), x\rangle  \\
		=&\ \langle  a\otimes b,(S\otimes \mathrm{id})\Delta(N(x))+(\mathrm{id}\otimes N^{2})\Delta(x)-(S\otimes N)\Delta(x)-(\mathrm{id}\otimes N)\Delta(N(x))\rangle.
	\end{align*}
	
	Thus, Eq.~(\ref{eq-star}) holds if and only if Eq.~(\ref{eq30'}) holds. Therefore, $N^*$ is adjoint-admissible to $(L^*,\Delta^*, S^*, \alpha^*,\beta^*)$ if and only if Eq.~(\ref{eq30'}) holds. This completes the proof.
\end{proof}

From the above analysis, we obtain the following remark.

\begin{remark}\label{remk_2}
	If $(L,[-,-],\Delta,\alpha,\beta)$ is a BiHom-Lie bialgebra and linear map $S$ is adjoint-admissible to Nijenhuis BiHom-Lie algebra $(L,[-,-],N,\alpha,\beta)$, then $(L,[-,-],N,\Delta,S,\alpha,\beta)$ is a Nijenhuis BiHom-Lie bialgebra if and only if Eq.~(\ref{NBL-coalg}) and Eq.~(\ref{eq30'}) hold.
\end{remark}

With $\Delta_{r}$ given by Eq.~(\ref{eq-coboundary}), we further process Eq.~(\ref{NBL-coalg}) and Eq.~(\ref{eq30'}) as follows.

\begin{theorem}\label{thm_3}
	Let $S$ be a linear map that is adjoint-admissible to Nijenhuis BiHom-Lie algebra $(L,[-,-],N,\alpha,\beta)$ and define $\Delta_{r}: L\rightarrow L\otimes L$ by Eq.~(\ref{eq-coboundary}). Then
	\begin{enumerate}
		\item Eq.~(\ref{NBL-coalg}) holds if and only if for any $x\in L$,
		\begin{align*}
			&\ \big(S \mathrm{ad}_{x}\otimes \beta-\mathrm{ad}_{S(x)}\otimes \beta\big)(\mathrm{id}\otimes S-N\otimes \mathrm{id})(r)\\
			&\ +\big(\beta\otimes S \mathrm{ad}_{\alpha^{-1}\beta(x)}-\beta\otimes \mathrm{ad}_{S \alpha^{-1}\beta(x)}\big)(S\otimes \mathrm{id}-\mathrm{id}\otimes N)(r)=0.
		\end{align*}
		\item Eq.~(\ref{eq30'}) holds if and only if for any $x\in L$,
		\begin{align*}
			&\ \big(\mathrm{ad}_{N(x)}\otimes \beta+S \mathrm{ad}_{x}\otimes \beta+\beta\otimes \mathrm{ad}_{N\alpha^{-1}\beta(x)}-\beta\otimes N \mathrm{ad}_{\alpha^{-1}\beta(x)}\big)(S\otimes \mathrm{id}-\mathrm{id}\otimes N)(r)\\
			&\ -(\mathrm{ad}_{x}\otimes \beta)(S^{2}\otimes \mathrm{id}-\mathrm{id}\otimes N^{2})(r)=0.
		\end{align*}
	\end{enumerate}
\end{theorem}

\begin{proof}
	\begin{enumerate}
	\item First, we write the Eq.~(\ref{NBL-coalg}) as
	\[(S\otimes S)\Delta_{r}+\Delta_{r}\circ S^{2}-(S\otimes \mathrm{id}+\mathrm{id}\otimes S)\Delta_{r}\circ S=0.\]
	
\noindent Then, for any $x\in L$, we have
		\begin{align*}
			&\ (S\otimes S)\Delta_{r}(x)+\Delta_{r}\circ S^{2}(x)\\
			=&\  (S\otimes S)\big(\mathrm{ad}_{x}\otimes \beta+\beta\otimes \mathrm{ad}_{\alpha^{-1}\beta(x)}\big)(r)+\big(\mathrm{ad}_{S^{2}(x)}\otimes \beta+\beta\otimes \mathrm{ad}_{\alpha^{-1}\beta S^{2}(x)}\big)(r)\\
			&\hspace{1cm}(\text{by Eq.~(\ref{eq-coboundary})})\\
			=&\  \big(S \mathrm{ad}_{x}\otimes S\circ\beta+S\circ\beta\otimes S \mathrm{ad}_{\alpha^{-1}\beta (x)}+\mathrm{ad}_{ S^{2}(x)}\otimes \beta+\beta\otimes \mathrm{ad}_{\alpha^{-1}\beta S^{2}(x)}\big)(r).
		\end{align*}
\noindent and		
		\begin{align*}
			&\ (S\otimes \mathrm{id}+\mathrm{id}\otimes S)\Delta_{r}\circ S(x)\\
			=&\ (S\otimes \mathrm{id})(\mathrm{ad}_{ S(x)}\otimes \beta+\beta\otimes \mathrm{ad}_{\alpha^{-1}\beta S(x)})(r)+(\mathrm{id}\otimes S)(\mathrm{ad}_{ S(x)}\otimes \beta+\beta\otimes \mathrm{ad}_{\alpha^{-1}\beta S(x)})(r)\\
			=&\ S([ S(x),r_{1}])\otimes\beta(r_{2})+S\beta(r_{1})\otimes [\alpha^{-1}\beta S(x),r_{2}]+[ S(x),r_{1}]\otimes S\beta(r_{2})\\
			&+\beta(r_{1})\otimes S([\alpha^{-1}\beta S(x),r_{2}])\\
			=&\ -S([\alpha\beta(r_{1}),S\alpha\beta(x)])\otimes \beta(r_{2})+S\beta(r_{1})\otimes [S\alpha^{-1}\beta(x),r_{2}]+[ S(x),r_{1}]\otimes S\beta(r_{2})\\
			&-\beta(r_{1})\otimes S([\alpha\beta(r_{2}),S(x)])
			 \hspace{1cm}(\text{by Eq.~(\ref{eq-BH-anti}) and }\alpha^2=\beta^2=\mathrm{id}_{L})\\
			 =&\ \big([N\alpha\beta(r_{1}),S\alpha\beta(x)]-S([N\alpha\beta(r_{1}),\alpha\beta(x)])-[\alpha\beta(r_{1}),S^{2}\alpha\beta(x)]\big)\otimes \beta(r_{2})\\
			&+\beta(r_{1})\otimes \big([N\alpha\beta(r_{2}),S(x)]-S([N\alpha\beta(r_{2}),x])-[\alpha\beta(r_{2}),S^{2}(x)]\big)\\
			&+[S(x),r_{1}]\otimes S\beta(r_{2})+S\beta(r_{1})\otimes [S\alpha^{-1}\beta(x),r_{2}]\hspace{1cm}(\text{by Eq.~(\ref{eq-adjoint-admissible})})\\
			=&\ -[S(x),N(r_{1})]\otimes \beta(r_{2})+S([x,N(r_{1})])\otimes \beta(r_{2})+[ S^{2}(x),r_{1}]\otimes \beta(r_{2})\\
			&-\beta(r_{1})\otimes[S\alpha^{-1}\beta(x),N(r_{2})]+\beta(r_{1})\otimes S([\alpha\beta(x),N(r_{2})])+\beta(r_{1})\otimes [\alpha\beta S^{2}(x),r_{2}]\\
			&+S\beta(r_{1})\otimes [\alpha^{-1}\beta S(x),r_{2}]+[S(x),r_{1}]\otimes S\beta(r_{2}).\\
			&\hspace{1cm}(\text{by Eq.~(\ref{eq-BH-anti}) and }\alpha^{2}=\beta^{2}=\mathrm{id}_{L})
		\end{align*}

\noindent		Thus, Eq.~(\ref{NBL-coalg}) holds if and only if
		\begin{align}\label{eq-a-1}
			&\ S \mathrm{ad}_{x}(r_{1})\otimes S\beta(r_{2})+S\beta(r_{1})\otimes S \mathrm{ad}_{\alpha^{-1}\beta(x)}(r_{2})+\mathrm{ad}_{S(x)}N(r_{1})\otimes \beta(r_{2})\nonumber\\
			&-S \mathrm{ad}_{x}N(r_{1})\otimes \beta(r_{2})-S\beta(r_{1})\otimes \mathrm{ad}_{S\alpha^{-1}\beta(x)}(r_{2})-\mathrm{ad}_{S(x)}(r_{1})\otimes S\beta(r_{2})\nonumber\\
			&+\beta(r_{1})\otimes \mathrm{ad}_{\alpha^{-1}\beta S(x)}N(r_{2})-\beta(r_{1})\otimes S \mathrm{ad}_{\alpha^{-1}\beta(x)}N(r_{2})=0.
		\end{align}
		
	\noindent	Here, we compute the left-hand side (LHS) of Eq.~(\ref{eq-a-1}) as follows:
		\begin{align*}
			\text{LHS}=&\ (S \mathrm{ad}_{x}\otimes \beta)(\mathrm{id}\otimes S)(r)+(\beta\otimes S \mathrm{ad}_{\alpha^{-1}\beta(x)})(S\otimes \mathrm{id})(r)\\
			&+(\mathrm{ad}_{S(x)}\otimes \beta)(N\otimes \mathrm{id})(r)-(S \mathrm{ad}_{x}\otimes \beta)(N\otimes \mathrm{id})(r)\\
			&-(\beta\otimes \mathrm{ad}_{S\alpha^{-1}\beta(x)})(S\otimes \mathrm{id})(r)-(\mathrm{ad}_{S(x)}\otimes \beta)(\mathrm{id}\otimes S)(r)\\
			&+(\beta\otimes \mathrm{ad}_{\alpha\beta S(x)})(\mathrm{id}\otimes N)(r)-(\beta\otimes S \mathrm{ad}_{\alpha\beta(x)})(\mathrm{id}\otimes N)(r)\hspace{1cm}(\text{by }\alpha^{2}=\id_{L})\\
			=&\ (S \mathrm{ad}_{x}\otimes \beta-\mathrm{ad}_{ S(x)}\otimes \beta)(\mathrm{id}\otimes S)(r)\\
			&+(\beta\otimes S \mathrm{ad}_{\alpha^{-1}\beta(x)}-\beta\otimes \mathrm{ad}_{S\alpha^{-1}\beta(x)})(S\otimes \mathrm{id})(r)\\
			&+(\mathrm{ad}_{S(x)}\otimes \beta-S \ad_{x}\otimes \beta)(N\otimes \id)(r)\\
			&+(\beta\otimes \mathrm{ad}_{\alpha^{-1}\beta S(x)}-\beta\otimes S \mathrm{ad}_{\alpha^{-1}\beta(x)}(\mathrm{id}\otimes N)(r)\hspace{1cm}(\text{by }\alpha^{2}=\mathrm{id}_{L})\\
			=&\ (S \mathrm{ad}_{x}\otimes \beta-\mathrm{ad}_{ S(x)}\otimes \beta)(\mathrm{id}\otimes S-N\otimes \mathrm{id})(r)\\
			&+(\beta\otimes S \mathrm{ad}_{\alpha^{-1}\beta(x)}-\beta\otimes \mathrm{ad}_{S\alpha^{-1}\beta(x)})(S\otimes \mathrm{id}-\mathrm{id}\otimes N)(r).
		\end{align*}
	
\noindent	Hence, Eq.~(\ref{NBL-coalg}) holds if and only if
	\begin{align*}
		&\ \big(S \mathrm{ad}_{x}\otimes \beta-\mathrm{ad}_{S(x)}\otimes \beta\big)(\mathrm{id}\otimes S-N\otimes \mathrm{id})(r)\\
		&+\big(\beta\otimes S \mathrm{ad}_{\alpha^{-1}\beta(x)}-\beta\otimes \mathrm{ad}_{S \alpha^{-1}\beta(x)}\big)(S\otimes \mathrm{id}-\mathrm{id}\otimes N)(r)=0.
	\end{align*}

	\item
Here, for any $x\in L$, we write Eq.~(\ref{eq30'}) as
\begin{align*}
	(S\otimes \mathrm{id})\Delta_{r}\circ N(x)+(\mathrm{id}\otimes N^{2})\Delta_{r}(x)=(S\otimes N)\Delta_{r}(x)+(\mathrm{id}\otimes N)\Delta_{r}\circ N(x).
\end{align*}	

\noindent Owing to
\begin{align*}
	&\ (S\otimes \mathrm{id})\Delta_{r}\circ N(x)\\
	=&\ (S\otimes \mathrm{id})(\mathrm{ad}_{N(x)}\otimes \beta+\beta\otimes \mathrm{ad}_{\alpha^{-1}\beta N(x)})(r)
	\\
	=&\ S([N(x),r_{1}])\otimes \beta(r_{2})+S\beta(r_{1})\otimes [\alpha^{-1}\beta N(x),r_{2}]\\
	=&\ [N(x),S(r_{1})]\otimes \beta(r_{2})+S([x,S(r_{1})])\otimes \beta(r_{2})-[x,S^{2}(r_{1})]\otimes\beta(r_{2})\\
	&+S\beta(r_{1})\otimes[\alpha^{-1}\beta N(x),r_{2}]\hspace{1cm}(\text{by Eq.~(\ref{eq-adjoint-admissible})})\\
	=&\ (\mathrm{ad}_{N(x)}S\otimes\beta+S\mathrm{ad}_{x}S\otimes \beta-\mathrm{ad}_{x}S^{2}\otimes\beta+S\beta\otimes \mathrm{ad}_{\alpha^{-1}\beta N(x)})(r),
\end{align*}
\begin{align*}
&\ (\mathrm{id}\otimes N)\Delta_{r}\circ N(x)\\
=&\ (\mathrm{id}\otimes N)(\mathrm{ad}_{N(x)}\otimes \beta+\beta\otimes \mathrm{ad}_{\alpha^{-1}\beta N(x)})(r)\\
=&\ [N(x),r_{1}]\otimes N\beta(r_{2})+\beta(r_{1})\otimes N([N\alpha^{-1}\beta(x),r_{2}])\\
=&\ [N(x),r_{1}]\otimes N\beta(r_{2})+\beta(r_{1})\otimes[N\alpha^{-1}\beta(x),N(r_{2})]+\beta(r_{1})\otimes N^{2}([\alpha^{-1}\beta(x),r_{2}])\\
&-\beta(r_{1})\otimes N([\alpha^{-1}\beta(x),N(r_{2})])\hspace{1cm}(\text{by Eq.~(\ref{Nij-identity})})\\
=&\ (\mathrm{ad}_{N(x)}\otimes N\beta+\beta\otimes \mathrm{ad}_{N\alpha^{-1}\beta(x)}N+\beta\otimes N^{2} \mathrm{ad}_{\alpha^{-1}\beta(x)}-\beta\otimes N \mathrm{ad}_{\alpha^{-1}\beta(x)}N)(r),
\end{align*}
as well as
\begin{align*}
	&\ (\mathrm{id}\otimes N^{2})\Delta_{r}(x)-(S\otimes N)\Delta_{r}(x)\\
	=&\ (\mathrm{id}\otimes N^{2})(\mathrm{ad}_{x}\otimes \beta+\beta\otimes \mathrm{ad}_{\alpha^{-1}\beta(x)})(r)-(S\otimes N)(\mathrm{ad}_{x}\otimes \beta+\beta\otimes \mathrm{ad}_{\alpha^{-1}\beta(x)})(r)\\
	=&\ (\mathrm{ad}_{x}\otimes N^{2}\beta+\beta\otimes N^{2} \mathrm{ad}_{\alpha^{-1}\beta(x)}-S \mathrm{ad}_{x}\otimes N\beta-S\beta\otimes N \mathrm{ad}_{\alpha^{-1}\beta(x)})(r).
\end{align*}

\noindent Thus, Eq.~(\ref{eq30'}) holds if and only if
\begin{align}\label{eq-b-1}
	&\ (\mathrm{ad}_{N(x)}S\otimes \beta+S \mathrm{ad}_{x}S\otimes \beta+S\beta\otimes \mathrm{ad}_{\alpha^{-1}\beta N(x)}- \mathrm{ad}_{x} S^{2}\otimes \beta+\beta\otimes N \mathrm{ad}_{\alpha^{-1}\beta(x)}N)(r)\nonumber\\
	&+(\mathrm{ad}_{x}\otimes N^2\beta-S \mathrm{ad}_{x}\otimes N\beta-S\beta\otimes N\mathrm{ad}_{\alpha^{-1}\beta(x)}-\mathrm{ad}_{N(x)}\otimes N\beta-\beta\otimes \mathrm{ad}_{N\alpha^{-1}\beta(x)}N)(r)=0.
\end{align}	

\noindent Then, we deal with the left-hand side (LHS) of Eq.~(\ref{eq-b-1}) as follows.
\begin{align*}
	\text{LHS}=&\ (\mathrm{ad}_{N(x)}\otimes \beta)(S\otimes \mathrm{id})(r)+(S \mathrm{ad}_{x}\otimes \beta)(S\otimes \mathrm{id})(r)-(\mathrm{ad}_{x}\otimes \beta)(S^{2}\otimes \mathrm{id})(r)\\
	&+(\beta\otimes \mathrm{ad}_{\alpha^{-1}\beta N(x)})(S\otimes \mathrm{id})(r)+(\mathrm{ad}_{x}\otimes \beta)(\mathrm{id}\otimes N^{2})(r)-(S \mathrm{ad}_{x}\otimes \beta)(\mathrm{id}\otimes N)(r)\\
	&-(\beta\otimes N \mathrm{ad}_{\alpha^{-1}\beta(x)})(S\otimes \mathrm{id})(r)-(\mathrm{ad}_{N(x)}\otimes \beta)(\mathrm{id}\otimes N)(r)\\
	&-(\beta\otimes \mathrm{ad}_{N \alpha^{-1}\beta(x)})(\mathrm{id}\otimes N)(r)+(\beta\otimes N \mathrm{ad}_{\alpha^{-1}\beta(x)})(\mathrm{id}\otimes N)(r)\\
	=&\ (\mathrm{ad}_{N(x)}\otimes\beta+S \mathrm{ad}_{x}\otimes\beta+\beta\otimes \mathrm{ad}_{\alpha^{-1}\beta N(x)}-\beta\otimes N \mathrm{ad}_{\alpha^{-1}\beta(x)})(S\otimes \mathrm{id})(r)\\
	&+(\beta\otimes N \mathrm{ad}_{\alpha^{-1}\beta(x)}-S \mathrm{ad}_{x}\otimes \beta-\mathrm{ad}_{N(x)}\otimes \beta-\beta\otimes \mathrm{ad}_{N\alpha^{-1}\beta(x)})(\mathrm{id}\otimes N)(r)\\
	&+(\mathrm{ad}_{x}\otimes \beta)(\mathrm{id}\otimes N^{2}-S^{2}\otimes \mathrm{id})(r)\\
	=&\ (\mathrm{ad}_{N(x)}\otimes \beta+S \mathrm{ad}_{x}\otimes \beta+\beta\otimes \mathrm{ad}_{\alpha^{-1}\beta N(x)}-\beta\otimes N \mathrm{ad}_{\alpha^{-1}\beta(x)})(S\otimes \mathrm{id}-\mathrm{id}\otimes N)(r)\\
	&-(\mathrm{ad}_{x}\otimes \beta)(S^{2}\otimes \mathrm{id}-\mathrm{id}\otimes N^{2})(r).
\end{align*}
	
\noindent	Hence, Eq.~(\ref{eq30'}) holds if and only if
	\begin{align*}
		&\ \big(\mathrm{ad}_{N(x)}\otimes \beta+S \mathrm{ad}_{x}\otimes \beta+\beta\otimes \mathrm{ad}_{N\alpha^{-1}\beta(x)}-\beta\otimes N \mathrm{ad}_{\alpha^{-1}\beta(x)}\big)(S\otimes \mathrm{id}-\mathrm{id}\otimes N)(r)\\
		&-(\mathrm{ad}_{x}\otimes \beta)(S^{2}\otimes \mathrm{id}-\mathrm{id}\otimes N^{2})(r)=0.
	\end{align*}
	\end{enumerate}
	
	This completes the proof.
\end{proof}

\begin{lemma}\label{lemma_4}
	Given two linear maps $S,\,N:L\rightarrow L$, if $r=-\tau(r)\in L\otimes L$ and $S\otimes \mathrm{id}-\mathrm{id}\otimes N=0$, then
	\begin{align}
		(\mathrm{id}\otimes S-N\otimes \mathrm{id})(r)=&\ 0,\label{eq-S-id-N}\\
		\big(\mathrm{ad}_{\alpha^{-1}\beta(x)}\otimes \beta\big)(S^{2}\otimes \mathrm{id}-\mathrm{id}\otimes N^{2})=&\ 0.\nonumber
	\end{align}
\end{lemma}
\begin{proof}
	First, for Eq.~(\ref{eq-S-id-N}), we have
	\begin{align*}
		&\ (\mathrm{id}\otimes S-N\otimes \mathrm{id})(r)\\
		=&\ -(\mathrm{id}\otimes S-N\otimes \mathrm{id})\tau(r)\hspace{1cm}(\text{by }r=-\tau(r))\\
		=&\ \tau\big((\mathrm{id}\otimes N-S\otimes \mathrm{id})(r)\big)\\
		=&\ \tau(0)\\
		=&\ 0.
	\end{align*}
	
	Second, for any $x\in L$, we have
	\begin{align*}
		&\ \big(\mathrm{ad}_{\alpha^{-1}\beta(x)}\otimes \beta\big)(S^{2}\otimes \mathrm{id}-\mathrm{id}\otimes N^{2})\\
		=&\ \big(\mathrm{ad}_{\alpha^{-1}\beta(x)}\otimes \beta\big)\big((S\otimes \mathrm{id})(\mathrm{id}\otimes N)-(\mathrm{id}\otimes N)(S\otimes \mathrm{id})\big)\\
		&\hspace{1cm}(\text{by }S\otimes \mathrm{id}=\mathrm{id}\otimes N)\\
		=&\ \big(\mathrm{ad}_{\alpha^{-1}\beta(x)}\otimes \beta\big)(S\otimes N-S\otimes N)\\
		=&\ 0.
	\end{align*}
	
	This completes the proof.
\end{proof}

Combining Remark~\ref{remk_2}, Theorem~\ref{thm_3}, and Lemma~\ref{lemma_4}, we directly obtain the Nijenhuis BiHom-Lie bialgebra as follows.

\begin{theorem}\label{thm-BHLbi-NBHLbi}
	Let $(L,[-,-],\Delta_{r},\alpha,\beta)$ be a BiHom-Lie bialgebra with $\Delta_{r}$ defined by Eq.~(\ref{eq-coboundary}), and let $S$ be a  linear map that is adjoint-admissible to Nijenhuis BiHom-Lie algebra $(L,[-,-],N,\alpha,\beta)$. If
	\[r=-\tau(r),\]
	\[(S\otimes \id-\id\otimes N)(r)=0, \]
	then $(L,[-,-],N,\Delta_{r},S,\alpha,\beta)$ is a Nijenhuis BiHom-Lie bialgebra.
\end{theorem}

\begin{defn}\label{def-Nij-BiHom-YBE}
	Let $S$ be a  linear map that is adjoint-admissible to Nijenhuis BiHom-Lie algebra $(L,[-,-],N,\alpha,\beta)$. A CBHYBE satisfying
		\[(S\otimes \id-\id\otimes N)(r)=0\]
		is called an {\bf $S$-Nijenhuis classical BiHom Yang-Baxter equation~(abbr.~SN-CBHYBE)} in $(L \, , $ \\ $[-,-],N,\alpha,\beta)$.
\end{defn}

Thus, we arrive at the main result of this section: the correspondence between the Nijenhuis BiHom-Lie bialgebras and the SN-CBHYBEs.

\begin{theorem}\label{thm-NBHYBE-NBHLbi}
	Let $S$ be a  linear map that is adjoint-admissible to Nijenhuis BiHom-Lie algebra $(L,[-,-],N,\alpha,\beta)$. If $r$ is a solution of the $\operatorname{SN-CBHYBE}$ and satisfies
		\[r=-\tau(r),\quad (\alpha\otimes\alpha)r=r, \quad (\beta\otimes\beta)r=r,\]
	 then $(L,[-,-],N,\Delta_{r},S,\alpha,\beta)$ is a Nijenhuis BiHom-Lie bialgbera, where $\Delta_{r}$ is defined by Eq.~(\ref{eq-coboundary}).
\end{theorem}

\begin{proof}
	It follows from Definition~\ref{def-Nij-BiHom-YBE}, Theorem~\ref{thm-BHYBE-BHLbi} and Theorem~~\ref{thm-BHLbi-NBHLbi}.
\end{proof}

\section{Differential Lie algebras and differential Lie bialgebras}\label{sec-3}

\subsection{Differential Lie algebras and their representations}
In this subsection, we mainly study the representation of differential Lie algebras.

In the beginning, we recall the concept of differential Lie algebras.
\begin{defn}~\cite{diffLie}\label{def-diffLie}
	Let $L$ be a vector space, $d: L\rightarrow L$ and $[-,-]: L\otimes L\rightarrow L$ be two linear maps, then we call $(L,[-,-],d)$ a {\bf differential Lie algebra of weight $\lambda$} if $(L,[-,-])$ is a Lie algebra and
	\begin{align}\label{eq-diff}
		d([x,y])=[d(x),y]+[x,d(y)]+\lambda[d(x),d(y)],
	\end{align}
for any $x,y\in L,$ and $\lambda\in {\bf k}$.
\end{defn}

In the following, we give the definition of representations of differential Lie algebras.

\begin{defn}
	A {\bf representation of differential Lie algebra $(L,[-,-],d)$} is a triple $(V,\rho,\xi)$, where $(V,\rho)$ is a representation of the Lie algebra $(L,[-,-])$ and $\xi: V\rightarrow V$ is a linear map on $V$ such that
	\begin{align}\label{eq-repre-diffLie}
		\xi\big(\rho(x)v\big)=\rho\big(d(x)\big)v+\rho(x)\xi(v)+\lambda\rho\big(d(x)\big)\xi(v),
	\end{align}
for any $x\in L,\, v\in V$.
\end{defn}

\begin{exam}
	Any differential Lie algebra $(L,[-,-],d)$ is a representation on itself through its multiplications. In this case, we call $(L,[-,-],d)$ the adjoint representation.
\end{exam}

Next, we characterize representations of differential Lie algebras in terms of semi-direct product.

\begin{prop}
	Let $(L,[-,-],d)$ be a differential Lie algebra of weight $\lambda$ and $V$ be a vector space. If $\rho:L\rightarrow \mathfrak{gl}(V)$ and $\xi: V\rightarrow V$ are two linear maps, then $(V,\rho,\xi)$ is a representation of $(L,[-,-],d)$ if and only if $(L\oplus V,[-,-]_{L\oplus V},d+\xi)$ is a differential Lie algebra of weight $\lambda$, where
	\[[x+a,y+b]_{L\oplus V}:=[x,y]+\rho(x)b-\rho(y)a,\quad \text{ for any }x,y\in L,\; a,b\in V.\]
	In this case, we say that $(L\oplus V,[-,-]_{L\oplus V},d+\xi)$ is a semi-direct product of $(L,[-,-],d)$ by $(V,\rho,\xi)$.
\end{prop}

\begin{proof}
	This proof is similar to Proposition~\ref{prop-semi-direct}.
\end{proof}

In the following, we study the dual representations of differential Lie algebras.

\begin{prop}\label{prop-diff-dualrepre}
	Let $(L,[-,-],d)$ be a differential Lie algebra of weight $\lambda$. Then $(V^*,\rho^*,\zeta^*)$ is a representation of differential Lie algebra $(L,[-,-],d)$ if and only if $(V,\rho)$ is a representation of Lie algebra $(L,[-,-])$ and
	\begin{align}\label{eq-admiss-diff}
		\rho(x)\zeta(b)=\rho(d(x))b+\zeta\big(\rho(x)b\big)+\lambda\zeta\rho\big(d(x)\big)b
	\end{align}
for any $x\in L,\; b\in V$.
\end{prop}

\begin{proof}
By Corollary~\ref{coro-dual-repre-Lie}, we only need to consider the relation between Eq.~(\ref{eq-repre-diffLie}) and Eq.~(\ref{eq-admiss-diff}). For any $x\in L,\; a^*\in V^*,\; b\in V$, we have
	\begin{align*}
		&\langle  \zeta^*\big(\rho^*(x)a^*\big)-\rho^*\big(d(x)\big)a^*-\rho^*(x)\zeta^*(a^*)-\lambda\rho^*\big(d(x)\big)\zeta^*(a^*),b\rangle  \\
		=&\ \langle  \rho^*(x)a^*,\zeta(b)\rangle  +\langle  a^*,\rho\big(d(x)\big)b\rangle  +\langle  \zeta^*(a^*),\rho(x)b\rangle  +\lambda\langle  \zeta^*(a^*),\rho\big(d(x)\big)b\rangle  \\
		=&\ -\langle  a^*,\rho(x)\zeta(b)-\rho\big(d(x)\big)b-\zeta\big(\rho(x)b\big)-\lambda\zeta\rho\big(d(x)\big)b\rangle.
	\end{align*}
	Thus, Eq.~(\ref{eq-repre-diffLie}) holds for $(V^*,\rho^*,\zeta^*)$ if and only if Eq.~(\ref{eq-admiss-diff}) holds. This completes the proof.
\end{proof}

\begin{coro}
	Let $(L,[-,-],d)$ be a differential Lie algebra of weight $\lambda$ and $\pi: L\rightarrow L$ be a linear map. Then $(L^*,\mathrm{ad}^*,\pi^*)$ is a representation of $(L,[-,-],d)$ if and only if
	\begin{align}\label{eq-adjoint-diff}
		[x,\pi(y)]=[d(x),y]+\pi\big([x,y]\big)+\lambda\pi\big([d(x),y]\big),
	\end{align}
for any $x,y\in L$.
\end{coro}
\begin{proof}
	It follows from Proposition~\ref{prop-diff-dualrepre} by taking the adjoint representation.
\end{proof}

Now we give the concept of (adjoint)-admissible to differential Lie algebras.

\begin{defn}
	Let $(V,\rho)$ be a representation of Lie algebra $(L,[-,-])$.
	\begin{enumerate}
		\item A linear map $\zeta:V\rightarrow V$ is {\bf admissible to the differential Lie algebra $(L,[-,-],d)$ on $(V,\rho)$} if $(V^*,\rho^*,\zeta^*)$ is a representation of $(L,[-,-],d)$, i.e., Eq.~(\ref{eq-admiss-diff}) holds.
	\item A linear map $\pi: L\rightarrow L$ is {\bf adjoint-admissible to differential Lie algebra $(L,[-,-],d)$} if $(L^*,\mathrm{ad}^*,\pi^*)$ is a representation of $(L,[-,-],d)$, i.e., Eq.~(\ref{eq-adjoint-diff}) holds.
\end{enumerate}
\end{defn}

By~\cite{double-B}, we know that an invariant bilinear form $\mathcal{B}$ on Lie algebra $(L,[-,-])$ is a bilinear form $\mathcal{B}$ satisfying
\[\mathcal{B}\big([x,y],z\big)=\mathcal{B}\big(x,[y,z]\big),\]
for any $x,y,z\in L$.

Furthermore, linear map $\hat{d}: L\rightarrow L$ is an {\bf adjoint linear map of $d$ with respect to $\mathcal{B}$} if $\hat{d}$ satisfies Eq.~(\ref{eq-adjoint-map}).

\begin{prop}
	Let $(L,[-,-],d)$ be a differential Lie algebra of weight $\lambda$ and $\mathcal{B}$ be a nondegenerate invariant bilinear form on the Lie algebra $(L,[-,-])$. If $\hat{d}$ is an adjoint linear map of $d$ with respect to $\mathcal{B}$, then $\hat{d}$ is adjoint-admissible to $(L,[-,-],d)$.
\end{prop}

\begin{proof}
	The proof is similar to~\cite[Proposition~2.7]{RBLbialg}.
\end{proof}

\subsection{Manin triples and matched pairs of differential Lie algebras}\label{sec-4.2}
In this subsection, we mainly investigate the Manin triples of differential Lie algebras, differential Lie bialgebras, and matched pairs of differential Lie algebras.

\begin{defn}
	A {\bf Manin triple of differential Lie algebras} with respect to the symmetric invariant bilinear form $\mathcal{B}_{m}$ is a triple $\big((L\oplus L^*,[-,-]_{\oplus},d+D^*),(L,[-,-],d),(L^*,[-,-]_{*},D^*)\big)$, where $(L,[-,-],d)$ and $(L^*,[-,-]_{*},D^*)$ are differential Lie subalgebras of differential Lie algebra $(L\oplus L^*,[-,-]_{\oplus},d+D^*)$ and the natural nondegenerate symmetric bilinear form $\mathcal{B}_{m}$ on $L\oplus L^*$ is defined by Eq.~(\ref{eq-B-d}).
\end{defn}

For convenience, we write $\big((L\oplus L^*,[-,-]_{\oplus},d+D^*),(L,[-,-],d),(L^*,[-,-]_{*},D^*)\big)$ as $\big((L\oplus L^*,[-,-]_{\oplus},d+D^*),L,L^*\big)$.

Next, we introduce the dual concept of differential Lie algebras, that is differential Lie coalgebras.

\begin{defn}\label{def-diff-co}
	A {\bf differential Lie coalgebra of weight $\lambda$} is a triple $(C,\delta,d)$, where $(C,\delta)$ is a Lie coalgebra and $d: L\rightarrow L$ is a linear map satisfying
	\begin{align}\label{eq-diff-co}
		\delta d(x)=(d\otimes \mathrm{id})\delta(x)+(\mathrm{id}\otimes d)\delta(x)+\lambda(d\otimes d)\delta(x),
	\end{align}
for any $x\in L$.
\end{defn}

Now, we give the relation between differential Lie algebras and differential Lie coalgebras by the following proposition.

\begin{prop}\label{prop-diffLiealg-coalg}
	Let linear map $\delta:L \rightarrow L\otimes L$ be the linear dual of $[-,-]_{*}: L^*\otimes L^*\rightarrow L^*$. Then $(L^*,[-,-]_{*},d^*)$ is a differential Lie algebra of weight $\lambda$ if and only if $(L,\delta, d)$ is a differential Lie coalgebra of weight $\lambda$.
\end{prop}

\begin{proof}
	For any $x\in L,\;a^*,\,b^*\in L^*$, we have
	\begin{align*}
		&\ \langle  x,d^{*}[a^*,b^*]_{*}-[d^*(a^*),b^*]_{*}-[a^*,d^*(b^*)]_{*}-\lambda[d^*(a^*),d^*(b^*)]_{*}\rangle  \\
		=&\ \langle  d(x),[a^*,b^*]_{*}\rangle  -\langle  \delta(x),d^*(a^*)\otimes b^*+a^*\otimes d^*(b^*)+\lambda d^*(a^*)\otimes d^*(b^*)\rangle  \\
		=&\ \langle \delta d(x)-(d\otimes \mathrm{id}+\mathrm{id} \otimes d+\lambda d\otimes d)\delta(x), a^*\otimes b^*\rangle  \\
		=&\ \langle  \delta d(x)-(d\otimes \mathrm{id})\delta(x)-(\mathrm{id}\otimes d)\delta(x)-\lambda(d\otimes d)\delta(x), a^*\otimes b^*\rangle  .
	\end{align*}

Hence the triple $(L^*,[-,-]_{*},d^*)$ satisfies Eq.~(\ref{eq-diff}) if and only if $(L,\delta,d)$ satisfies Eq.~(\ref{eq-diff-co}). By Definition~\ref{def-diffLie} and Definition~\ref{def-diff-co}, we complete this proof.
\end{proof}

Let $d: L\rightarrow L$ be a linear map. We now study the condition under which $d^*$ is adjoint-admissible to the differential Lie algebra $(L^*,[-,-]_{*},D^*)$. By Proposition~\ref{prop-diffLiealg-coalg}, for any $x\in L,\; a^*,b^*\in L^*$, we have
\begin{align*}
	&\ \langle  [a^*,d^*(b^*)]_{*}-[D^*(a^*),b^*]_{*}-d^*([a^*,b^*])-\lambda d^*([D^*(a^*),b^*]_{*}),x\rangle  \\
	=&\ \langle  a^*\otimes d^*(b^*),\delta(x)\rangle  -\langle  D^*(a^*)\otimes b^*,\delta(x)\rangle  -\langle  [a^*,b^*]_{*},d(x)\rangle  -\lambda\langle  [D^*(a^*),b^*]_{*},d(x)\rangle  \\
	=&\ \langle  a^*\otimes b^*,(\mathrm{id}\otimes d)\delta(x)-(D\otimes \mathrm{id})\delta(x)-\delta d(x)-\lambda(D\otimes \mathrm{id})\delta\big(d(x)\big)\rangle.
\end{align*}

Hence, linear map $d^*$ is adjoint-admissible to differential Lie algebra $(L^*,[-,-]_{*},D^*)$ if and only if
\begin{align*}
	\delta d(x)+(D\otimes \mathrm{id}-\mathrm{id}\otimes d)\delta(x)+\lambda(D\otimes \mathrm{id})\delta d(x)=0,
\end{align*}
for any $x\in L$, where $\delta:L \rightarrow L\otimes L$ be the linear dual of $[-,-]_{*}: L^*\otimes L^*\rightarrow L^*$.

\begin{lemma}
	Let $\big((L\oplus L^*,[-,-]_{\oplus},d+D^*),L,L^*\big)$ be a Manin triple of differential Lie algebras. Then
	\begin{enumerate}
		\item the adjoint linear map $\widehat{d+D^*}$ of linear map $d+D^*$ with respect to $\mathcal{B}_{m}$ is $D+d^*$. Furthermore, $D+d^*$ is adjoint-admissible to differential Lie algebra $(L\oplus L^*,[-,-]_{\oplus},d+D^*)$;
		\item the linear map $D$ is adjoint-admissible to differential Lie algebra $(L,[-,-],d)$;
		\item the linear map $d^*$ is adjoint-admissible to differential Lie algebra $(L^*,[-,-]_{*},D^*)$.
	\end{enumerate}
\end{lemma}
\begin{proof}
	This proof is similar to Lemma~\ref{lemm-**}.
\end{proof}

Now we give the concept of differential Lie bialgebras as follows.

\begin{defn}\label{def-diff-Liebialg}
	A {\bf differential Lie bialgebra} is a $5$-tuple $(L,[-,-],d,\delta,D)$, where $L$ is a vector space, $[-,-]:L\otimes L\rightarrow L,\; d,D:L\rightarrow L$ are linear maps satisfying
	\begin{enumerate}
		\item the triple $(L,[-,-],\delta)$ is a Lie bialgebra;
		\item the triple $(L,[-,-],d)$ is a differential Lie algebra of weight $\lambda$;
		\item the triple $(L,\delta,D)$ is a differential Lie coalgebra of weight $\lambda$;
		\item linear map $D$ is adjoint-admissible to $(L,[-,-],d)$;
		\item linear map $d^*$ is adjoint-admissible to $(L^*,\delta^*,D^*)$.
	\end{enumerate}
\end{defn}

Next, for any $x,y\in L,\; a^*,b^*\in L^*$, we define the operation $[-,-]_{\oplus}^{\circ}: L\oplus L^{*}\rightarrow L\oplus L^{*}$ by
\[[x+a^*,y+b^*]_{\oplus}^{\circ}:=[x,y]+\mathfrak{ad}^*(a^*)y-\mathfrak{ad}^*(b^*)x+[a^*,b^*]_{*}+\mathrm{ad}^*(x)b^*-\mathrm{ad}^*(y)a^*.\]

\begin{theorem}
	Let $(L,[-,-],d)$ and $(V,[-,-]_{V},\xi)$ be two differential Lie algebras of weight $\lambda$. If $\rho: L\rightarrow \mathfrak{gl}(V)$ and $h: V\rightarrow \mathfrak{gl}(V)$ are two linear maps such that $(V,\rho,\xi)$ is a representation of $L$ and $(L,h,d)$ is a representation of $V$ satisfying
	\begin{align}
		[y,h(c)x]-[x,h(c)y]-h(\rho(y)c)(x)+h(\rho(x)c)(y)+h(c)([x,y])=0,\label{*-1}\\
		[b,\rho(z)a]_{V}-[a,\rho(z)b]_{V}-\rho(h(b)z)(a)+\rho(h(a)z)(b)+\rho(z)([a,b]_{V})=0,\label{*-2}
	\end{align}
for any $x,y\in L,\, a,b\in V$, then there is a differential Lie algebra structure $(L\oplus V,[-,-]_{\oplus}^{\circ},d+\xi)$ on the direct sum $L\oplus V$ as follows:
\[(d+\xi)(x+a):=d(x)+\xi(a),\]
	\[[x+a,y+b]_{\oplus}^{\circ}=[x,y]+h(a)(y)-h(b)(x)+[a,b]_{V}+\rho(x)b-\rho(y)a.\]
We denote this differential Lie algebra by $L\bowtie_{h}^{\,\rho}V.$
\end{theorem}

\begin{proof}
	Similar to the proof of Theorem~\ref{th-mp}. We only need to prove Eq.~(\ref{eq-diff}) on $L\oplus V$. On the one hand, we have
\begin{align*}
(d+\xi)([x+a,y+b]_{\oplus}^{\circ})=&\ (d+\xi)\big([x,y]+h(a)y-h(b)x+[a,b]_{V}+\rho(x)b-\rho(y)a\big)\\
=&\ d([x,y])+d(h(a)y)-d(h(b)x)+\xi([a,b]_{V})+\xi(\rho(x)b)-\xi(\rho(y)a).
\end{align*}
On the other hand, we have
\begin{align*}
&\ [(d+\xi)(x+a),y+b]_{\oplus}^{\circ}+[x+a,(d+\xi)(y+b)]_{\oplus}^{\circ}+\lambda[(d+\xi)(x+a),(d+\xi)(y+b)]_{\oplus}^{\circ}\\
=&\ [d(x)+\xi(a),y+b]_{\oplus}^{\circ}+[x+a,d(y)+\xi(b)]_{\oplus}^{\circ}+\lambda[d(x)+\xi(a),d(y)+\xi(b)]_{\oplus}^{\circ}\\
=&\ [d(x),y]+h(\xi(a))y-h(b)d(x)+[\xi(a),b]_{V}+\rho(d(x))b-\rho(y)\xi(a)+[x,d(y)]+h(a)d(y)\\
&-h(\xi(b))x+[a,\xi(b)]_{V}+\rho(x)\xi(b)-\rho(d(y))a+\lambda\big([d(x),d(y)]+h(\xi(a))d(y)-h(\xi(b))d(x)\\
&+[\xi(a),\xi(b)]_{V}+\rho(d(x))\xi(b)-\rho(d(y))\xi(a)\big)\\
=&\ d([x,y])+d(h(a)y)-d(h(b)x)+\xi([a,b]_{V})+\xi(\rho(x)b)-\xi(\rho(y)a).\hspace{1cm}(\text{by Eqs.~(\ref{eq-diff})-(\ref{eq-repre-diffLie})})
\end{align*}
So we get
\[(d+\xi)([x+a,y+b]_{\oplus}^{\circ})=[(d+\xi)(x+a),y+b]_{\oplus}^{\circ}+[x+a,(d+\xi)(y+b)]_{\oplus}^{\circ}+\lambda[(d+\xi)(x+a),(d+\xi)(y+b)]_{\oplus}^{\circ}.\]
This completes the proof.
\end{proof}

\begin{remark}
	\begin{enumerate}
		\item It is easy to see that Eqs.~(\ref{M-3})-(\ref{M-4}) reduce to Eqs.~(\ref{*-1})-(\ref{*-2}) when $\alpha=\beta=\mathrm{id}_{L}$ and $p=q=\mathrm{id}_{V}$.
		\item The 4-tuple $(L,V,h,\rho)$ is a matched pair of Lie algebras $(L,[-,-])$ and $(V,[-,-]_{V})$ if $(V,\rho)$, $(L,h)$ are representations of Lie algebras $(L,[-,-])$ and $(V,[-,-]_{V})$, respectively and Eqs.~(\ref{*-1})-(\ref{*-2}) hold.
	\end{enumerate}
\end{remark}

Now, we give the definition of the matched pair on differential Lie algebras as follows.

\begin{defn}
	{\bf A matched pair of differential Lie algebras $(L,[-,-],d)$ and $(V,[-,-]_{V},\xi)$} is the 6-tuple $(L,V,d,\xi,h,\rho)$, where
	\begin{enumerate}
		\item $(L,V,h,\rho)$ is a matched pair of Lie algebras $(L,[-,-])$ and $(V,[-,-]_{V})$;
		\item $(V,\rho,\xi)$ is a representation of differential Lie algebra $(L,[-,-],d)$;
		\item $(L,h,d)$ is a representation of differential Lie algebra $(V,[-,-]_{V},\xi)$.
	\end{enumerate}
\end{defn}

Next, we give the relations among Manin triples of differential Lie algebras, differential Lie bialgebras and matched pairs of differential Lie algebras.

\begin{theorem}\label{th-equiva-diff-Lie}
Let $(L,[-,-],d)$ and $(L^{*},[-,-]_{*},D^{*})$ be two differential Lie algebras. Then the following conditions are equivalent.
\begin{enumerate}
\item\label{i} The triple $\big( (L\oplus L^{*},[-,-]_{\oplus}^{\circ},d+D^{*}), L, L^{*}\big)$ is a Manin triple of differential Lie algebras.
\item\label{ii} The 5-tuple $(L,[-,-],d,\delta, D)$ is a differential Lie bialgebra, where $\delta: L\rightarrow L\otimes L$ is the linear dual of $[-,-]_{*}$.
\item\label{iii} The 6-tuple $(L,L^{*},d,D^{*},\mathfrak{ad}^{*},\mathrm{ad}^{*})$ is a matched pair of differential Lie algebras.
\end{enumerate}
\end{theorem}

\begin{proof}
The proof of the equivalence of ~\ref{i} and ~\ref{ii} is similar to that of Theorem~\ref{th}, and the proof of the equivalence of ~\ref{ii} and ~\ref{iii} is similar to that of Theorem~\ref{th-23}.
\end{proof}

\subsection{Differential Lie bialgebras and DD-CYBEs}

In this subsection, we first introduce the $D$-differential classical Yang-Baxter equations~(abbr.~DD-CYBEs), then we give the relation between differential Lie bialgebras introduced in Section~\ref{sec-4.2} and DD-CYBEs.

Let $(L,[-,-],d)$ be a differential Lie algebra of weight $\lambda$, for $r=r_{1}\otimes r_{2}\in L\otimes L$, we define $\delta_{r}: L\rightarrow L\otimes L$ by Eq.~(\ref{eq-delta-r}). Now, we consider the equivalent conditions for $(L,\delta_{r},D)$ to be a differential Lie coalgebra of weight $\lambda$ and linear map $d^*$ to be adjoint-admissible to $(L^*,\delta_{r}^*,D^*)$, respectively.
\begin{prop}\label{prop-equi-diff}
	Let $D$ be adjoint-admissible to differential Lie algebra $(L,[-,-],d)$ of weight $\lambda$. Then
	\begin{enumerate}
		\item $(L,\delta_{r},D)$ is a differential Lie coalgebra of weight $\lambda$ if and only if
		\begin{align*}
			&\ (\mathrm{ad}_{x}\otimes \mathrm{id}+\lambda D \mathrm{ad}_{x}\otimes \mathrm{id})(d\otimes \mathrm{id}-\mathrm{id}\otimes D)(r)\\
			&+(\mathrm{id}\otimes \mathrm{ad}_{x}+\mathrm{id}\otimes \lambda D \mathrm{ad}_{x})(\mathrm{id}\otimes d-D\otimes \mathrm{id})(r)=0.
		\end{align*}
		\item linear map $d^*:L^*\rightarrow L^*$ is adjoint-admissible to $(L^*,\delta_{r}^*,D^*)$ if and only if
		 \begin{align*}
		 	(\mathrm{ad}_{x}\otimes \mathrm{id}+\mathrm{id}\otimes \mathrm{ad}_{x}+\mathrm{id}\otimes \lambda \mathrm{ad}_{d(x)})(\mathrm{id}\otimes d-D\otimes \mathrm{id})(r)=0.
		 \end{align*}
	\end{enumerate}
\end{prop}

\begin{proof}
	The result follows by a similar argument to that in Lemma~\ref{lemma_1}, combined with the method of Theorem~\ref{thm_3}.
\end{proof}

\begin{remark}
Given two linear map $D,\, d : L\rightarrow L$, for $r\in L\otimes L$, if $D\otimes \mathrm{id}-\mathrm{id}\otimes d=0$, then $d\otimes \mathrm{id}-\mathrm{id}\otimes D=0$.
\end{remark}

\begin{defn}\label{def-diff-YBE}
	Let $D$ be a linear map that is adjoint-admissible to differential Lie algebra $(L,[-,-], d)$. A CYBE satisfying
	\[(D\otimes \mathrm{id}-\mathrm{id}\otimes d)(r)=0\]
	is called a {\bf $D$-differential classical Yang-Baxter equation~(abbr.~DD-CYBE)} in $(L,[-,-],d)$.
	\end{defn}

\begin{theorem}\label{thm-DYBE-diffLiebi}
	Let $D$ be a linear map that is adjoint-admissible to differential Lie algebra $(L,[-,-], d)$. If $r$ is a solution of the DD-CYBE, and satisfies
	\[r=-\tau(r),\quad (\alpha\otimes\alpha)r=r, \quad (\beta\otimes\beta)r=r,\]
	then $(L,[-,-],d,\delta_{r},D)$ is a differential Lie bialgebra, where $\delta_{r}:L\rightarrow L\otimes L$ is defined by Eq.~(\ref{eq-delta-r}).
\end{theorem}
\begin{proof}
	By Definition~\ref{CYBE-bi} and Definition~\ref{def-diff-YBE}, we obtain that $(L,[-,-],\delta_{r})$ is a Lie bialgebra. Then by the Definition~\ref{def-diff-Liebialg}, we only need to prove
	\begin{enumerate}
		\item the triple $(L,\delta_{r},D)$ is a differential Lie coalgebra of weight $\lambda$;
		\item linear map $d^*$ is adjoint-admissible to $(L^*,\delta_{r}^*,D^*)$.
	\end{enumerate}
	Thus, we complete this proof by Proposition~\ref{prop-equi-diff}.
\end{proof}

\section{Relations}\label{sec-4}

In this section, we investigate the relationships among various related structures considered in this paper. As a preliminary step, we establish the connection between differential Lie algebras and Nijenhuis Lie algebras.

	\begin{prop}\label{prop-diff-Nij}
		Let $(L,[-,-],d)$ be a differential Lie algebra of weight $\lambda$. If $\lambda\neq 0$, and we define the operator
		\[N:=d:=-\frac{\mathrm{id}}{\lambda}: L\rightarrow L,\]
		then $(L,[-,-],N)$ is a Nijenhuis Lie algebra.
	\end{prop}
\begin{proof}
	We only need to prove Eq.~(\ref{Nij-identity}). 	For any $x,y\in L$, we have
	\begin{align*}
		&\ [N(x),N(y)]+N^{2}([x,y])-N([N(x),y])-N([x,N(y)])\\
		=&\ [d(x),d(y)]+d^{2}([x,y])-d([d(x),y])-d([x,d(y)])\\
		=&\ [d(x),d(y)]+d\Big([d(x),y]+[x,d(y)]+\lambda[d(x),d(y)]\Big)-d[d(x),y]-d[x,d(y)]\\
		&\hspace{1cm}(\text{by Eq.~(\ref{eq-diff})})\\
		=&\ (1+\lambda d)[d(x),d(y)]\\
		=&\ 0.
	\end{align*}
	This completes the proof.
\end{proof}

\begin{remark}
When $d=\mathrm{id}_{L}$ and $\lambda=-1$, differential Lie algebra $(L,[-,-],d)$ of weight $\lambda$ reduces to a Lie algebra $(L,[-,-])$.
\end{remark}

Now, we summarize the relationships among some concepts defined in our paper as follows.
	
 First, we present the connections between matched pairs, Lie bialgebras, Manin triples and their differential and Nijenhuis BiHom analogues.

\[\xymatrix{
	\text{matched pairs of NBLA}\ar[r]^-{\text{Theorem~\ref{th-23}}}\ar@<.5ex>[d]^{\text{reduce}} & \text{NBLBA}\ar@<.5ex>[d]^{\text{reduce}}\ar[l] \ar[r]^-{\text{Theorem~\ref{th}}}&\text{manin triples of NBLA}\ar@<.5ex>[d]^{\text{reduce}}\ar[l] \\
	\text{matched pairs of Lie algebras}\ar[r]^-{\text{Theorem~\ref{th-eq-liebialg}}}& \text{\,\,Lie bialgebras}\ar[r]^-{\text{Theorem~\ref{th-eq-liebialg}}}\ar[l]& \text{\,\,Manin triples of Lie algebras}\ar[l]\\
	\text{matched pairs of DLA}\ar[u]^-{\text{$d$=id,$\lambda$=-1}}\ar[r]^-{\text{Theorem~\ref{th-equiva-diff-Lie}}}&\text{DLBA}\ar@<.5ex>[u]^{\text{$d$=id,$\lambda$=-1}}\ar[l]\ar[r]^-{
		\text{Theorem~\ref{th-equiva-diff-Lie}}}& \text{manin triples of DLA.}\ar[l]\ar[u]^{\text{$d$=id,$\lambda$=-1}}
}\]\\
Here, we denote Nijenhuis BiHom-Lie algebras by NBLA, differential Lie algebras by DLA, Nijenhuis BiHom-Lie bialgebras by NBLBA, differential Lie bialgebras by DLBA.

Next, we describe the relations between the various Yang-Baxter equations considered in this paper and their related structures.\\

\begin{center}
	\begin{tikzcd}[row sep=1.8em, column sep=2.8em, font=\small]
		\text{CYBE}
		\arrow[r, "\text{+BiHom}"]
		\arrow[d, "\text{+differential operator}"']
		& \text{CBHYBE}
		\arrow[d, "\text{+Nijenhuis operator}"]
		\\
		\text{DD-CYBE}
		\arrow[d, "\text{Theorem~\ref{thm-DYBE-diffLiebi}}"']
		& \text{SN-CBHYBE}
		\arrow[d,"\text{Theorem~\ref{thm-NBHYBE-NBHLbi}}"]
		\\
		\begin{array}{c}
			\text{differential Lie bialgebras} \\[1pt]
			(L,[-,-],d,\delta,D)
		\end{array}
		\arrow[d, "\text{reduce}"']
		&
		\begin{array}{c}
			\text{Nijenhuis BiHom-Lie bialgebras} \\[1pt]
			(L,[-,-],\Delta,N,S,\alpha,\beta)
		\end{array}
		\arrow[d, "\text{reduce}"]
		\\
		\begin{array}{c}
			\text{differential Lie algebras} \\[1pt]
			(L, [-,-], d)
		\end{array}
		\arrow[r, "\text{Proposition~\ref{prop-diff-Nij}}"']
		&
		\begin{array}{c}
			\text{Nijenhuis Lie algebras} \\[1pt]
			(L, [-,-], N).
		\end{array}
	\end{tikzcd}
\end{center}

\smallskip
\noindent
{\bf Acknowledgments.} This work is  partially supported by Natural Science Foundation of China (No. 12101183). Y. Y. Zhang is also sponsored by Natural Science Foundation of Henan (No.2623\\00421229), the Postdoctoral Fellowship Program of CPSF under Grant Number (No.GZC2024040\\6) and the Henan Provincial Selective Research Funding Program for Returned Scholars Studying Abroad (No.HNLX202613).

\smallskip
\noindent
{\bf Data availability} All datasets underlying the conclusions of the paper are available to readers.

\smallskip
\noindent
{\bf Conflict of interest} The authors declare that there is no conflict of interest.



\begin{thebibliography}{99}
     \bibitem{extending}
     A. L. Agore and G. Militaru, Extending structures for Lie algebras, \emph{Monatsh. Math.} {\bf 174} (2014), no. 2, 169-193.

     \bibitem{bhPL}
     H. Adimi, T. Chtioui, S. Mabrouk and S. Massoud, Constructions and representation theory of BiHom-post-Lie algebras, {\em Rend. Circ. Mat. Palermo (2)} \textbf{72} (2023), 2137-2157.

     \bibitem{double-B}
     C. M. Bai, Double constructions of Frobenius algebras, Connes cocycles and their duality, {\em J. Noncommut. Geom.} {\bf 4} (2010), no. 4, 475-530.

     \bibitem{RBLbialg}
     C. M. Bai, L. Guo, G. L. Liu and T. S. Ma, Rota-Baxter Lie bialgebras, classical Yang-Baxter equations and special L-dendriform bialgebras, \emph{Algebr. Represent. Theory} {\bf 27} (2024), no. 2, 1347-1372.

     \bibitem{preLie}
     C. Bai, Left-symmetric bialgebras and an analogue of the classical Yang-Baxter equation, {\em Commun. Contemp. Math.} {\bf 10} (2008), 221-260.


     \bibitem{YBE1}
     R. J. Baxter, Partition function of the Eight-Vertex lattics model, {\em Ann. Physics}, {\bf 70} (1972), 193-228.

     \bibitem{YBE2}
     R. J. Baxter, Exactly solved models in statistical mechanics, {\em Academic Press, Inc. London}, 1982. xii+486 pp.


     \bibitem{BHL}
     Y. S. Cheng and H. G. Qi, Representations of Bihom-Lie Algebras, {\em Algebra Colloq.} {\bf 29} (2022), no.1, 125-142.

     \bibitem{liebialg}
     V. Chari, A. Pressley, A Quide to Quantum Groups, \emph{Cambridge University Press, Cambridge} (1994).

     \bibitem{double-C}
     Y. Y. Chen, H. H. Zheng and L. Y. Zhang, Double Hom-associative algebra and double Hom-Lie bialgebra, \emph{Adv. Appl. Clifford Algebr.} {\bf 30} (2020), no.1, Paper No. 8, 25 pp.

     \bibitem{Coho-Nij-Lie}
     A. Das, Cohomology Theory of Nijenhuis Lie Algebras and (Generic) Nijenhuis Lie bialgebras, arXiv:math.RA/2502.16257.

     \bibitem{D1983}
    V. G. Drinfeld, Hamiltonian structures on Lie groups, Lie bialgebras and the geometric meaning of classical Yang-Baxter equations. {\em Dokl. Akad. Nauk SSSR}, {\bf 268} (1983), no. 2, 285-287.

     \bibitem{gmmp}
	 G. Graziani, A. Makhlouf, C. Menini and F. Panaite, BiHom-associative algebras, BiHom-Lie algebras and BiHom-bialgebras, \emph{SIGMA Symmetry Integrability Geom. Methods Appl.} \textbf{11} (2015), Paper 086, 34 pp.

     \bibitem{Nijprbialg}
     L. Guo and T. S. Ma, Nijenhuis pre-Lie bialgebras, Nijenhuis Lie bialgebras and $S$-equation. arXiv:math.RA/2508.02983.

     \bibitem{diff-RBalg}
    L. Guo and W. Keigher, On differential Rota-Baxter algebras, {\em J. Pure Appl. Algebra} {\bf 212} (2008), no. 3, 522-540.


    \bibitem{homlie}
    J. T. Hartwig, D. Larsson, S. D. Silvestrov, Deformations of Lie algebras using $\sigma$-derivations, \emph{J. Algebra} {\bf 295} (2006), 314-361.	


    \bibitem{diffalg-1}
    E. R. Kolchin, Differential Algebra and Algebraic Groups, {\em Academic Press}, New York, 1973.

    \bibitem{NIj}
    Y. Kosmann-Schwarzbach, From Poisson algebras to Gerstenhaber algebras, {\em Ann. Inst. Fourier (Grenoble)} {\bf 46} (1996), no. 5, 1243-1274.

    \bibitem{NijLie-1}
    Y. Kosmann-Schwarzbach and F. Magri, Poisson-Nijenhuis structures, {\em Ann. Inst. H. Poincar\'{e} Phys. Th\'{e}or.} {\bf 53} (1990), no. 1, 35-81.

     \bibitem{BH-Nij-operator}
    J. Li, L. Y. Chen, B. Sun, BiHom-Nijienhuis operators and $T^*$-extensions of BiHom-Lie superalgebras, {\em Hacet. J. Math. Stat.} {\bf 48} (2019), no. 3, 785-799.

    \bibitem{bhpL}
	L. Liu, A. Makhlouf, C. Menini and F. Panaite, BiHom-pre-Lie algebras, BiHom-Leibniz algebras and Rota-Baxter operators on BiHom-Lie algebras, {\em Georgian Math. J.} \textbf{28} (2021), 581-594.

    \bibitem{bhtriden}	
	L. Liu, A. Makhlouf, C. Menini and F. Panaite, Rota-Baxter operators on BiHom-associative algebras and related structures, {\em Colloq. Math.} \textbf{161} (2020), 263-294.


     \bibitem{Liealg}
    S. Lie, $\mathrm{\ddot{U}}$ber Gruppen von Transformationen, (1993), pp 8-15.


    \bibitem{diffLie}
    Y. Z. Li and D. G. Wang, Lie algebras with differential operators of any weights, {\em Electron. Res. Arch.} {\bf 31} (2023), no. 3, 1195-1211.

    \bibitem{NijYBE}
    H. Y. Li and T. S. Ma, Classical Yang-Baxter equations and Nijenhuis operators for Lie algebras, arXiv:math.RA/2502.18717.



    \bibitem{NijLie-2}
    F. Magri and C. Morosi, A geometrical characterization of integrable Hamiltonian systems through the theory of Poisson-Nijenhuis
    manifolds, Quaderno S/19, Milan, 1984. Re-issued: Universit\`a di Milano Bicocca, Quaderno 3, 2008.

    \bibitem{CYBE}
    S. Majid, Foundations of quantum group theory, Cambridge University Press, Cambridge, 1995.

   \bibitem{assoBHCYBE}
   T. S. Ma, B. Li, H. H. Zheng and X. F. Zhao, Antisymmetric infinitesimal bialgebras of any weight and related associative classical Yang-Baxter equations, {\em Algebra Colloq.}  {\bf 32} (2025), no. 4, 663-684.

    \bibitem{RBbialg}
    T. S. Ma and L. L. Liu, Rota-Baxter coalgebras and Rota-Baxter bialgebras, {\em Linear Multilinear Algebra} {\bf 64} (2016), no. 5, 968-979.

    \bibitem{Nij}
    A. Nijenhuis, Jacobi-type identities for bilinear differential concomitants of certain tensor fields. I, II. {\em Nederl. Akad. Wetensch. Proc. Ser. A} {\bf 58}. {\em Indag. Math.} {\bf 17} (1955), 390-397, 398-403.


    \bibitem{Nijoperator}
    A. Nijenhuis, $X_{n-1}$-forming sets of eigenvectors, {\em Indag. Math.} {\bf 13} (1951), 200-212.

    \bibitem{diffalg-2}
    J. F. Ritt, Differential Algebra, {\em American Mathematical Society} New York, 1950.


    \bibitem{Nij-song}
    C. Song, K. Wang, Y. Y. Zhang and G. D. Zhou, Deformations of Nijenhuis Lie algebras and Nijenhuis Lie algebroids, arXiv:math.DG/2503.22157.


    \bibitem{CYBE3}
    L. A. Tahtad\v{z}jan and L. D. Faddeev, The quantum method for the inverse problem and the $XYZ$ Heisenberg model, {\em Uspekhi Mat. Nauk} {\bf 34} (1979), no. 5 (209), 13-63, 256.

    \bibitem{BHL-coalg}
    L. L. Wu, Classification of 3-dimensional Rota-Baxter algebras (coalgebras) and BiHom-Lie bialgebras, {\em M. Sc. Thesis}, Henan University, 2018.


    \bibitem{YBE3}
    C. N. Yang, Some exact results for the many-body problem in one dimension with repulsive delta-function interaction, {\em Phys. Rev. Lett.}, {\bf 19} (1967), 1312-1315.


    \bibitem{Hom-Liebialg}
    D. Yau, The classical Hom-Yang-Baxter equation and Hom-Lie bialgebras, {\em Int. Electron. J. Algebra} {\bf 17} (2015), 11-45.


    \bibitem{diff-dendri}
    Y. Y. Zhang, H. H. Zhang, T. Z. Wu and X. Gao, Weighted differential ($q$-tri)dendriform algebras, {\em J. Algebra Appl.} {\bf 25} (2026), no. 1, Paper No. 2550304.

\end{thebibliography}
\end{document}